\documentclass[preprint,10pt]{elsarticle}
\usepackage[showframe = false, headheight = 2cm, margin = 2cm, bottom = 3cm]{geometry}
\usepackage{amssymb}
\usepackage{bigints}
\usepackage{mathrsfs}
\usepackage{calc}
\usepackage{tikz}
\usepackage{xstring}
\usepackage{tkz-euclide}
\usepackage{tikz-cd}
\usepackage{graphics}
\usepackage{multicol}
\usepackage{adjustbox}
\usepackage{subcaption}
\usepackage{multirow}
\usepackage{etoolbox}
\usepackage{datetime}
\usepackage{amsmath, accents}
\usepackage{amsthm}
\usepackage{amssymb}
\usepackage{afterpage}
\usepackage{acro}
\usepackage{float}
\usepackage{relsize}
\usepackage{array}
\usepackage{textcomp}
\usepackage{comment}
\usepackage{diagbox}
\usepackage{listings}
\usepackage{tabularx}
\usepackage{ragged2e}
\usepackage{lscape}
\usepackage{enumitem}
\usepackage{eurosym}
\usepackage{footmisc}
\usepackage{todonotes}
\usepackage{mwe}
\usepackage{gensymb}
\usepackage{wasysym}
\usepackage{graphicx}
\usepackage{blindtext}
\usepackage{appendix}
\usepackage{multicol}
\usepackage{hyperref}
\usepackage{cleveref}
\usepackage{cancel}
\usepackage{bbm}
\usepackage{mathtools}
\usepackage{nomencl}
\usepackage{accents}
\usepackage{calligra}
\usepackage{miama}
\usepackage{empheq}
\usepackage{xcolor}
\usepackage{soul}
\usepackage[normalem]{ulem}
\usepackage{algorithm}
\usepackage{algpseudocode}

\definecolor{codegreen}{rgb}{0,0.6,0}
\usetikzlibrary{patterns.meta}
\journal{}

\theoremstyle{definition}

\theoremstyle{remark}
\newtheorem{remark}{Remark}[section]

\makeatletter
\DeclareRobustCommand{\pbar}{\mathord{%
  \text{$\m@th\mkern-2mu\raisebox{-1.5ex}[0pt][0pt]{$\mathchar'26$}\mkern-7mu p$}%
}}
\makeatother

\begin{document}

\newcommand{\parddx}[2]{\dfrac{\partial #1}{\partial #2}}
\newcommand{\parddxdx}[2]{\dfrac{\partial^2 #1}{\partial #2^2}}
\newcommand{\E}[2]{\mathbb{E}^{#1, #2}}
\newcommand{\tE}[2]{\widetilde{\mathbb{E}}^{#1, #2}}
\newcommand{\tETranspose}[2]{\widetilde{\mathbb{E}}^{#1, #2^{T}}}
\newcommand{\ETranspose}[2]{\mathbb{E}^{#1, #2^{T}}}
\newcommand{\MassM}[1]{\mathbb{M}^{(#1)}}
\newcommand{\MassMP}[2]{\mathbb{M}^{(#1, #2)}}
\newcommand{\MassP}[3]{\mathbb{M}^{(#1, #2)}_{#3}}
\newcommand{\MassMinv}[1]{\mathbb{M}^{(#1)^{-1}}}
\newcommand{\weakForm}[2]{\int_{\Omega} #1 \wedge \star \: #2}
\newcommand{\boldCal}[1]{\boldsymbol{\mathcal{#1}}}
\newcommand{\bilnear}[2]{\left(#1, #2\right)_{L^2(\Omega)}}
\newcommand{\bilnearM}[3]{\left(#1, #2\right)_{L^2(#3)}}
\usetikzlibrary{math}
\begin{frontmatter}



\title{Optimal solutions employing an algebraic Variational Multiscale approach \\ Part II: Application to Navier-Stokes}

\author[1,2]{Suyash Shrestha}
\ead{s.shrestha@upm.es; s.shrestha-1@tudelft.nl}
\author[2]{Marc Gerritsma}
\author[1,3]{Gonzalo Rubio}
\author[2]{Steven Hulshoff}
\author[1,3]{Esteban Ferrer}
% \author[2,3]{Joey Dekker}
% \ead{joey.dekker.98@gmail.com}
% \ead{M.I.Gerritsma@tudelft.nl}
% \ead{S.J.Hulshoff@tudelft.nl}
% \author[4]{Ido Akkerman}
% \ead{I.Akkerman@tudelft.nl}
\affiliation[1]{organization={ETSIAE-UPM-School of Aeronautics, Universidad Politécnica de Madrid},
            addressline={Plaza Cardenal Cisneros 3},
            city={Madrid},
            postcode={E-28040},
            country={Spain}}

\affiliation[2]{organization={Delft University of Technology, Faculty of Aerospace Engineering},
            addressline={Kluyverweg 1},
            city={Delft},
            postcode={2629 HS},
            country={The Netherlands}}
\affiliation[3]{organization={Center for Computational Simulation, Universidad Politécnica de Madrid, Campus de Montegancedo},
            addressline={Boadilla del Monte},
            city={Madrid},
            postcode={E-28660},
            country={Spain}}

% \affiliation[2]{organization={Delft University of Technology, Delft Institute of Applied Mathematics},
%             addressline={Mekelweg 4},
%             city={Delft},
%             postcode={2628 CD},
%             country={The Netherlands}}
% \affiliation[4]{organization={Delft University of Technology, Faculty of Mechanical Engineering},
%             addressline={Mekelweg 2},
%             city={Delft},
%             postcode={2628 CD},
%             country={The Netherlands}}

% \author{}

% \affiliation{organization={},%Department and Organization

\begin{abstract}
This work presents a non-linear extension of the high-order discretisation framework based on the Variational Multiscale (VMS) method previously introduced for steady linear problems. We build on the concept of an optimal projector defined via the symmetric part of the governing operator. Using this idea, we generalise the formulation to the two-dimensional incompressible Navier–Stokes equations. The approach maintains a clear separation between resolved and unresolved scales, with the fine-scale contribution approximated through the approximate Fine-Scale Greens' operator of the associated symmetric operator. This enables a consistent variational treatment of non-linearity while preserving high-order accuracy. We show that the method yields numerical solutions that closely approximate the optimal projection of the continuous/highly-resolved solution and inherits desirable conservation properties. Particularly, the formulation guarantees discrete conservation of mass, energy, and vorticity, where enstrophy conservation is also achieved when exact or over-integration is employed. Numerical results confirm the methodology’s robustness and accuracy, while also demonstrating its computational cost advantage compared to the baseline Galerkin approach for the same accuracy.

\end{abstract}


% \begin{highlights}
%     \item Construction and use of algebraic dual basis functions w.r.t. $L^2$ and $H_0^1$ for the projection onto the resolved scales
%     \item Explicit formulation of the Fine-Scale Greens' function in terms of $L^2$ and $H_0^1$ dual bases
% \end{highlights}

\begin{keyword}
Optimal projections; Variational Multiscale; Fine-Scale Greens' function; Spectral Element Method; Stokes flow; Navier-Stokes;
\end{keyword}

\end{frontmatter}

\tableofcontents


\section{Introduction}
This work continues the development initiated in \cite{Shrestha2025OptimalProblems}, where we proposed a discretisation methodology grounded in the Variational Multiscale (VMS) paradigm \cite{Hughes1998TheMechanics, ShakibFarzin1989FiniteEquations} for steady linear problems. In that first part, we introduced an abstract VMS formulation at the discrete level, centred around the concept of an \emph{optimal projector}, and demonstrated how to recover projected solutions by numerically approximating the Fine-Scale Greens' function \cite{Shrestha2024ConstructionScales, Hughes2007VariationalMethods} of the symmetric part of the governing differential operator. The resulting scheme yielded exponential convergence toward optimal projections and was verified on steady advection–diffusion problems in both direct and mixed formulations.

In the present work, we extend this methodology to more general unsteady and non-linear problems, specifically the incompressible Navier–Stokes equations. This extension builds on the abstract formulation from Part I \cite{Shrestha2025OptimalProblems}, incorporating non-linear residuals into the coarse-scale–fine-scale interaction. The goal remains consistent with Part I: to compute a numerical solution that corresponds to an a priori chosen projection of the infinite-dimensional solution. 
% The extension to non-linear settings necessitates a careful treatment of the coarse–fine scale interaction, especially as the unresolved scales now appear non-linearly in the governing equations. 
We demonstrate how the previously developed methodology can be systematically generalised to such cases while preserving its variational character and high-order accuracy. This paper specifically employs the Mimetic Spectral Element Method (MSEM) \cite{Jain2021ConstructionMeshes}, but the concepts are readily extendable to other Finite/Spectral Element, or Isogeometric frameworks. In particular, the extension is straightforward for discretisations that respect the discrete de Rham complex, as the underlying structural properties are naturally preserved. For methods that do not satisfy the de Rham sequence, some reformulations may be required, including the verification of stability properties such as the inf-sup condition and appropriate reconstruction of compatible function spaces.

The VMS framework, introduced in \cite{hughes_1995, Brezzi1997BG, Hughes1998TheMechanics}, interprets stabilisation through a systematic decomposition into coarse (resolved) and fine (unresolved) scales. A central development for the present work is the projection-based
reformulation of \cite{Hughes2007VariationalMethods}, which formalised Fine-Scale Greens' functions. In this formulation, scale separation follows directly from a projector, and the coarse-scale solution is precisely the chosen projection of the exact solution onto the resolved space. Following this development, VMS has had a significant impact across applications, most notably in stabilisation and Large Eddy Simulation (LES), where numerous studies \cite{munts_hulshoff_de_borst_2004, bazilevs_calo_cottrell_hughes_reali_scovazzi_2007, Holmen_Hughes_Oberai_Wells_2004, stoter_turteltaub_hulshoff_schillinger_2018, ten_eikelder_bazilevs_akkerman_2020} have shown how it enables LES formulations driven by residual-based fine-scale closures rather than heuristic eddy viscosities. These advances span incompressible and compressible flow \cite{Holmen2004SensitivityFlow, Bazilevs2007VariationalFlows, Koobus2004AShedding} and include high-order spectral element extensions \cite{Bazilevs2010LargeMethod, Wang2010SpectralMethod}. Additional advancements include orthogonal-subscale formulations \cite{codina2000stabilization, codina2002analysis, codina2008analysis}, as well as methods for rarefied gas dynamics \cite{Baidoo2024ExtensionsEquation}. A comprehensive review of the advancements is found in \cite{Rasthofer2018RecentFlow, Ahmed2017AFlows}. Beyond turbulence modelling, VMS principles have also been employed in numerical homogenisation and adaptive finite element methods. Numerical homogenisation is often realised through local fine-scale problems \cite{efendiev2013generalized, maalqvist2020numerical, altmann2021numerical}, notably in the influential Local Orthogonal Decomposition approach \cite{henning2013oversampling, hauck2023super, Henning2025StableSpaces} where these local fine-scales are used to enhance approximation qualities of the coarse scales. In adaptive strategies, VMS ideas have similarly informed a range of developments, including those in \cite{larson2009adaptive, Giraldo2023ANorms} where the fine-scales guide the mesh refinement procedure.

Alongside the development of VMS, substantial progress has been made on structure-preserving formulations of the Navier–Stokes equations, aimed to retain key invariants of the continuous system at the discrete level. Early work demonstrated conservation of momentum and kinetic energy in covolume schemes \cite{Perot2000ConservationSchemes} and preservation of circulation and kinetic energy via the rotational form \cite{Zhang2002AccuracyDynamics}. Subsequent developments introduced discrete diffusion operators \cite{Perot2007DiscreteDiffusion}, Discrete Exterior Calculus (DEC)-based conservative schemes \cite{Elcott2007StableFluids, Mullen2009Energy-preservingAnimation, Mohamed2016DiscreteMeshes}, and velocity–vorticity approaches \cite{Rebholz2007AnEquations, Olshanskii2010VelocityvorticityhelicityEquations}. A particularly influential line of work is the Mimetic Energy Enstrophy Vorticity-Conserving (MEEVC) formulation, which enforces conservation of mass, energy, enstrophy, and vorticity \cite{Palha2017AEquations, deDiego2019InclusionScheme, Zhang2024AConditions}. The most recent of said formulation \cite{Zhang2024AConditions} introduces a mixed spectral element method based on the rotational form that preserves these invariants under general boundary conditions.

Despite these advances, structure-preserving and multiscale modelling techniques have largely evolved independently. Only a few studies \cite{vanOpstal2017IsogeometricFlows, coley2017residual} have attempted to combine VMS ideas with conservative formulations. A related development appears in \cite{Budninskiy2019Operator-adaptedForms}, which introduces operator-adapted wavelets that bear strong similarities to the VMS framework. To the best of our knowledge, however, no existing method integrates a high-order structure-preserving approach, such as the MEEVC formulation, with a VMS formulation that aims to recover a target projection of the solution. This is precisely the main contribution of the present work. We formulate a VMS methodology compatible with the MEEVC scheme, enabling the computation of projected solutions that preserve key physical invariants at the discrete level. To maintain these conservation properties, we deliberately avoid localising the fine-scale problem. While localisation techniques are highly advantageous, they introduce additional approximations that would compromise the MEEVC properties guaranteed by our construction of the spaces and projector. Nonetheless, we demonstrate that the proposed methodology remains computationally efficient. This synthesis therefore provides a structure-preserving approach that yields optimal solutions, with optimality determined by the choice of projector in the VMS formulation.

We take a moment to clarify the intended meaning of the terms ``optimal solution'' and ``optimal projector'' as used in this work. The notion of optimality here is specific to the projector we construct, which is chosen based on its alignment with the structure of the governing PDE. This should not be interpreted as a claim of universal superiority over all other possible projectors or solutions. Rather, it reflects a deliberate design choice guided by the goals and formulation of the present study. While alternative projectors may be better suited in different contexts or under different norms, our aim is not to exhaustively compare all such options. Instead, we focus on demonstrating how the proposed methodology can recover solutions that are optimal with respect to the chosen projection.

This paper is structured as follows. In \Cref{sec:review}, we briefly revisit and summarise the VMS methodology from Part I. In \Cref{sec:stokes} we present a mixed formulation for Stokes flow and define the core structure of our \emph{optimal projector}. Following that, we present the projection-based VMS formulation for the 2D incompressible Navier-Stokes equations in \Cref{sec:navier_stokes}. Furthermore, we discuss the cost estimates of solving the VMS formulation along with the MEEVC properties of the overall scheme. A convergence analysis along with a numerical illustration of the proposed methodology is presented in \Cref{sec:num_test} to verify the method. Finally, we conclude in \Cref{sec:summary} with a summary of the presented work.

\section{Background and theoretical preliminaries}
\label{sec:review}

\subsection{The abstract problem}
Consider a generic elliptic problem defined as follows. Let $\Omega \subset \mathbb{R}^n$ be a contractible open bounded domain with a Lipschitz continuous boundary $\partial \Omega$. Let $V$ be a Hilbert space in $\Omega$ endowed with a norm $\lVert\cdot \rVert_V$ and an inner product $(\cdot,\cdot)_V$. Let $\mathcal{L}$ be a symmetric positive-definite operator which is taken to be a linear isomorphism between $V$ and its dual space $V^*$, $\mathcal{L}: V \to V^*$. We state the generic problem as: given ${f} \in V^*$, find $\phi \in V$ such that
\begin{align}
    \mathcal{L} \phi &= {f}, \quad \text{in } \Omega .\label{eq:pde}
    % \mathcal{B} \underline{u} &= g, \quad \text{on } \partial \Omega,
    % \label{eq:pde_bc}
\end{align}
We may variationally characterise the solution to such a problem as the critical (saddle) point of the energy functional $\mathcal{I}$ associated with $\mathcal{L}$ \cite{Arnold2010FINITESTABILITY, Boffi2013MixedApplications, Marsden2008TheMethods}. This ultimately leads to a weak formulation of the original problem. The weak formulation can be discretised using the Galerkin approach with the restriction of the test and trial spaces to a \emph{suitably} chosen set of finite-dimensional spaces \cite{Arnold2010FINITESTABILITY}. Provided that we choose these spaces suitably, we can show that the corresponding discrete solution is in fact the $\mathcal{L}$-orthogonal projection of the continuous (infinite-dimensional) solution onto the finite-dimensional space where the projector is defined using the energy functional $\mathcal{I}$. Hence, the discrete solution is optimal in the error norm embedded in the structure of $\mathcal{I}$. We demonstrated this concept in \cite{Shrestha2025OptimalProblems} for a simple Poisson problem, and we shall demonstrate it for the Stokes flow problem in this paper in \Cref{sec:stokes}. 

The problem in \eqref{eq:pde} possesses the advantageous variational structure when the problem consists solely of a symmetric positive-definite operator $\mathcal{L}$. This structure, however, is lost for problems involving additional skew-symmetric operators. To be more concrete, consider a generic operator problem
\begin{equation}
    \mathcal{L}\phi + \mathcal{C} \phi = f, \quad \text{in }\Omega,
    \label{eq:pde_skew}
\end{equation}
with $\mathcal{C}$ being a skew-symmetric, possibly non-linear operator. In this case, the solution can no longer be characterised as the critical (saddle) point of an energy functional. Consequently, a standard Galerkin discretisation of this problem does \emph{not} yield a solution that is the projection of the continuous solution. In other words, the Galerkin solution is sub-optimal with respect to the norm induced by $\mathcal{I}$, which is associated with the underlying symmetric positive-definite problem. This loss of optimality for the Galerkin scheme for generic operators does not necessarily affect the stability of the discretisation. What it does imply is that the solutions on coarse spaces may share a weaker likeness to the infinite-dimensional (or highly resolved) solutions, attributed to a larger error in the optimal norm. 

\subsection{Review of Part I}
\begin{figure}[htp]
    \centering
    \begin{tikzpicture}
    % Outer dotted circle
    \draw[dotted, thick] (0,0) circle (3cm);
    
    \filldraw[fill=gray!1, pattern={Lines[
                  distance=2mm,
                  angle=-45,
                  line width=0.05mm
                 ]}] (0,0) circle (2cm);
    
    % Inner darker shaded circle
    \filldraw[fill=gray!70] (0,0) circle (1cm);

    \node at (2.5, -2.5) {$V$};

    \node at (0, 0) {$\bar{V}$};

    \node at (1, 1) {${V}^{\prime}_k$};
\end{tikzpicture}
    \caption{Abstract VMS framework at the discrete level}
    \label{fig:VMS_frame}
\end{figure}
In Part~I~\cite{Shrestha2025OptimalProblems}, we introduced an algebraic VMS methodology tailored to \emph{linear} advection-diffusion problems. The method is based on our \emph{optimal projector} \(\mathcal{P}: V \to \bar{V}\), which maps onto the finite-element subspace \(\bar{V}\) constructed from \(\mathcal{L}\), the symmetric positive-definite part of the differential operator. The null space of $\mathcal{P}$, consisting of all components orthogonal to $\bar{V}$, forms a closed subspace $V'$, thereby establishing the scale separation $V = \bar{V} \oplus V' \Rightarrow \phi = \bar{\phi} + \phi'$. The problem in \eqref{eq:pde_skew} is then cast in a split form, noting the linearity of $\mathcal{C}$ as
\begin{equation}
    \begin{array}{ccccc}
        \bar{\mathcal{L}}\bar{\phi} + \bar{\mathcal{C}}\bar{\phi} & + & \bar{\mathcal{C}}\phi' &=& \bar{f}, \\[0.2cm]
        \mathcal{C}'\bar{\phi} & + & \mathcal{L}'\phi' + \mathcal{C}'\phi' &=& f',
    \end{array}
    \label{eq:sys_0}
\end{equation}
where, noting the $\mathcal{L}$-orthogonality between the scales allows us to eliminate some terms.

With the scale separation defined, the split form is solved using a refined-mesh approach with a hierarchical polynomial basis. The resolved (coarse) scales are approximated on $\bar{V}$, while the unresolved (fine) scales are represented on a separate $p$-refined space $V'_k$. The fine scales are constructed explicitly via an approximate Fine-Scale Greens' operator~\cite{Shrestha2024ConstructionScales} ($\mathcal{G}'$), requiring only the Greens' operator of the \emph{symmetric} part of the operator ($\mathcal{G}' = (\mathcal{L}')^{-1}$). This abstraction is illustrated in \Cref{fig:VMS_frame}. Formally, an infinite-dimensional Hilbert space $V$ is decomposed into $\bar{V}$ and its orthogonal complement $V'$, with $V'_k$ approximating $V'$ by enriching $\bar{V}$ with a polynomial increment $k$. Thus, $\bar{V}$ uses degree $p$, while $V'_k$ uses $p+k$, with the projection onto $\bar{V}$ removed.

This formulation was verified in~\cite{Shrestha2025OptimalProblems} on 1D and 2D linear advection-diffusion benchmarks, yielding coarse-scale solutions closely matching the projected infinite-dimensional solution, while also approximating the unresolved scales. The solution strategy solves the full PDE on the coarse space and only solves for the Fine-Scale Greens' function of the underlying symmetric problem on the refined space. The present work extends this algebraic VMS methodology to the non-linear Navier-Stokes equations, detailing the discrete spaces and optimal projector derived from the symmetric Stokes problem.

\section{Stokes flow}
\label{sec:stokes}

\subsection{Continuous formulation}
Consider a Stokes flow problem stated as follows
\begin{align}
    -\Delta \underline{u} + \nabla \pbar &= \underline{f}, \quad \text{in } \Omega \\
    \nabla \cdot \underline{u} &= 0, \quad \text{in } \Omega \\
    \underline{u} \cdot \hat{n} = u_{\perp}, \; \text{on } \Gamma_n, &\quad \quad \nabla \times \underline{u} = \hat{\omega}, \quad \text{on} \: \Gamma_{\omega}  \\
    \pbar = \hat{\pbar}, \quad \text{on} \: \Gamma_p, &\quad \quad {\underline{u}} \times \hat{n} = u_{\parallel}, \quad \text{on} \: \Gamma_t \\
    \text{with } \Gamma_n \cup \Gamma_p = \Gamma_{\omega} \cup \Gamma_t &= \partial \Omega, \quad \text{and} \: \Gamma_n \cap \Gamma_p = \Gamma_{\omega} \cap \Gamma_t = \varnothing,
\end{align}
where we seek a velocity field $\underline{u}$ and a pressure field $\pbar$ which acts as a Lagrange multiplier to enforce $\underline{u}$ to be divergence-free. We rewrite the problem in a mixed formulation by noting the identity of the vector Laplacian $\Delta \underline{u} = \nabla \nabla \cdot \underline{u} - \nabla \times \nabla \times \underline{u}$. We thus have the following mixed formulation
\begin{align}
    \omega - \nabla \times \underline{u} &= 0, \quad \text{in } \Omega \\
    \nabla \times \omega + \nabla \pbar &= \underline{f}, \quad \text{in } \Omega \\
    \nabla \cdot \underline{u} &= 0, \quad \text{in } \Omega,
    % \underline{u} \cdot \hat{n} = g, \quad \text{on} \: \Gamma_n \\
    % \omega \times \hat{n} = h, \quad \text{on} \: \Gamma_t \\
    % \text{with } \Gamma_n \cup \Gamma_t = \partial \Omega,
\end{align}
with $\omega$ being the vorticity. We define the three quantities $\omega$, $\underline{u}$ and $\pbar$ in separate compatible spaces. Since we limit the focus of this work to 2D, we work with the following Sobolev spaces
\begin{align*}
    L^2(\Omega) &:= \{\eta \: | \: \bilnear{\eta}{\eta} < +\infty\} \\
    H(\mathrm{curl}, \Omega) &:= \{\varepsilon \in L^2(\Omega) \: | \: \nabla \times \varepsilon \in [L^2(\Omega)]^2 \} \\
    H(\mathrm{div}, \Omega) &:= \{\underline{v} \in [L^2(\Omega)]^2 \: | \: \nabla \cdot \underline{v} \in L^2(\Omega) \},
\end{align*}
with each space equipped with an inner product and a naturally induced norm, thus making them Hilbert spaces. In 2D, the $H(\mathrm{curl}, \Omega)$ space contains all vector fields that solely point in the out-of-plane direction, whose curl is square integrable. The $H(\mathrm{div}, \Omega)$ space contains all vector fields whose vector components are square integrable and whose divergence is square integrable. We seek $\omega \in H(\mathrm{curl}, \Omega)$, $\underline{u} \in H(\mathrm{div}, \Omega)$, and $\pbar \in L^2(\Omega)$, and we take three finite-dimensional spaces, $\bar{U}$ as the subspace of $H(\mathrm{curl}, \Omega)$, $\bar{V}$ as the subspace of $H(\mathrm{div}, \Omega)$ and $\bar{W}$ as the subspace of $L^2(\Omega)$. 
\[
\bar{U} \subset  H(\mathrm{curl}, \Omega), \quad \bar{V} \subset  H(\mathrm{div}, \Omega), \quad \bar{W} \subset L^2(\Omega)
\]
The three spaces form a de~Rham complex (or Hilbert sub-complex)~\cite{Arnold2010FINITESTABILITY}, and their finite-dimensional counterparts are constructed to preserve this structure. Following the MSEM framework~\cite{jain_zhang_palha_gerritsma_2021}, the construction of these finite-dimensional spaces ensures that $\mathrm{Range}(\nabla \times) \subseteq \mathrm{Ker}(\nabla\cdot) \subset \bar{V}$ and that $\mathrm{Range}(\nabla\cdot) \subset \bar{W}$, thereby forming a de~Rham complex. Moreover, we consider contractable domains in this work, and hence have an exact complex with $\mathrm{Range}(\nabla \times) = \mathrm{Ker}(\nabla\cdot)$.
% spaces are defined as follows:
% \begin{itemize}
%     \item $\bar{U}$ consists of scalar functions constructed as a tensor product of two 1D polynomials of degree $p$
%     \item $\bar{V}$ consists of vector fields whose first component has degree $(p,\, p-1)$ in $x$ and $y$, and whose second component has degree $(p-1,\, p)$
%     \item $\bar{W}$ consists of scalar functions with polynomial degree $p-1$ in both $x$ and $y$.
% \end{itemize}
% This construction ensures that $\mathrm{Range}(\nabla \times) \subseteq \mathrm{Ker}(\nabla\cdot) \subset \bar{V}$ and that $\mathrm{Range}(\nabla\cdot) \subset \bar{W}$, thereby forming a de~Rham complex. Moreover, we consider contractable domains in this work, and hence have an exact complex with $\mathrm{Range}(\nabla \times) = \mathrm{Ker}(\nabla\cdot)$.

% \[ \mathbb{R} \hookrightarrow H(\mathrm{curl}, \Omega) \xrightarrow{\;\;\nabla \times\;\;} H(\mathrm{div}, \Omega) \xrightarrow{\;\;\nabla \cdot\;\;} L^2(\Omega)\longrightarrow 0\;. \] 
% \[ \mathbb{R} \hookrightarrow \bar{U} \xrightarrow{\;\;\nabla_h \times\;\;} \bar{V} \xrightarrow{\;\;\nabla_h \cdot\;\;} \bar{W}\longrightarrow 0\;. \] 
\[
\begin{tikzcd}
    \mathbb{R} \arrow[r, hook] & H(\mathrm{curl}, \Omega) \arrow[r,"\nabla \times"] & H(\mathrm{div}, \Omega) \arrow[r,"\nabla \cdot"]  & L^2(\Omega) \arrow[r] & 0 \\
     \mathbb{R} \arrow[r, hook] & \bar{U} \arrow[r,"\nabla \times"] & \bar{V} \arrow[r,"\nabla \cdot"] & \bar{W} \arrow[r] & 0
\end{tikzcd}
\]
\begin{remark}
    While we denote the discretisation order of this sequence by a single polynomial degree $p$, the actual spaces differ slightly: $\bar{U}$ consists of polynomials of degree $p$, $\bar{V}$ consists of polynomials of degrees $p$ and $p-1$, and $\bar{W}$ consists of polynomials of degree $p-1$.
    \label{rem:degree}
\end{remark}
The choice to define the three physical quantities in these separate spaces ensures well-posedness of the mixed formulation at the discrete level, as will be discussed later. 
% The energy functional for Stokes flow is given by
% \begin{gather}
%     J(\underline{u}) := \frac{1}{2} \left\lVert \sqrt{\nu} \: \nabla \times \underline{u} \right\rVert_{L^2}^2 = \frac{1}{2} \int_{\Omega} \nu |\omega|^2 \: \mathrm{d}\Omega.
% \end{gather}
The mixed Stokes flow problem can be variationally characterised with the following constrained energy functional
\begin{gather}
    \mathcal{I}(\omega, \underline{u}, \pbar; \underline{f}, u_{||}, \hat{\pbar}) := \int_{\Omega} \frac{1}{2} |\omega|^2 \: \mathrm{d}\Omega - \int_{\Omega} \underline{u} \cdot \nabla \times \omega \: \mathrm{d}\Omega + \int_{\Omega} \pbar \nabla \cdot \underline{u} \: \mathrm{d}\Omega + \int_{\Omega} \underline{u} \underline{f}  \: \mathrm{d}\Omega - \int_{\partial \Omega} \omega u_{||} \: \mathrm{d} \Gamma_t - \int_{\partial \Omega} \underline{u} \hat{\pbar} \: \mathrm{d} \Gamma_p,
\end{gather}
where the solution $(\omega, \underline{u}, \pbar) \in H(\mathrm{curl}, \Omega) \times H(\mathrm{div}, \Omega) \times L^2(\Omega)$ is the saddle (critical) point of the functional. Consequently, we arrive at a weak form of the problem by taking variations, which yields the following system
% \begin{gather}
%     (\omega, \underline{{\omega}})_{L^2(\Omega)} - (\nabla \times \omega,  \underline{u})_{L^2(\Omega)} = -\int_{\partial \Omega} \omega u_{||} \: \mathrm{d} \Gamma_t, \quad \forall \omega \in H(\mathrm{curl}, \Omega) \label{eq:stokes_weak_1}\\
%     -(\underline{v}, \nabla \times \omega)_{L^2(\Omega)} + (\nabla \cdot \underline{v}, p)_{L^2(\Omega)} = (\underline{v}, \underline{f})_{L^2(\Omega)}, \quad \forall \underline{v} \in H(\mathrm{div}, \Omega) \label{eq:stokes_weak_2}\\
%     (\eta, \nabla \cdot \underline{\bar{u}})_{L^2(\Omega)} = 0, \quad \forall \eta \in L^2(\Omega). \label{eq:stokes_weak_3}
% \end{gather}
\begin{equation}
    \begin{array}{cccccccl}
(\varepsilon, {{\omega}})_{L^2(\Omega)} & - & (\nabla \times \varepsilon,  \underline{u})_{L^2(\Omega)} & & & = &\displaystyle\int_{\partial \Omega} \varepsilon u_{||} \: \mathrm{d} \Gamma_t,  & \forall \varepsilon \in H(\mathrm{curl}, \Omega) \\[1ex]
-(\underline{v}, \nabla \times \omega)_{L^2(\Omega)} 
& & + & & (\nabla \cdot \underline{v}, \pbar)_{L^2(\Omega)} 
&= &\displaystyle\int_{\partial \Omega} \underline{v} \hat{\pbar} \: \mathrm{d} \Gamma_p - (\underline{v}, \underline{f})_{L^2(\Omega)} , 
& \forall \underline{v} \in H(\mathrm{div}, \Omega) \\[1ex]
& & (\eta, \nabla \cdot \underline{{u}})_{L^2(\Omega)} & & &= & 0, 
&  \forall \eta \in L^2(\Omega)
\end{array}.
\label{eq:stokes_weak}
\end{equation}
Note that this is still a continuous weak form as we are considering the solution to be in infinite-dimensional function spaces.

\subsection{Projector definition}
We are naturally interested in the discretisation of the problem where we seek $(\bar{\omega}, \underline{\bar{u}}, \bar{\pbar}) \in \bar{U} \times \bar{V} \times \bar{W}$. The Galerkin approach would be to simply restrict both the test and trial spaces to corresponding finite-dimensional spaces and arrive at a discrete system. In line with our prior work in \cite{Shrestha2025OptimalProblems}, we show that the same Galerkin algebraic system can be obtained by invoking the notion of what we define as our \emph{optimal projector} $\mathcal{P} : H(\mathrm{curl}, \Omega) \times H(\mathrm{div}, \Omega) \times L^2(\Omega) \rightarrow \bar{U} \times \bar{V} \times \bar{W}$. We define this projector as the saddle point of the functional $\mathcal{I}(\bar{\omega} - \omega, \underline{\bar{u}} - \underline{u}, \bar{\pbar} - \pbar; 0, 0, 0)$
% \begin{gather}
%     \mathcal{P} \underline{u} = \argmin_{\underline{\bar{u}} \in \bar{V}} \left\{\frac{1}{2} \left\lVert \nabla \times \underline{\bar{u}} - \nabla \times \underline{u} \right\rVert_{L^2}^2 \right\}, \quad \quad s.t \: \nabla \cdot \underline{\bar{u}} = 0.
%     \label{eq:stokes_mix_proj}
% \end{gather}
\begin{gather}
    \mathcal{I}(\bar{\omega} - \omega, \underline{\bar{u}} - \underline{u}, \bar{\pbar} - \pbar; 0, 0, 0) := \int_{\Omega}  \frac{1}{2} |\bar{\omega} - \omega|^2 \: \mathrm{d}\Omega - \int_{\Omega}(\underline{\bar{u}} - \underline{u}) \cdot \nabla \times (\bar{\omega} - \omega) \: \mathrm{d}\Omega + \int_{\Omega} (\bar{\pbar} - \pbar) \nabla \cdot (\underline{\bar{u}} - \underline{{u}}) \: \mathrm{d} \Omega,
    \label{eq:stokes_energy}
\end{gather}
where $(\omega, \underline{u}, p)$ are the solution to the continuous weak form in \eqref{eq:stokes_weak}. Taking variations with respect to $(\bar{\omega}, \bar{\underline{u}}, \bar{\pbar})$ gives the following system that defines the projector
\begin{equation}
    \begin{array}{ccccccl}
        \bilnear{\varepsilon^h}{\bar{\omega} - {\omega}} & - & \bilnear{\nabla \times \varepsilon^h}{\bar{\underline{u}} - {\underline{u}}} & & = &0,  & \forall \varepsilon^h \in \bar{U} \\[1.5ex]
        -\bilnear{\underline{v}^h}{\nabla \times (\bar{\omega} - {\omega})}
        & & + & \bilnear{\nabla \cdot \underline{v}^h}{\bar{\pbar} - {\pbar}}
        &= & \underline{0}, 
        & \forall \underline{v}^h \in \bar{V} \\[1.5ex]
        & & \bilnear{\eta^h}{\nabla \cdot (\bar{\underline{u}} - \underline{{u}})} & & = & 0, 
        &  \forall \eta^h \in \bar{W}
    \end{array}.
    \label{eq:stokes_projector}
\end{equation}

If we now consider the first equation of \eqref{eq:stokes_projector}, we have
% \begin{gather}
%     (\varepsilon^h, \underline{\tilde{\omega}})_{L^2(\Omega)} + (\nabla \times \varepsilon^h,  \underline{\tilde{u}})_{L^2(\Omega)} = 0, \quad \forall \varepsilon^h \in H(\mathrm{curl}, \Omega).
% \end{gather}
% Take $\underline{\tilde{\omega}} = \bar{\omega} - \omega$ and $\underline{\tilde{u}} = \underline{\bar{u}} - \underline{u}$
\begin{gather}
    (\varepsilon^h, \bar{\omega})_{L^2(\Omega)} - (\nabla \times \varepsilon^h, \underline{\bar{u}})_{L^2(\Omega)} = (\varepsilon^h, {{\omega}})_{L^2(\Omega)} - (\nabla \times \varepsilon^h, \underline{{u}})_{L^2(\Omega)}, \quad \forall \varepsilon^h \in \bar{U},
\end{gather}
where the right-hand side (containing the continuous solution) may be exactly re-expressed purely in terms of the problem data by restricting the test space in the first equation of \eqref{eq:stokes_weak} to $\bar{U}$
\begin{gather}
    % (\varepsilon^h, \bar{\omega})_{L^2(\Omega)} + (\nabla \times \varepsilon^h, \underline{\bar{u}})_{L^2(\Omega)} = (\varepsilon^h, \underline{{\omega}})_{L^2(\Omega)} - (\varepsilon^h, \nabla \times \underline{{u}})_{L^2(\Omega)} + \int_{\partial \Omega} \varepsilon^h \nabla \times \underline{u} \times \hat{n} \: \mathrm{d} \Gamma, \quad \forall \varepsilon^h \in H(\mathrm{curl}, \Omega) \\
    (\varepsilon^h, \bar{\omega})_{L^2(\Omega)} - (\nabla \times \varepsilon^h, \underline{\bar{u}})_{L^2(\Omega)} = \int_{\partial \Omega} \varepsilon^h u_{||} \: \mathrm{d} \Gamma_t, \quad \forall \varepsilon^h \in \bar{U}.
\end{gather}
Similarly, we can restrict the test space to $\bar{V}$ in the second equation of \eqref{eq:stokes_weak} and re-write the second equation of \eqref{eq:stokes_projector} as 
% \begin{gather}
%     (\underline{v}^h,  \nabla \times \bar{\omega})_{L^2(\Omega)} + (\nabla \cdot \underline{v}^h, \bar{p})_{L^2(\Omega)} = (\underline{v}^h,  \nabla \times \underline{{\omega}})_{L^2(\Omega)} + (\nabla \cdot \underline{v}^h, {p})_{L^2(\Omega)}, \quad \forall \underline{v}^h \in \bar{V}.
% \end{gather}
% Substitute $\underline{\tilde{\omega}} = \bar{\omega} - \omega$ and $\tilde{p} = \bar{p} - p$
% \begin{gather}
%     (v^h,  \nabla \times \bar{\omega})_{L^2(\Omega)} + (\nabla \cdot v^h, \bar{p})_{L^2(\Omega)} = (v^h,  \nabla \times \underline{{\omega}})_{L^2(\Omega)} + (\nabla \cdot v^h, {p})_{L^2(\Omega)}, \quad \forall v^h \in H(\mathrm{div}, \Omega).
% \end{gather}
% Apply integration by parts on the right gives
\begin{gather}
    % (v^h,  \nabla \times \bar{\omega})_{L^2(\Omega)} + (\nabla \cdot v^h, \bar{p})_{L^2(\Omega)} = (v^h,  \nabla \times \underline{{\omega}})_{L^2(\Omega)} - (v^h, \nabla {p})_{L^2(\Omega)}, \quad \forall v^h \in H(\mathrm{div}, \Omega) \\
    -(\underline{v}^h,  \nabla \times \bar{\omega})_{L^2(\Omega)} + (\nabla \cdot \underline{v}^h, \bar{\pbar})_{L^2(\Omega)} = \int_{\partial \Omega} \underline{v}^h \hat{\pbar} \: \mathrm{d} \Gamma_p - (\underline{v}^h, \underline{f})_{L^2(\Omega)}, \quad \forall \underline{v}^h \in \bar{V}.
\end{gather}
% where we ignore the boundary term for the pressure as the pressure is simply a Lagrange multiplier for the velocity and as such its true value is not required. Hence the pressure is unique up to a constant. 
Lastly, we apply the same recipe for the third equation of \eqref{eq:stokes_projector} and we get
\begin{gather}
    (\eta^h, \nabla \cdot \underline{\bar{u}})_{L^2(\Omega)} = (\eta^h, \nabla \cdot \underline{{u}})_{L^2(\Omega)} = 0, \quad \forall \eta^h \in \bar{W}.
\end{gather}
We thus arrive at the final discrete system for the Stokes flow problem
\begin{equation}
    \begin{array}{cccccccl}
(\varepsilon^h, \bar{\omega})_{L^2(\Omega)} & - & (\nabla \times \varepsilon^h,  \underline{\bar{u}})_{L^2(\Omega)} & & & = &\displaystyle\int_{\partial \Omega} \varepsilon^h u_{||} \: \mathrm{d} \Gamma_t,  & \forall \varepsilon^h \in \bar{U} \\[1ex]
-(\underline{v}^h, \nabla \times \bar{\omega})_{L^2(\Omega)} 
& & + & & (\nabla \cdot \underline{v}^h, \bar{\pbar})_{L^2(\Omega)} 
&= &\displaystyle\int_{\partial \Omega} \underline{v}^h \hat{\pbar} \: \mathrm{d} \Gamma_p - (\underline{v}^h, \underline{f})_{L^2(\Omega)}, 
& \forall \underline{v}^h \in \bar{V} \\[1ex]
& & (\eta^h, \nabla \cdot \underline{\bar{u}})_{L^2(\Omega)} & & &= & 0, 
&  \forall \eta^h \in \bar{W}
\end{array},
\label{eq:stokes_weak_discrete}
\end{equation}
where the tangential velocity ($u_{||}$) and pressure ($\hat{p}$) boundary conditions are imposed weakly, and the normal velocity ($u_{\perp}$) and vorticity boundary conditions are imposed strongly. For this mixed formulation to be well-posed, we require that the bilinear forms satisfy the following two conditions \cite{Boffi2013MixedApplications, Ern2021FiniteI}
\begin{enumerate}
    \item the bilinear form $(\varepsilon^h, \bar{\omega})_{L^2(\Omega)}$ is bounded and coercive in the null space of the curl
    \item the bilinear forms $(\nabla \times \varepsilon^h,  \underline{\bar{u}})_{L^2(\Omega)}$ and $(\nabla \cdot \underline{v}^h, \bar{\pbar})_{L^2(\Omega)}$ satisfy the \emph{inf-sup} conditions
    \[\adjustlimits\inf_{\underline{\bar{u}} \in \bar{V}} \sup_{\varepsilon^h \in \bar{U}} 
     \frac{(\nabla \times \varepsilon^h,  \underline{\bar{u}})_{L^2(\Omega)}}{\lVert \varepsilon^h \rVert_{H(\mathrm{curl})} \lVert \underline{\bar{u}}\rVert_{H(\mathrm{div})}} > \beta_{\omega} > 0, \quad \quad \quad \quad 
     \adjustlimits\inf_{{\bar{\pbar}} \in \bar{W}} \sup_{\underline{v}^h \in \bar{V}} 
     \frac{(\nabla \cdot \underline{v}^h, \bar{\pbar})_{L^2(\Omega)}}{\lVert \underline{v}^h \rVert_{H(\mathrm{div})} \lVert {\bar{\pbar}}\rVert_{L^2}} > \beta_{\underline{u}} > 0
     \]
\end{enumerate}
We provide the formal proofs in \ref{sec:well_posed}. In summary, the mixed formulation is shown to be well posed on a contractible domain with $\hat{\pbar}$ prescribed on a portion of the boundary. However, if we only prescribe the normal velocity on the boundary and $\Gamma_p = \varnothing$, then the second \emph{inf-sup} condition is only satisfied for the pressure in $\bar{W}\backslash \mathbb{R}$. Nevertheless, the velocity and vorticity solutions are unique and bounded in either scenario.

As the projector arises from a saddle-point formulation, the projected field inherits a natural orthogonal decomposition structure. Let $\bar{\omega} \in \bar{U}$ denote the projection of $\omega \in H(\mathrm{curl}, \Omega)$ defined through~\eqref{eq:stokes_projector}. This projection admits a discrete Helmholtz decomposition,
\[
\bar{\omega} = \bar{\omega}_{0} + \bar{\omega}_{\perp},
\qquad
\bar{\omega}_{0} \in \mathrm{Ker}^{h}(\nabla \times),
\quad
\bar{\omega}_{\perp} \in \mathrm{Ker}^{h}(\nabla \times)^{\perp}.
\]
The component $\bar{\omega}_{0}$ is the $L^{2}$-best approximation of $\omega$ within $\mathrm{Ker}^{h}(\nabla \times)$,
\[
\|\bar{\omega}_{0} - \omega\|_{L^{2}}
=
\inf_{\varepsilon^{h}_{0} \in \mathrm{Ker}^{h}(\nabla \times)}
\|\varepsilon^{h}_{0} - \omega\|_{L^{2}},
\]
while $\bar{\omega}_{\perp}$ provides the $L^{2}$-best approximation of $\omega$ in the orthogonal complement, characterised by
\[
\|\nabla \times \bar{\omega}_{\perp} - \nabla \times \omega\|_{L^{2}}
=
\inf_{\varepsilon^{h}_{\perp} \in \mathrm{Ker}^{h}(\nabla \times)^{\perp}}
\|\nabla \times \varepsilon^{h}_{\perp} - \nabla \times \omega\|_{L^{2}}.
\]
A complete proof of this optimality property is provided in \ref{app:optimal}.

Using this projector, we decompose the ambient spaces into {resolved} and {unresolved} scales. The resolved scales correspond to the range of the projector,
\(
\mathrm{Range}(\mathcal{P}) = \bar{U} \times \bar{V} \times \bar{W},
\)
while the unresolved scales lie in its null space,
\(
\mathrm{Ker}(\mathcal{P}) = U' \times V' \times W',
\quad
(\omega',\underline{u}',p') := (\omega-\bar{\omega},\, \underline{u}-\underline{\bar{u}},\, p-\bar{p}).
\)
Furthermore, the resolved and unresolved scales are orthogonal in the following form
\begin{equation}
    \begin{array}{cccccl}
        \bilnear{\varepsilon^h}{{\omega'}} & - \bilnear{\nabla \times \varepsilon^h}{{\underline{u}'}} & & = &0,  & \forall \varepsilon^h \in \bar{U}, \forall (\omega', \underline{u}') \in U' \times V' \\[1.5ex]
        -\bilnear{\underline{v}^h}{\nabla \times {\omega'}}
        &+ & \bilnear{\nabla \cdot \underline{v}^h}{{\pbar}'}
        &= & \underline{0}, 
        & \forall \underline{v}^h \in \bar{V},\forall (\omega', \pbar') \in U' \times W' \\[1.5ex]
        &  \bilnear{\eta^h}{\nabla \cdot \underline{{u}}'} & & = & 0, 
        &  \forall \eta^h \in \bar{W},\forall \underline{u}' \in V'
    \end{array},
    \label{eq:stokes_projector_ortho}
\end{equation}
Thus, the continuous spaces decompose as
\[
H(\mathrm{curl},\Omega) = \bar{U} \oplus U', \qquad
H(\mathrm{div},\Omega) = \bar{V} \oplus V', \qquad
L^{2}(\Omega) = \bar{W} \oplus W'.
\]

\section{Navier-Stokes}
\label{sec:navier_stokes}
Having established the structure of our optimal projector for the Stokes problem, we now introduce the VMS formulation of the Navier–Stokes equations built around a projector that preserves this optimal Stokes structure. 
\subsection{Continuous mixed formulation}
Much like the Stokes flow problem, we express the incompressible Navier-Stokes equations in a mixed formulation by defining the vorticity $\omega$ and further rewriting the convective term in the Lamb (Rotational) form. The continuous problem in a domain $\Omega$ over a time period $\mathcal{T} = (0, T]$ reads
\begin{align}
    {\omega} - \nabla \times \underline{u} &= 0, \quad \text{in } \Omega \times \mathcal{T}  \\
    \parddx{\underline{u}}{t} + \omega \times \underline{u} + \frac{1}{Re} \nabla \times {\omega} + \nabla P &= \underline{0}, \quad \text{in } \Omega \times \mathcal{T} \label{eq:NS_momentum}\\
    \nabla \cdot \underline{u} &= 0, \quad \text{in } \Omega \times \mathcal{T} \\
    \underline{u} \cdot \hat{n} = u_{\perp}, \quad \text{on} \: \Gamma_n \times \mathcal{T}, &\quad \quad \omega = \hat{\omega}, \quad \text{on} \: \Gamma_{\omega}  \times \mathcal{T} \\
    P = \hat{P}, \quad \text{on} \: \Gamma_P \times \mathcal{T}, &\quad \quad {\underline{u}} \times \hat{n} = u_{\parallel}, \quad \text{on} \: \Gamma_t \times \mathcal{T} \\
    \text{with } \Gamma_n \cup \Gamma_P = \Gamma_{\omega} \cup \Gamma_t &= \partial \Omega, \quad \text{and} \: \Gamma_n \cap \Gamma_P = \Gamma_{\omega} \cap \Gamma_t = \varnothing,
\end{align}
where $P := \pbar + \frac{1}{2} \underline{u}^2$ is the total or Bernoulli pressure. The system is supplemented with an initial condition for the velocity and vorticity. 

At the continuous level, the formulation exhibits strong conservation and dissipation properties. Mass conservation is enforced by the divergence-free condition. Kinetic energy $\mathcal{K} := \frac{1}{2} \bilnear{\underline{u}}{\underline{u}}$ and enstrophy $\mathcal{E} := \frac{1}{2} \bilnear{\omega}{\omega}$ evolve in time through viscous dissipation: the kinetic energy decays proportionally to the enstrophy, and the enstrophy in turn dissipates via palinstrophy $\mathscr{P} := \frac{1}{2} \bilnear{\nabla \times \omega}{\nabla \times \omega}$. Furthermore, kinetic energy and enstrophy are exactly conserved in the inviscid limit ($Re \rightarrow 
\infty$). These relations emerge naturally from the structure of the equations and are captured through balance laws derived under suitable assumptions (e.g., periodic domains or vanishing boundary fluxes). Moreover, the total vorticity $\mathcal{W} := \bilnear{1}{\omega}$ is also conserved, a result directly linked to the circulation theorem. The four conserved quantities are listed as follows (see \cite{Zhang2024AConditions} for further details)
\[\nabla \cdot \underline{u} = 0, \quad \quad \parddx{\mathcal{K}}{t} + \frac{2}{Re} \mathcal{E} = 0, \quad\quad \parddx{\mathcal{E}}{t} + \frac{2}{Re} \mathscr{P} = 0, \quad\quad \mathcal{W} = \int_{\partial \Omega} u_{||} \: \mathrm{d} \Gamma_t.\]

% The MEEVC discretisation \cite{Zhang2024AConditions} is a mixed spectral element scheme constructed to preserve the structure and invariants of the continuous problem. The three physical quantities are defined in the aforementioned Sobolev spaces as previously shown for Stokes flow. This choice ensures that discrete differential operators mimic their continuous counterparts, thereby enabling the conservation of mass, energy, enstrophy, and vorticity at the discrete level.
The essence of the MEEVC scheme is to conserve these symmetries at the discrete level. We propose to do the same while simultaneously supplementing the scheme with VMS analysis. 
\begin{remark}
    \label{rmk:under_int}
    % While the MEEVC scheme \cite{Zhang2024AConditions} is designed to conserve enstrophy at the discrete level, this property depends critically on the accuracy of numerical integration. In particular, the non-linear convective term and all other bilinear forms must be integrated exactly to ensure enstrophy conservation. In practice, this requires either exact integration or the use of sufficiently high-order quadrature rules. If standard quadrature is used (i.e., of the same polynomial order as the basis functions), the conservation of enstrophy is no longer guaranteed. 
    To recover all the conservation properties of the MEEVC scheme, the bilinear and trilinear terms in the weak form must be evaluated using exact or over-integration. However, when these terms are under-integrated, enstrophy conservation is no longer guaranteed. 
    % In this work, we deliberately use Gauss–Lobatto quadrature rules of the same order as the basis functions, which generally results in under-integration. As a result, the presented simulations (both the basic Galerkin and VMS formulations) do not strictly conserve enstrophy. We make this choice intentionally, as the under-integrated case presents a more critical and practically relevant scenario. If strict enstrophy conservation is required, the methodology remains fully applicable with a simple modification to over-integrate the bilinear and trilinear forms. 
\end{remark}

In line with the MEEVC scheme, we spatially define the three physical quantities in the aforementioned Sobolev spaces as previously shown for Stokes flow. 
However, we first perform a temporal discretisation to reduce the time-dependent problem to an effective steady problem that can be solved every time step. We employ the simple second-order accurate symplectic time integrator, namely the Crank-Nicolson scheme
\begin{gather}
    \frac{\underline{{u}}_{n + 1} - \underline{{u}}_n}{\Delta t} = -{\omega}_{n + \frac{1}{2}} \times {\underline{u}}_{n + \frac{1}{2}} - \frac{1}{Re} \nabla \times {\omega}_{n + \frac{1}{2}} - \nabla {P}_{n + \frac{1}{2}}, \quad \quad \text{where } [\cdot]_{n + \frac{1}{2}} = \frac{1}{2} ([\cdot]_{n + 1} + [\cdot]_n).
\end{gather}
We henceforth directly work with the time-discretised equations, with the effective steady problem to be solved every time step given as
\begin{align}
    {\omega}_{n + 1} - \nabla \times {\underline{u}}_{n + 1} &= 0 \\
    \frac{-1}{2 Re} \nabla \times {\omega}_{n + 1} - \frac{1}{\Delta t} {\underline{u}}_{n + 1} + \nabla {P}_{n + \frac{1}{2}} &= \frac{1}{2 Re} \nabla \times {\omega}_{n} - \frac{1}{\Delta t} {\underline{u}}_{n} + {\omega}_{n + \frac{1}{2}} \times {\underline{u}}_{n + \frac{1}{2}} \\
    \nabla \cdot {\underline{u}}_{n + 1} &= 0.
\end{align}
We derive an infinite-dimensional weak form of the problem as follows
% \begin{gather}
%     \bilnear{\omega}{{\omega}_{n + 1}} - \bilnear{\nabla \times \omega}{{\underline{u}}_{n + 1}} = \int_{\partial \Omega} \omega u_{\parallel} \: \mathrm{d}\Gamma_t, \quad \forall \omega \in H(\text{curl}, \Omega)  \\
%     \begin{split}
%         \frac{1}{\Delta t} \bilnear{\underline{v}}{{\underline{u}}_{n + 1}} + \frac{1}{2 Re} \bilnear{\underline{v}}{\nabla \times {\omega}_{n + 1}} + \bilnear{\nabla \cdot \underline{v}}{{P}_{n + \frac{1}{2}}} \\= \frac{1}{\Delta t} \bilnear{\underline{v}}{{\underline{u}}_{n}} - \frac{1}{2 Re} \bilnear{\underline{v}}{\nabla \times {\omega}_{n}} + a({\underline{u}}_{n + \frac{1}{2}}, {\omega}_{n + \frac{1}{2}}, \underline{v}) + \int_{\partial \Omega} \hat{P}_t \underline{v} \: \mathrm{d} \Gamma_P
%     \end{split}, \quad \forall \underline{v} \in H(\text{div}, \Omega) \\
%     \bilnear{\eta}{\nabla \cdot {\underline{u}}_{n + 1}} = 0, \quad \forall \eta \in L^2(\Omega).
% \end{gather}
{\small
\begin{equation}
    \begin{gathered}
    \begin{array}{cccccl}
        \bilnear{\varepsilon}{{\omega}_{n + 1}} & -\bilnear{\nabla \times \varepsilon}{{\underline{u}}_{n + 1}} & & = &\displaystyle\int_{\partial \Omega} \varepsilon u_{||} \: \mathrm{d} \Gamma_t,  & \forall \omega \in H(\text{curl}, \Omega) \\[1.5ex]
        \displaystyle\frac{-1}{2 Re} \bilnear{\underline{v}}{\nabla \times {\omega}_{n + 1}} 
        & -\displaystyle\frac{1}{\Delta t} \bilnear{\underline{v}}{{\underline{u}}_{n + 1}} & + \bilnear{\nabla \cdot \underline{v}}{{P}_{n + \frac{1}{2}}}
        &= & \begin{aligned}
            \mathcal{F}(\underline{v}, \underline{{u}}_n, {{\omega}}_n, \hat{P}) + \\a(\underline{v}, {\underline{u}}_{n + \frac{1}{2}}, {\omega}_{n + \frac{1}{2}})
        \end{aligned} \:, 
        & \forall \underline{v} \in H(\text{div}, \Omega) \\[4ex]
        &  \bilnear{\eta}{\nabla \cdot {\underline{u}}_{n + 1}} & & = & 0, 
        &  \forall \eta \in L^2(\Omega)
    \end{array} \\ \\
\text{where } \quad \begin{aligned}
                &\mathcal{F}(\underline{v}, \underline{{u}}_n, {{\omega}}_n, \hat{P}) := \frac{-1}{\Delta t} \bilnear{\underline{v}}{{\underline{u}}_{n}} + \frac{1}{2 Re} \bilnear{\underline{v}}{\nabla \times {\omega}_{n}} + \int_{\partial \Omega} \underline{v} \hat{P} \: \mathrm{d} \Gamma_P, \\
                &a(\underline{v}, {\underline{u}}_{n + \frac{1}{2}}, {\omega}_{n + \frac{1}{2}}) := \bilnear{\underline{v}}{{\omega}_{n + \frac{1}{2}} \times {\underline{u}}_{n + \frac{1}{2}}}
            \end{aligned}
\end{gathered}
\label{eq:naver_stokes_weak_discrete}
\end{equation}
}
The system undoubtedly shares similarities with the system in \eqref{eq:stokes_weak}. However, we cannot variationally characterise the full system through a saddle point problem as we did for Stokes, primarily due to the presence of the non-linear trilinear form $a(\underline{v}, {\underline{u}}_{n + \frac{1}{2}}, {\omega}_{n + \frac{1}{2}})$. Consequently, we note that the basic Galerkin discretisation of this system will \emph{not} yield a desirable projection of the exact solution. Hence, we employ VMS analysis wherein we leverage the parts that do fit in the aforementioned variational formulation and formulate an approach that does yield a close approximation of the projected solution.

\subsection{Variational Multiscale formulation}
We define our \emph{optimal projector} $\mathcal{P} : H(\mathrm{curl}, \Omega) \times H(\mathrm{div}, \Omega) \times L^2(\Omega) \rightarrow \bar{U} \times \bar{V} \times \bar{W}$ for Navier-Stokes at an arbitrary time level, by slightly augmenting the the energy functional from \eqref{eq:stokes_energy} as follows
\begin{gather}
    \begin{split}
        \mathcal{I_{NS}}({\bar{\omega}} - {\omega}, \underline{\bar{u}} - \underline{{u}}, \bar{P} - {P}; 0, 0, 0) := \int_{\Omega} \frac{1}{2 Re} \left( \frac{1}{2} |{\bar{\omega}} - {{\omega}}|^2 - (\underline{\bar{u}} - \underline{{u}}) \cdot \nabla \times ({\bar{\omega}} - {{\omega}}) \right) + \\\left(\frac{-1}{2 \Delta t}|\underline{\bar{u}} - \underline{{u}}|^2 + (\bar{P} - {P}) \nabla \cdot (\underline{\bar{u}} - \underline{{u}}) \right) \: \mathrm{d}\Omega,
    \end{split}
\end{gather}
where the discrete solution $({\bar{\omega}}, \underline{\bar{u}}, \bar{P}) \in \bar{U} \times \bar{V} \times \bar{W}$ is a stationary point of the functional for the given time discretised exact spatial solution $({{\omega}}, \underline{{u}}, {P}) \in H(\mathrm{curl}, \Omega) \times H(\mathrm{div}, \Omega) \times L^2(\Omega)$. In this context, the term optimal denotes a choice that is appropriately aligned with the structure of the governing PDE. We subsequently show that this choice allows the formulation to inherit the structure-preserving character of the discretisation. Taking variations gives
% \begin{gather}
%     \bilnear{\varepsilon^h}{\bar{\omega} - {\omega}} - \bilnear{\nabla \times \varepsilon^h}{\bar{\underline{u}} - {\underline{u}}} = 0, \quad \forall \varepsilon^h \in \bar{U} \subset H(\mathrm{curl}, \Omega) \\
%     \frac{1}{\Delta t} \bilnear{\underline{v}^h}{\underline{\bar{u}} -{\underline{u}}} + \frac{1}{2 Re} \bilnear{\underline{v}^h}{\nabla \times (\bar{\omega} - {\omega})} + \bilnear{\nabla \cdot \underline{v}^h}{\bar{P} - {P}} = 0, \quad \forall \underline{v}^h \in  \bar{V} \subset H(\text{div}, \Omega) \\
%     \bilnear{\eta^h}{\nabla \cdot (\bar{\underline{u}} - \underline{{u}})} = 0, \quad \forall \eta^h \in \bar{W} \subset L^2(\Omega)
% \end{gather}
\begin{equation}
    \begin{array}{cccccl}
        \bilnear{\varepsilon^h}{\bar{\omega} - {\omega}} & -\bilnear{\nabla \times \varepsilon^h}{\bar{\underline{u}} - {\underline{u}}} & & = &0,  & \forall \varepsilon^h \in \bar{U} \\[1.5ex]
        \displaystyle\frac{-1}{2 Re} \bilnear{\underline{v}^h}{\nabla \times (\bar{\omega} - {\omega})}
        & -\displaystyle\frac{1}{\Delta t} \bilnear{\underline{v}^h}{\underline{\bar{u}} -{\underline{u}}} & + \bilnear{\nabla \cdot \underline{v}^h}{\bar{P} - {P}}
        &= & \underline{0}, 
        & \forall \underline{v}^h \in \bar{V} \\[3ex]
        &  \bilnear{\eta^h}{\nabla \cdot (\bar{\underline{u}} - \underline{{u}})} & & = & 0, 
        &  \forall \eta^h \in \bar{W}
    \end{array},
    \label{eq:NS_proj}
\end{equation}
where well posedness is again guaranteed if the conditions mentioned in \Cref{sec:stokes} are satisfied along with an additional condition on the boundedness and coercivity of the bilinear form $\bilnear{\underline{v}^h}{\underline{\bar{u}} -{\underline{u}}}$ in $\mathrm{Ker}^h(\nabla \cdot)$, which are trivially satisfied. For this projector, we have the following orthogonality conditions
{\small
\begin{equation}
    \begin{array}{ccllll}
        \bilnear{\varepsilon^h}{{\omega}'} & -\bilnear{\nabla \times \varepsilon^h}{{\underline{u}}'} & & = &0,  & \forall \varepsilon^h \in \bar{U}, \forall (\omega', \underline{u}') \in U' \times V' \\[1.5ex]
        \displaystyle\frac{-1}{2 Re} \bilnear{\underline{v}^h}{\nabla \times {\omega'}}
        & -\displaystyle\frac{1}{\Delta t} \bilnear{\underline{v}^h}{\underline{{u}}'} & + \bilnear{\nabla \cdot \underline{v}^h}{{P}'}
        &= & \underline{0}, 
        & \forall \underline{v}^h \in \bar{V}, \forall (\omega', \underline{u}', P') \in U' \times V' \times W'  \\[3ex]
        &  \bilnear{\eta^h}{\nabla \cdot {\underline{u}'}} & & = & 0, 
        &  \forall \eta^h \in \bar{W}, \forall \underline{u}' \in V'
    \end{array}.
    \label{eq:NS_ortho}
\end{equation}
}

Noting the projector in \eqref{eq:NS_proj}, we derive the resolved-scale equations, by restricting the test space $\varepsilon = \varepsilon^h \in \bar{U}, \underline{v} = \underline{v}^h \in \bar{V}, \eta = \eta^h \in \bar{W}$ and decomposing the solution as ${\omega} = \bar{\omega} + \omega', {\underline{u}} = \bar{\underline{u}} + \underline{u}', {P} = \bar{P} + P'$. Filling in this decomposition into \eqref{eq:naver_stokes_weak_discrete} and employing the restricted test space leads to
% \begin{gather}
%     \bilnear{\varepsilon^h}{\bar{\omega}_{n + 1}} - \bilnear{\nabla \times \varepsilon^h}{\bar{\underline{u}}_{n + 1}} = \int_{\partial \Omega} \varepsilon^h u_{\parallel} \: \mathrm{d}\Gamma_t, \quad \forall \varepsilon^h \in \bar{U} \\
%     \begin{split}
%         \frac{1}{\Delta t} \bilnear{\underline{v}^h}{\bar{\underline{u}}_{n + 1}} + \frac{1}{2 Re} \bilnear{\underline{v}^h}{\nabla \times \bar{\omega}_{n + 1}} + \bilnear{\nabla \cdot \underline{v}^h}{\bar{P}_{n + \frac{1}{2}}} = \\\frac{1}{\Delta t} \bilnear{\underline{v}^h}{\bar{\underline{u}}_{n} + \underline{u}'_{n}} - \frac{1}{2 Re} \bilnear{\underline{v}^h}{\nabla \times (\bar{\omega}_{n} + \omega'_n)} + a({\underline{u}}_{n + \frac{1}{2}}, {\omega}_{n + \frac{1}{2}}, \underline{v}^h) + \int_{\partial \Omega} \hat{P}_t \underline{v}^h \hat{n} \: \mathrm{d} \Gamma_P
%     \end{split}, \quad \forall \underline{v}^h \in \bar{V} \\
%     \bilnear{\eta^h}{\nabla \cdot \bar{\underline{u}}_{n + 1}} = 0, \quad \forall \eta^h \in \bar{W}
% \end{gather}
{\small
\begin{equation}
    \begin{array}{cccccl}
        \bilnear{\varepsilon^h}{\bar{\omega}_{n + 1}} & -\bilnear{\nabla \times \varepsilon^h}{\bar{\underline{u}}_{n + 1}} & & = &\displaystyle\int_{\partial \Omega} \varepsilon^h u_{||} \: \mathrm{d} \Gamma_t,  & \forall \varepsilon^h \in \bar{U} \\[1.5ex]
        \displaystyle\frac{-1}{2 Re} \bilnear{\underline{v}^h}{\nabla \times \bar{\omega}_{n + 1}}
        & -\displaystyle\frac{1}{\Delta t} \bilnear{\underline{v}^h}{\bar{\underline{u}}_{n + 1}} & + \bilnear{\nabla \cdot \underline{v}^h}{\bar{P}_{n + \frac{1}{2}}}
        &= & \begin{aligned}
            \mathcal{F}(\underline{v}^h, \underline{{u}}_n, {{\omega}}_n, \hat{P}) + \\a(\underline{v}^h, {\underline{u}}_{n + \frac{1}{2}}, {\omega}_{n + \frac{1}{2}})
        \end{aligned}, 
        & \forall \underline{v}^h \in \bar{V} \\[4ex]
        &  \bilnear{\eta^h}{\nabla \cdot {\underline{\bar{u}}}_{n + 1}} & & = & 0, 
        &  \forall \eta^h \in \bar{W}
    \end{array},
    \label{eq:NS_resolved}
\end{equation}
}
where we have eliminated the fine-scale terms using the orthogonality conditions in \eqref{eq:NS_ortho}. Similarly, by restricting the test space to $\varepsilon = \varepsilon' \in {U}', \underline{v} = \underline{v}' \in {V}', \eta = \eta' \in {W}'$ we get the unresolved-scale equations
% \begin{gather}
%     \bilnear{\omega'}{{\omega}'_{n + 1}} - \bilnear{\nabla \times \omega'}{{\underline{u}}'_{n + 1}} = \int_{\partial \Omega} \omega' u_{\parallel} \: \mathrm{d}\Gamma_t, \quad \forall \omega' \in {U}' \\
%     \begin{split}
%         \frac{1}{\Delta t} \bilnear{\underline{v}'}{{\underline{u}}'_{n + 1}} + \frac{1}{2 Re} \bilnear{\underline{v}'}{\nabla \times {\omega}'_{n + 1}} + \bilnear{\nabla \cdot \underline{v}'}{{P}'_{n + \frac{1}{2}}} = \\\frac{1}{\Delta t} \bilnear{\underline{v}'}{\bar{\underline{u}}_{n} + \underline{u}'_{n}} - \frac{1}{2 Re} \bilnear{\underline{v}'}{\nabla \times (\bar{\omega}_{n} + \omega'_n)} + a({\underline{u}}_{n + \frac{1}{2}}, {\omega}_{n + \frac{1}{2}}, \underline{v}') + \int_{\partial \Omega} \hat{P}_t \underline{v}' \hat{n} \: \mathrm{d} \Gamma_P
%     \end{split}, \quad \forall \underline{v}' \in {V}' \\
%     \bilnear{\eta'}{\nabla \cdot {\underline{u}}'_{n + 1}} = 0, \quad \forall \eta' \in {W}'
% \end{gather}
{\small
\begin{equation}
    \begin{array}{cccccl}
        \bilnear{\varepsilon'}{{\omega}'_{n + 1}} & -\bilnear{\nabla \times \varepsilon'}{{\underline{u}}'_{n + 1}} & & = &\displaystyle\int_{\partial \Omega} \varepsilon' u_{\parallel} \: \mathrm{d}\Gamma_t,  & \forall \varepsilon' \in {U}' \\[1.5ex]
        \displaystyle\frac{-1}{2 Re} \bilnear{\underline{v}'}{\nabla \times {\omega}'_{n + 1}}
        & -\displaystyle\frac{1}{\Delta t} \bilnear{\underline{v}'}{{\underline{u}}'_{n + 1}} & + \bilnear{\nabla \cdot \underline{v}'}{{P}'_{n + \frac{1}{2}}}
        &= & \begin{aligned}
            \mathcal{F}(\underline{v}', \underline{{u}}_n, {{\omega}}_n, \hat{P}) + \\a(\underline{v}', {\underline{u}}_{n + \frac{1}{2}}, {\omega}_{n + \frac{1}{2}})
        \end{aligned}, 
        & \forall \underline{v}' \in {V}' \\[4ex]
        &  \bilnear{\eta'}{\nabla \cdot {\underline{u}}'_{n + 1}} & & = & 0, 
        &  \forall \eta; \in {W}'
    \end{array}.
    \label{eq:NS_unresolved}
\end{equation}
}
Note that the terms $\mathcal{F}(\underline{v}, \underline{{u}}_n, \underline{{\omega}}_n, \hat{P}) + a(\underline{v}, {\underline{u}}_{n + \frac{1}{2}}, {\omega}_{n + \frac{1}{2}})$ %= \mathcal{F}(\underline{v}, \underline{\bar{u}}_n + \underline{{u}}'_n, \underline{\bar{\omega}}_n + \underline{{\omega}}'_n, \hat{P}) + a(\underline{v}, {\underline{\bar{u}}}_{n + \frac{1}{2}} + {\underline{u}}'_{n + \frac{1}{2}}, {\underline{\bar{\omega}}}_{n + \frac{1}{2}} + {\omega'}_{n + \frac{1}{2}})$
that appear on the right-hand side of \eqref{eq:NS_resolved} and \eqref{eq:NS_unresolved} include both the resolved and unresolved components, hence the two equations are coupled through the non-linear term.

\subsection{Properties of the formulation}\label{sec:VMS_prop}
We summarise the full structure of this VMS formulation through an illustration of the Hilbert (sub)complex in \Cref{fig:deRham}, defined over a 2D contractible domain. At the continuous level, the spaces that form the de Rham complex are split into closed subspaces corresponding to the null spaces ($\mathrm{Ker}(\cdot)$)/ranges ($\mathrm{Range}(\cdot)$) of the differential operators, and their open $L^2$-orthogonal complements ($\mathrm{Ker}(\cdot)^{\perp}$ or $\mathrm{Range}(\cdot)^{\perp}$) \cite{Arnold2010FINITESTABILITY} \cite[\S 2.1.8]{PalhadaSilvaClerigo2013HighQuadrilaterals}. The finite-dimensional subspaces ($\bar{U}, \bar{V}, \bar{W}$) are constructed to contain a discrete null space $\mathrm{Ker}^h(\cdot)$ and its discrete orthogonal complement $\mathrm{Ker}^h(\cdot)^{\perp}$. We note that the finite-dimensional null space is contained within the infinite-dimensional subspace $\mathrm{Ker}^h(\cdot) \subset \mathrm{Ker}(\cdot)$. However, this is not true for the orthogonal complement. Formally, the complementary open sets (${U}', {V}', {W}'$) include everything outside the sets ($\bar{U}, \bar{V}, \bar{W}$) and hence represent the unresolved parts of the full Hilbert spaces not captured by the coarse-scale approximation. To approximate these fine-scale components, we introduce (${U}_k', {V}_k', {W}_k'$) which are finite-dimensional subspaces constructed to be exactly orthogonal to ($\bar{U}, \bar{V}, \bar{W}$) in the norm defined by \eqref{eq:NS_ortho}. These fine-scale spaces are constructed via $p$-refinement of the resolved spaces, with the polynomial enrichment indicated by the parameter $k$. As $k$ increases, the fine-scale approximation space more fully captures the structure of the infinite-dimensional unresolved subspace (see \cite{Shrestha2025OptimalProblems} for further details). 
\begin{figure}[htp]
    \centering
    \begin{tikzpicture}

% Define colors for sections
\definecolor{darkgray}{RGB}{120,120,120}
\definecolor{gray}{RGB}{211,211,211}
\definecolor{lightgray}{RGB}{255,255,255}
\definecolor{black}{RGB}{0,0,0}
\tikzmath{\spacing = 6; \radiOut = 1.8; \radiMid = 1; \radiIn = 1; \ellipR = 1.6;}

% Draw ellipses for each space
\foreach \i/\N in {0/{U}, 1/{V}, 2/{W}} {
    % Closed space
    \filldraw[pattern={Lines[distance=0.8mm, angle=-45, line width=0.05mm]}, opacity = 0.4] (\i*\spacing + 0.6*\radiMid, 0) ellipse[x radius=0.8*\radiMid*1.4, y radius=0.8*\radiMid*1.2];
    \draw[thick] (\i*\spacing + 0.6*\radiMid, 0) ellipse[x radius=0.8*\radiMid*1.4, y radius=0.8*\radiMid*1.2];
    \fill[gray] (\i*\spacing + 0.6*\radiMid, 0) ellipse[x radius=0.4*\radiMid*1.8, y radius=0.4*\radiMid*1.2];
    \draw[thick] (\i*\spacing + 0.6*\radiMid, 0) ellipse[x radius=0.4*\radiMid*1.8, y radius=0.4*\radiMid*1.2];
    \node at (\i*\spacing + 0.6*\radiMid, 0) {$\bar{\N}$};
    \node at (\i*\spacing + 0.5*\radiMid, -0.75) {${\N}'_k$};
}

\foreach \i in {0, 1, 2} {
    % Outer ellipse
    \draw[thick, dashed, opacity = 0.5] (\i*\spacing, 0) ellipse[x radius=\radiOut, y radius=\radiOut*\ellipR];
    % Inner ellipse
    % \fill[gray] (\i*\spacing, 0) ellipse[x radius=\radiMid, y radius=\radiMid*\ellipR];
    % \draw[thick] (\i*\spacing, 0) ellipse[x radius=\radiMid, y radius=\radiMid*\ellipR];
    % Core ellipse
    % \fill[gray, opacity = 0.5] (\i*\spacing, 0) ellipse[x radius=\radiIn, y radius=\radiIn*\ellipR];
    \draw[thick, opacity = 0.2] (\i*\spacing, 0) ellipse[x radius=\radiIn, y radius=\radiIn*\ellipR];
    % Center label
    \node at (\i*\spacing + 0.15, -0.1) {\scriptsize $0$};
    \fill[black, opacity = 0.5] (\i*\spacing, 0) circle[radius = 0.08];
}
\node at (0*\spacing, 0.75*\ellipR*\radiMid) {\scriptsize $\mathrm{Ker}(\nabla \times)$};
\node at (1*\spacing, 0.75*\ellipR*\radiMid) {\scriptsize $\mathrm{Ker}(\nabla \cdot)$};
\node at (2*\spacing, 0.75*\ellipR*\radiMid) {\scriptsize $\mathrm{Range}(\nabla \cdot)$};

\node at (0*\spacing, 1.4*\ellipR*\radiMid) {\scriptsize $\mathrm{Ker}(\nabla \times)^{\perp}$};
\node at (1*\spacing, 1.4*\ellipR*\radiMid) {\scriptsize $\mathrm{Ker}(\nabla \cdot)^{\perp}$};
\node at (2*\spacing, 1.4*\ellipR*\radiMid) {\scriptsize $\mathrm{Range}(\nabla \cdot)^{\perp}$};

\node at (-\spacing*0.4, \radiOut*\ellipR + 0.3*\radiOut*\ellipR) {$\mathbb{R}$};
\node at (0*\spacing, \radiOut*\ellipR + 0.3*\radiOut*\ellipR) {$H(\mathrm{curl}, \Omega)$};
\node at (1*\spacing, \radiOut*\ellipR + 0.3*\radiOut*\ellipR) {$H(\mathrm{div}, \Omega)$};
\node at (2*\spacing, \radiOut*\ellipR + 0.3*\radiOut*\ellipR) {$L^2(\Omega)$};
\node at (2*\spacing + \spacing*0.4,  \radiOut*\ellipR + 0.3*\radiOut*\ellipR) {$0$};

\draw[right hook->, thick] (-\spacing*0.4 + 0.4, \radiOut*\ellipR + 0.3*\radiOut*\ellipR) -- (-\spacing*0.2, \radiOut*\ellipR + 0.3*\radiOut*\ellipR);
\draw[->, thick] (0*\spacing + 0.2*\spacing, \radiOut*\ellipR + 0.3*\radiOut*\ellipR) -- (0*\spacing + \spacing - 0.2*\spacing, \radiOut*\ellipR + 0.3*\radiOut*\ellipR) node[midway, above] {$\nabla \times$};
\draw[->, thick] (1*\spacing + 0.2*\spacing, \radiOut*\ellipR + 0.3*\radiOut*\ellipR) -- (1*\spacing + \spacing - 0.15*\spacing, \radiOut*\ellipR + 0.3*\radiOut*\ellipR) node[midway, above] {$\nabla \cdot$};
\draw[->, thick] (2*\spacing + 0.8, \radiOut*\ellipR + 0.3*\radiOut*\ellipR) -- (2*\spacing + \spacing*0.35, \radiOut*\ellipR + 0.3*\radiOut*\ellipR);

% Draw arrows between spaces
\foreach \i in {0, 1} {
    \draw[-, thin, dashed, opacity = 0.5] (\i*\spacing, \radiMid*\ellipR) -- (\i*\spacing + \spacing, 0) {};
    \draw[-, thin, dashed, opacity = 0.5] (\i*\spacing, -\radiMid*\ellipR) -- (\i*\spacing + \spacing, 0) {};
    \draw[-, thin, opacity = 0.5] (\i*\spacing, \radiOut*\ellipR) -- (\i*\spacing + \spacing, \radiIn*\ellipR) {};
    \draw[-, thin, opacity = 0.5] (\i*\spacing, -\radiOut*\ellipR) -- (\i*\spacing + \spacing, -\radiIn*\ellipR) {};
}

% Legend
% \node[rectangle, draw, thick, fill=gray, minimum width=0.6cm, minimum height=0.4cm] (legend1) at (2*\spacing + 0.5*\spacing, \radiOut*\ellipR + 1) {};
% \node[right] at (legend1.east) {$\mathcal{Z}$};
% % \node[rectangle, draw, thick, fill=gray, minimum width=0.6cm, minimum height=0.4cm] (legend2) at (2*\spacing + 0.5*\spacing, \radiOut*\ellipR + 0.25) {};
% % \node[right] at (legend2.east) {$\mathcal{H}^k(\mathcal{M})$};
% \node[rectangle, draw, thick, fill=lightgray, minimum width=0.6cm, minimum height=0.4cm] (legend3) at (2*\spacing + 0.5*\spacing, \radiOut*\ellipR + 0.25) {};
% \node[right] at (legend3.east) {$\mathcal{Z}^{\perp}$};

% Brace and final label
% \draw[decorate, decoration={brace, amplitude=10pt}] (14.5, 3.5) -- (14.5, 0.5) node[midway, right, xshift=5pt] {$\mathcal{Z}(d; \Lambda^k)$};

\end{tikzpicture}
    \caption{Pictorial description of the de Rham complex on a 2D contractible domain, including the depiction of the finite-dimensional resolved and unresolved subspaces}
    \label{fig:deRham}
\end{figure}

As we are working with the MEEVC scheme, we now prove that the proposed VMS formulation with our optimal projector satisfies the MEEVC properties.
Specifically, the proposed VMS formulation exactly satisfies conservation of mass of both resolved and unresolved scales
    \begin{equation}
    \nabla \cdot \bar{\underline{u}}_{n + 1} = 0, \quad \nabla \cdot {\underline{u}}'_{n + 1} = 0,
\end{equation}
it exactly satisfies conservation of total vorticity
\begin{equation}
     \underbrace{(1, \bar{\omega}_{n + 1} + \omega'_{n + 1})_{L^2(\Omega)}}_{\mathcal{W}_{n + 1}} = \int_{\partial \Omega} u_{||} \: \mathrm{d} \Gamma_t
\end{equation}
it satisfies conservation of total kinetic energy
\begin{gather}
    \begin{split}
        \underbrace{\left(\frac{1}{2} \bilnear{\underline{\bar{u}}_{n + 1} + \underline{u}'_{n + 1}}{\underline{\bar{u}}_{n + 1} + \underline{u}'_{n + 1}} \right.}_{\mathcal{K}_{n + 1}} - \underbrace{\left.\frac{1}{2}\bilnear{\underline{\bar{u}}_{n} + \underline{u}'_{n}}{\underline{\bar{u}}_{n} + \underline{u}'_{n}} \right)}_{\mathcal{K}_n} + \\\frac{ \Delta t}{Re}\underbrace{\bilnear{{\bar{\omega}}_{n + \frac{1}{2}} + \omega'_{n + \frac{1}{2}}}{{\bar{\omega}}_{n + \frac{1}{2}} + \omega'_{n + \frac{1}{2}}}}_{2\mathcal{E}_{n + \frac{1}{2}}} = 0,
    \end{split}
\end{gather}
moreover, if an exact or high-order integration rule is used, it satisfies conservation of total enstrophy
\begin{gather}
    \begin{split}
        \underbrace{\left(\frac{1}{2} \bilnear{{\bar{\omega}}_{n + 1} + {\omega}'_{n + 1}}{{\bar{\omega}}_{n + 1} + {\omega}'_{n + 1}} \right.}_{\mathcal{E}_{n + 1}} - \underbrace{\left.\frac{1}{2}\bilnear{{\bar{\omega}}_{n} + {\omega}'_{n}}{{\bar{\omega}}_{n} + {\omega}'_{n}} \right)}_{\mathcal{E}_n} + \\\frac{ \Delta t}{Re}\underbrace{\bilnear{\nabla \times ({\bar{\omega}}_{n + \frac{1}{2}} + \omega'_{n + \frac{1}{2}})}{ \nabla \times ({\bar{\omega}}_{n + \frac{1}{2}} + \omega'_{n + \frac{1}{2}})}}_{2\mathscr{P}_{n + \frac{1}{2}}} = 0,
    \end{split}
\end{gather}

The conservation properties of the proposed VMS approach are linked to the orthogonality conditions between the resolved and unresolved scales. We recall the orthogonality conditions imposed by our chosen projector for the Navier-Stokes equations from \eqref{eq:NS_ortho} which holds for any arbitrary time level $n$. 
    
Conservation of mass is trivially satisfied as we explicitly enforce $\nabla \cdot \bar{\underline{u}}$ and $\nabla \cdot \bar{\underline{u}}$ to zero through the means of Lagrange multipliers $\bar{P}$ and $P'$ respectively.

Considering the conservation of total vorticity, we find
\begin{align}
    \mathcal{W}_{n + 1} &= \bilnear{1}{\bar{\omega}_{n + 1} + \omega'_{n + 1}} \\
    \mathcal{W}_{n + 1} &= \bilnear{1}{\bar{\omega}_{n + 1}} + \bilnear{1}{{\omega}'_{n + 1}} = \int_{\partial \Omega} u_{||} \: \mathrm{d} \Gamma_t
    %\mathcal{W}_{n + 1} &= \bilnear{1}{\bar{\omega}_{n + 1}} = 0.
\end{align}
The term $\bilnear{1}{{\omega}'_{n + 1}}$ vanishes as dictated by the first equation of \eqref{eq:NS_ortho} with a particular choice of the test function $\varepsilon^h = 1$. We are then left with $ \bilnear{1}{\bar{\omega}_{n + 1}}$ with equates to the boundary integral as per \eqref{eq:NS_resolved}. 

To show the conservation of energy, we start by noting the following two relations for the kinetic energy $\mathcal{K}$ and enstrophy $\mathcal{E}$
\begin{gather}
    \frac{1}{\Delta t} \bilnear{\underline{u}_{n + \frac{1}{2}}}{\underline{u}_{n + 1} - \underline{u}_{n}} = \frac{1}{2 \Delta t} \left(\bilnear{\underline{u}_{n + 1}}{\underline{u}_{n + 1}} - \bilnear{\underline{u}_{n}}{\underline{u}_{n}}\right) = \frac{1}{\Delta t} (\mathcal{K}_{n + 1} - \mathcal{K}_n) \label{eq:ke_def} \\
    \bilnear{\underline{u}_{n + \frac{1}{2}}}{\nabla \times \omega_{n + \frac{1}{2}}} = \bilnear{\omega_{n + \frac{1}{2}}}{\omega_{n + \frac{1}{2}}} = 2 \mathcal{E}_{n + \frac{1}{2}}. \label{eq:enstro_def} 
\end{gather}
What we want to show is that the total kinetic energy is dissipated with the enstrophy, and that energy is exactly conserved in the inviscid limit. Specifically, we need to show 
\begin{gather}
    \begin{split}
        \frac{1}{\Delta t}((\bar{\mathcal{K}}_{n + 1} + {\mathcal{K}}'_{n + 1}) - (\bar{\mathcal{K}}_n + {\mathcal{K}}'_{n}))  + \frac{2}{Re} (\bar{\mathcal{E}}_{n + \frac{1}{2}} + \mathcal{E}'_{n + \frac{1}{2}}) = \\ \frac{1}{\Delta t} \bilnear{\underline{\bar{u}}_{n + \frac{1}{2}}}{\underline{u}'_{n} - \underline{u}'_{n + 1}} - \frac{1}{Re} \bilnear{\underline{\bar{u}}_{n + \frac{1}{2}}}{\nabla \times \omega'_{n + \frac{1}{2}}} + \\ \frac{1}{\Delta t} \bilnear{\underline{{u}}'_{n + \frac{1}{2}}}{\underline{\bar{u}}_{n} - \underline{\bar{u}}_{n + 1}} - \frac{1}{Re} \bilnear{\underline{{u}}'_{n + \frac{1}{2}}}{\nabla \times \bar{\omega}_{n + \frac{1}{2}}},
    \end{split}
    \label{eq:KE_target}
\end{gather}
which is what we get when filling in $\underline{u} = \bar{\underline{u}} + {\underline{u}}'$ and $\omega = \bar{\omega} + \omega'$ in \eqref{eq:ke_def} and \eqref{eq:enstro_def} where we define resolved and unresolved kinetic energy and enstrophy
\[ \bar{\mathcal{K}} := \frac{1}{2}\bilnear{\underline{\bar{u}}}{\underline{\bar{u}}}, \quad  \mathcal{K}' := \frac{1}{2}\bilnear{\underline{{u}}'}{\underline{{u}}'}, \quad  \bar{\mathcal{E}} := \frac{1}{2}\bilnear{{\bar{\omega}}}{{\bar{\omega}}}, \quad  \mathcal{E}' := \frac{1}{2}\bilnear{{{\omega}'}}{{{\omega}'}} .\]

We consider the resolved and unresolved momentum evolution equations from \eqref{eq:NS_resolved} and \eqref{eq:NS_unresolved} where we take $\underline{v}^h = \underline{\bar{u}}_{n + \frac{1}{2}}$ and $\underline{v}' = \underline{{u}}'_{n + \frac{1}{2}}$ as the test functions. Noting that $\nabla \cdot \bar{\underline{u}} = 0$ and $\nabla \cdot {\underline{u}}' = 0$, we get the following equations for these particular choices of test functions
\begin{gather}
    \begin{split}
        \frac{1}{\Delta t} \bilnear{\underline{\bar{u}}_{n + \frac{1}{2}}}{\bar{\underline{u}}_{n + 1} - \bar{\underline{u}}_{n}} + \frac{1}{Re} \bilnear{\underline{\bar{u}}_{n + \frac{1}{2}}}{\nabla \times \bar{\omega}_{n + \frac{1}{2}}} = \\\frac{1}{\Delta t} \bilnear{\underline{\bar{u}}_{n + \frac{1}{2}}}{\underline{u}'_{n}} - \frac{1}{2 Re} \bilnear{\underline{\bar{u}}_{n + \frac{1}{2}}}{\nabla \times \omega'_n} + a(\tilde{\underline{u}}_{n + \frac{1}{2}}, \tilde{\omega}_{n + \frac{1}{2}}, \underline{\bar{u}}_{n + \frac{1}{2}})
    \end{split} \\
    \begin{split}
        \frac{1}{\Delta t} \bilnear{\underline{{u}}'_{n + \frac{1}{2}}}{{\underline{u}}'_{n + 1} - {\underline{u}}'_{n}} + \frac{1}{Re} \bilnear{\underline{{u}}'_{n + \frac{1}{2}}}{\nabla \times {\omega}'_{n + \frac{1}{2}}} = \\\frac{1}{\Delta t} \bilnear{\underline{{u}}'_{n + \frac{1}{2}}}{\underline{\bar{u}}_{n}} - \frac{1}{2 Re} \bilnear{\underline{{u}}'_{n + \frac{1}{2}}}{\nabla \times \bar{\omega}_n} + a(\tilde{\underline{u}}_{n + \frac{1}{2}}, \tilde{\omega}_{n + \frac{1}{2}}, \underline{{u}}'_{n + \frac{1}{2}})
    \end{split}
\end{gather}
We may simplify the expressions using \eqref{eq:ke_def} and \eqref{eq:enstro_def} along with the definitions of the resolved and unresolved kinetic energy and enstrophy
\begin{gather}
        \frac{1}{\Delta t}(\bar{\mathcal{K}}_{n + 1} - \bar{\mathcal{K}}_n)  + \frac{2}{Re} \bar{\mathcal{E}}_{n + \frac{1}{2}} = \frac{1}{\Delta t} \bilnear{\underline{\bar{u}}_{n + \frac{1}{2}}}{\underline{u}'_{n}} - \frac{1}{2 Re} \bilnear{\underline{\bar{u}}_{n + \frac{1}{2}}}{\nabla \times \omega'_n} + a(\tilde{\underline{u}}_{n + \frac{1}{2}}, \tilde{\omega}_{n + \frac{1}{2}}, \underline{\bar{u}}_{n + \frac{1}{2}})
\end{gather}
\begin{gather}
        \frac{1}{\Delta t} ({\mathcal{K}}'_{n + 1} - {\mathcal{K}}'_n) + \frac{2}{Re} \mathcal{E}'_{n + \frac{1}{2}} = \frac{1}{\Delta t} \bilnear{\underline{{u}}'_{n + \frac{1}{2}}}{\underline{\bar{u}}_{n}} - \frac{1}{2 Re} \bilnear{\underline{{u}}'_{n + \frac{1}{2}}}{\nabla \times \bar{\omega}_n} + a(\tilde{\underline{u}}_{n + \frac{1}{2}}, \tilde{\omega}_{n + \frac{1}{2}}, \underline{{u}}'_{n + \frac{1}{2}}).
\end{gather}
We can add the two equations together, noting the following property of the trilinear form \cite{Zhang2024AConditions}
\begin{gather}
    a(\tilde{\underline{u}}_{n + \frac{1}{2}}, \tilde{\omega}_{n + \frac{1}{2}}, \underline{\bar{u}}_{n + \frac{1}{2}}) + a(\tilde{\underline{u}}_{n + \frac{1}{2}}, \tilde{\omega}_{n + \frac{1}{2}}, \underline{{u}}'_{n + \frac{1}{2}}) = a(\tilde{\underline{u}}_{n + \frac{1}{2}}, \tilde{\omega}_{n + \frac{1}{2}}, \underline{\bar{u}}_{n + \frac{1}{2}} + \underline{{u}}'_{n + \frac{1}{2}}) = a(\tilde{\underline{u}}_{n + \frac{1}{2}}, \tilde{\omega}_{n + \frac{1}{2}}, \tilde{\underline{u}}_{n + \frac{1}{2}}) = 0,
\end{gather}
and we get
\begin{gather}
    \begin{split}
        &\frac{1}{\Delta t}((\bar{\mathcal{K}}_{n + 1} + {\mathcal{K}}'_{n + 1}) - (\bar{\mathcal{K}}_n + {\mathcal{K}}'_{n}))  + \frac{2}{Re} (\bar{\mathcal{E}}_{n + \frac{1}{2}} + \mathcal{E}'_{n + \frac{1}{2}}) = \\ \frac{1}{\Delta t} \bilnear{\underline{\bar{u}}_{n + \frac{1}{2}}}{\underline{u}'_{n}} - &\frac{1}{2 Re} \bilnear{\underline{\bar{u}}_{n + \frac{1}{2}}}{\nabla \times \omega'_n} + \frac{1}{\Delta t} \bilnear{\underline{{u}}'_{n + \frac{1}{2}}}{\underline{\bar{u}}_{n}} - \frac{1}{2 Re} \bilnear{\underline{{u}}'_{n + \frac{1}{2}}}{\nabla \times \bar{\omega}_n}
    \end{split}.
    \label{eq:energy_halfway}
\end{gather}
The left-hand side of the resulting expression is identical to the left-hand side of \eqref{eq:KE_target}. What remains is to show that the right-hand sides also match. To show this, we essentially add two \emph{zeros} to the right-hand side of \eqref{eq:energy_halfway}. The first zero we add is the second expression from \eqref{eq:NS_ortho} with $\underline{v}^h = \underline{\bar{u}}_{n + \frac{1}{2}}$ with the solution at the $n + 1$ time level. The second zero we add is the same expression but with $\underline{v}^h = \underline{\bar{u}}_{n + 1}$ as the test function while taking the solution at the $n + \frac{1}{2}$ time level.
\begin{gather}
    \begin{split}
        \frac{1}{\Delta t}((\bar{\mathcal{K}}_{n + 1} + {\mathcal{K}}'_{n + 1}) - (\bar{\mathcal{K}}_n &+ {\mathcal{K}}'_{n}))  + \frac{2}{Re} (\bar{\mathcal{E}}_{n + \frac{1}{2}} + \mathcal{E}'_{n + \frac{1}{2}}) = \\ \frac{1}{\Delta t} \bilnear{\underline{\bar{u}}_{n + \frac{1}{2}}}{\underline{u}'_{n}} - \frac{1}{2 Re} \bilnear{\underline{\bar{u}}_{n + \frac{1}{2}}}{\nabla \times \omega'_n} + &\frac{1}{\Delta t} \bilnear{\underline{{u}}'_{n + \frac{1}{2}}}{\underline{\bar{u}}_{n}} - \frac{1}{2 Re} \bilnear{\underline{{u}}'_{n + \frac{1}{2}}}{\nabla \times \bar{\omega}_n} \\
        \underbrace{-\frac{1}{\Delta t} \bilnear{\underline{\bar{u}}_{n + \frac{1}{2}}}{\underline{{u}}'_{n + 1}} - \frac{1}{2 Re} \bilnear{\underline{\bar{u}}_{n + \frac{1}{2}}}{\nabla \times {\omega}'_{n + 1}}}_{0} - &\underbrace{\frac{1}{\Delta t} \bilnear{\underline{\bar{u}}_{n + 1}}{\underline{{u}}'_{n + \frac{1}{2}}} - \frac{1}{2 Re} \bilnear{\underline{\bar{u}}_{n + 1}}{\nabla \times {\omega}'_{n + \frac{1}{2}}}}_{0}
    \end{split}
\end{gather}
Grouping things together, we get
\begin{gather}
    \begin{split}
        \frac{1}{\Delta t}((\bar{\mathcal{K}}_{n + 1} + &{\mathcal{K}}'_{n + 1}) - (\bar{\mathcal{K}}_n + {\mathcal{K}}'_{n}))  + \frac{2}{Re} (\bar{\mathcal{E}}_{n + \frac{1}{2}} + \mathcal{E}'_{n + \frac{1}{2}}) = \\ \frac{1}{\Delta t} \bilnear{\underline{\bar{u}}_{n + \frac{1}{2}}}{\underline{u}'_{n} - \underline{u}'_{n + 1}} - &\frac{1}{Re} \bilnear{\underline{\bar{u}}_{n + \frac{1}{2}}}{\nabla \times \omega'_{n + \frac{1}{2}}} + \\ \frac{1}{\Delta t} \bilnear{\underline{{u}}'_{n + \frac{1}{2}}}{\underline{\bar{u}}_{n} - \underline{\bar{u}}_{n + 1}} - &\frac{1}{2 Re} \bilnear{\underline{{u}}'_{n + \frac{1}{2}}}{\nabla \times \bar{\omega}_n} - \frac{1}{2 Re} \bilnear{\underline{\bar{u}}_{n + 1}}{\nabla \times {\omega}'_{n + \frac{1}{2}}}
    \end{split}
\end{gather}
We once again add another \emph{zero}, namely the first expression from \eqref{eq:NS_ortho} with $\varepsilon^h = \bar{\omega}_{n + \frac{1}{2}}$ as the test function
\begin{gather}
    \begin{split}
        & \quad\quad\quad\quad\quad\quad \frac{1}{\Delta t}((\bar{\mathcal{K}}_{n + 1} + {\mathcal{K}}'_{n + 1}) - (\bar{\mathcal{K}}_n + {\mathcal{K}}'_{n}))  + \frac{2}{Re} (\bar{\mathcal{E}}_{n + \frac{1}{2}} + \mathcal{E}'_{n + \frac{1}{2}}) = \\ \frac{1}{\Delta t} &\bilnear{\underline{\bar{u}}_{n + \frac{1}{2}}}{\underline{u}'_{n} - \underline{u}'_{n + 1}} - \frac{1}{Re} \bilnear{\underline{\bar{u}}_{n + \frac{1}{2}}}{\nabla \times \omega'_{n + \frac{1}{2}}} + \\ \frac{1}{\Delta t} &\bilnear{\underline{{u}}'_{n + \frac{1}{2}}}{\underline{\bar{u}}_{n} - \underline{\bar{u}}_{n + 1}} - \frac{1}{2 Re} \bilnear{\underline{{u}}'_{n + \frac{1}{2}}}{\nabla \times \bar{\omega}_n} - \frac{1}{2 Re} \bilnear{\underline{\bar{u}}_{n + 1}}{\nabla \times {\omega}'_{n + \frac{1}{2}}}  \\
        -\frac{1}{2 Re} &\underbrace{\left( \bilnear{\underline{{u}}'_{n + \frac{1}{2}}}{\nabla \times \bar{\omega}_{n + 1}} - \bilnear{\omega'_{n + \frac{1}{2}}}{\bar{\omega}_{n + 1}} \right)}_{0}
    \end{split}.
\end{gather}
Simplifying the expression yields
\begin{gather}
    \begin{split}
        \frac{1}{\Delta t}((\bar{\mathcal{K}}_{n + 1} + {\mathcal{K}}'_{n + 1}) - (\bar{\mathcal{K}}_n + {\mathcal{K}}'_{n}))  + \frac{2}{Re} (\bar{\mathcal{E}}_{n + \frac{1}{2}} + \mathcal{E}'_{n + \frac{1}{2}}) = \\ \frac{1}{\Delta t} \bilnear{\underline{\bar{u}}_{n + \frac{1}{2}}}{\underline{u}'_{n} - \underline{u}'_{n + 1}} - \frac{1}{Re} \bilnear{\underline{\bar{u}}_{n + \frac{1}{2}}}{\nabla \times \omega'_{n + \frac{1}{2}}} + \\ \frac{1}{\Delta t} \bilnear{\underline{{u}}'_{n + \frac{1}{2}}}{\underline{\bar{u}}_{n} - \underline{\bar{u}}_{n + 1}} - \frac{1}{Re} \bilnear{\underline{{u}}'_{n + \frac{1}{2}}}{\nabla \times \bar{\omega}_{n + \frac{1}{2}}} - \\ \cancel{\frac{1}{2 Re} \bilnear{\underline{\bar{u}}_{n + 1}}{\nabla \times {\omega}'_{n + \frac{1}{2}}} - \frac{1}{2 Re} \bilnear{\omega'_{n + \frac{1}{2}}}{\bar{\omega}_{n + 1}}}
    \end{split},
\end{gather}
with the final term dropping out due to the orthogonality between the scales. We finally end up with the same expression from \eqref{eq:KE_target}
\begin{gather}
    \begin{split}
        \frac{1}{\Delta t}((\bar{\mathcal{K}}_{n + 1} + {\mathcal{K}}'_{n + 1}) - (\bar{\mathcal{K}}_n + {\mathcal{K}}'_{n}))  + \frac{2}{Re} (\bar{\mathcal{E}}_{n + \frac{1}{2}} + \mathcal{E}'_{n + \frac{1}{2}}) = \\ \frac{1}{\Delta t} \bilnear{\underline{\bar{u}}_{n + \frac{1}{2}}}{\underline{u}'_{n} - \underline{u}'_{n + 1}} - \frac{1}{Re} \bilnear{\underline{\bar{u}}_{n + \frac{1}{2}}}{\nabla \times \omega'_{n + \frac{1}{2}}} + \\ \frac{1}{\Delta t} \bilnear{\underline{{u}}'_{n + \frac{1}{2}}}{\underline{\bar{u}}_{n} - \underline{\bar{u}}_{n + 1}} - \frac{1}{Re} \bilnear{\underline{{u}}'_{n + \frac{1}{2}}}{\nabla \times \bar{\omega}_{n + \frac{1}{2}}}
    \end{split}
\end{gather}
whereby we have proven conservation of energy for the VMS formulation. The proof for conservation of enstrophy follows the same recipe. In view of the conservation properties induced by the orthogonality conditions of the defined \emph{optimal projector}, this projector may be interpreted as a structure-preserving projector.

Having shown the conservation properties, we now make a brief assessment of the nature of our chosen projector from \eqref{eq:NS_proj}. Since the first equation of \eqref{eq:NS_proj} holds for all admissible $\varepsilon^h \in \bar{U}$, we can make a particular choice of $\varepsilon^h = \varepsilon^h_0 \in \mathrm{Ker}^h(\nabla \times)$ and get 
\[\bilnear{\varepsilon^h_0}{\bar{\omega}} = \bilnear{\varepsilon^h_0}{{\omega}}, \quad \forall \varepsilon^h_0 \in \mathrm{Ker}^h(\nabla \times),\]
which essentially matches the component of $\bar{\omega}$ in $\mathrm{Ker}^h(\nabla \times)$ with the exact $\omega$ we prescribe on the right-hand side. Similarly, we can make a particular choice of $\underline{v}^h = \underline{v}^h_0 \in \mathrm{Ker}^h(\nabla \cdot)$ in the second equation of \eqref{eq:NS_proj} and get
\[ \frac{-1}{2 Re} \bilnear{\underline{v}_0^h}{\nabla \times \bar{\omega}} - \frac{1}{\Delta t} \bilnear{\underline{v}_0^h}{\underline{\bar{u}}} = \frac{-1}{2 Re} \bilnear{\underline{v}_0^h}{\nabla \times {\omega}} - \frac{1}{\Delta t} \bilnear{\underline{v}_0^h}{\underline{{u}}}, \quad \forall \underline{v}_0^h \in \mathrm{Ker}^h(\nabla \cdot), \]
which determines the curl of vorticity in $\mathrm{Ker}^h(\nabla \cdot)$ or equivalently the vorticity in $\mathrm{Ker}^h(\nabla \times)^{\perp}$, and the velocity in $\mathrm{Ker}^h(\nabla \cdot)$. We now remark on a peculiar behaviour of the pressure. Since both $\nabla \times \bar{\omega}$ and $\underline{\bar{u}}$ lie in $\mathrm{Ker}^h(\nabla \cdot)$, the pressure field is entirely invisible to the vorticity and velocity. This means that any term of the form $\bilnear{\nabla \cdot \underline{v}^h}{P}$ may be prescribed on the right-hand side of the second equation in \eqref{eq:NS_proj}, and it would affect $\bar{P}$ without influencing either the velocity or the vorticity. This behaviour arises from the fact that pressure acts solely as a Lagrange multiplier to enforce the divergence-free constraint on the velocity field. The projection of the pressure from \eqref{eq:NS_proj} can be interpreted as matching the \emph{weak} gradient of the pressure with that of the exact pressure. This becomes evident by choosing a specific test function $\underline{v}^h = \underline{v}_{\perp}^h \in \mathrm{Ker}^h(\nabla \cdot)^{\perp}$ in the second equation of \eqref{eq:NS_proj}
\[
    \bilnear{\nabla \cdot \underline{v}_{\perp}^h}{\bar{P}} = \bilnear{\nabla \cdot \underline{v}_{\perp}^h}{{P}}, \quad \forall \, \underline{v}_{\perp}^h \in \mathrm{Ker}^h(\nabla \cdot)^{\perp},
\]
where the terms involving the curl of the vorticity and velocity vanish, since they lie in $\mathrm{Ker}^h(\nabla \cdot)$.

This relation is also reflected in the resolved and unresolved equations of the VMS formulation. By choosing test functions in $\mathrm{Ker}(\nabla \cdot)^{\perp}$ in \eqref{eq:NS_resolved} and \eqref{eq:NS_unresolved}, we obtain
\[
    \bilnear{\nabla \cdot \underline{v}_{\perp}^h}{\bar{P}_{n + \frac{1}{2}}} = a(\underline{v}_{\perp}^h, \underline{u}_{n + \frac{1}{2}}, \omega_{n + \frac{1}{2}}) + \int_{\partial \Omega} \underline{v}_{\perp}^h \hat{P} \, \mathrm{d}\Gamma_P, \quad \forall \, \underline{v}_{\perp}^h \in \mathrm{Ker}^h(\nabla \cdot)^{\perp},
\]
\[
    \bilnear{\nabla \cdot \underline{v}_{\perp}'}{{P}'_{n + \frac{1}{2}}} = a(\underline{v}_{\perp}', \underline{u}_{n + \frac{1}{2}}, \omega_{n + \frac{1}{2}}) + \int_{\partial \Omega} \underline{v}_{\perp}' \hat{P} \, \mathrm{d}\Gamma_P, \quad \forall \, \underline{v}_{\perp}' \in \mathrm{Ker}'(\nabla \cdot)^{\perp}.
\]
Here, the trilinear form $a(\underline{v}, \underline{u}_{n + \frac{1}{2}}, \omega_{n + \frac{1}{2}})$ represents an $L^2$ inner product between the test functions and the non-linear convective term with the full velocity and vorticity fields ($\bar{\underline{u}} + \underline{u}', \: \bar{\omega} + \omega'$). As such, these expressions can be viewed as $L^2$ projections of the convective term into the resolved space $\bar{V}$ and the unresolved space $V'$, respectively. Thus, the weak pressure gradients in both the resolved and unresolved equations are matched to the $L^2$ projections of the convective term. In the numerical examples considered in this work, we observed that the individual resolved and unresolved pressures exhibited nonphysical behaviour, likely due to the $L^2$ projections of the convective term being too weak to accurately represent the full convective dynamics. However, the full pressure field, given by the sum $\bar{P} + P'$, showed no such artefacts and yielded physically consistent results.

\begin{remark}
    We deliberately choose to project the pressure separately using the $L^2$ projector, since the pressure is naturally defined in $L^2(\Omega)$. Therefore, we disregard the value of $\bar{P}$ obtained from \eqref{eq:NS_proj} and instead compute it via an $L^2$ projection of the analytical pressure field. While the VMS formulation is solved using \eqref{eq:NS_resolved} and \eqref{eq:NS_unresolved} without modification, we post-process the computed pressure field by defining $\widetilde{P} := \bar{P} + P'$ and projecting it onto the coarse mesh via the $L^2$ projector: $\bar{P} := \mathcal{P}_{L^2} \widetilde{P}$ and $P' := (1 - \mathcal{P}_{L^2}) \widetilde{P}$.
\end{remark}
The aforementioned strange behaviour of the pressure warrants further investigation, possibly using stronger projectors for the convective term. Furthermore, there remains potential for future work to explore alternative treatments of the pressure term. One possible direction is to employ the pressure Poisson equation to reconstruct the pressure field. However, in this study, we restrict our attention to the simpler and more direct approach described above, where the pressure is projected separately using the $L^2$ projector.

\subsection{Solution algorithm and cost estimates}\label{sec:VMS_cost}
While the resolved-scale system in~\eqref{eq:NS_resolved} is finite-dimensional, it depends continuously on the fine-scale solution, which is governed by the infinite-dimensional system in~\eqref{eq:NS_unresolved}. We therefore require an effective and efficient strategy for solving these non-linear, coupled systems. For notational convenience, we write the two systems in a compact form as follows
\begin{equation}
    \begin{array}{ccccc}
        \bar{\mathcal{L}}\,\bar{\phi}_{n+1} & + & \bar{\mathscr{C}}(\bar{\phi}_{n+1}, \phi'_{n+1}) &=& \bar{f}_{n}, \\[0.2cm]
        \mathscr{C}'(\bar{\phi}_{n+1}, \phi'_{n+1}) & + & \mathcal{L}'\,\phi'_{n+1} &=& f'_{n},
    \end{array}
    \label{eq:sys}
\end{equation}
where $\phi_{n+1}$ denotes the solution vector $(\omega_{n+1}, \underline{u}_{n+1}, P_{n+1})$ at the new time level, $\mathcal{L}$ is the linear Stokes-like operator, $\mathscr{C}$ the non-linear convective term, and $f_n$ collects all quantities known from the previous time step.

To determine the fine scales, we employ the Fine-Scale Greens' operator, following~\cite{Shrestha2024ConstructionScales, Hughes2007VariationalMethods}. As shown in~\cite{Shrestha2025OptimalProblems}, it is sufficient to construct the Greens' operator associated with the \emph{symmetric} part of the full operator, namely the Stokes-like substructure $\mathcal{L}'$ appearing in the VMS formulation of the Navier-Stokes equations. This construction yields $(\mathcal{L}')^{-1} = \mathcal{G}'$, where the Fine-Scale Greens' operator is defined by
\[
\mathcal{G}' = \mathcal{G} - \mathcal{G}\mathcal{P}^{T}\left(\mathcal{P}\mathcal{G}\mathcal{P}^{T}\right)^{-1}\mathcal{P}\mathcal{G},
\]
with $\mathcal{G}$ the classical Greens' operator and $\mathcal{P}$ the projector onto the resolved space. The second term subtracts the resolved component of~$\mathcal{G}$, ensuring that $\mathcal{G}'$ lies entirely in the null space of the projector. In our approach, as highlighted in \Cref{sec:VMS_prop}, the classical Greens' operator is approximated on a finer mesh. More precisely, the resolved scales are discretised in a polynomial space of degree~$p$, while the classical Greens' operator associated with $\mathcal{L}$ is approximated on a richer polynomial space of degree~$p+k$. With this formulation, the coupled system may be written equivalently as
\begin{equation}
    \begin{aligned}
        \phi'_{n+1} = \mathcal{G}'\!\left(f'_{n} - \mathscr{C}'(\bar{\phi}_{n+1}, \phi'_{n+1})\right), \\
        \bar{\mathcal{L}}\,\bar{\phi}_{n+1} + \bar{\mathscr{C}}(\bar{\phi}_{n+1}, \phi'_{n+1}) = \bar{f}_{n}.
    \end{aligned}
    \label{eq:VMS_sys}
\end{equation}

We solve these equations using the iterative procedure outlined in Algorithm~\ref{alg:VMS}.
\begin{algorithm}[htp] \caption{Variational Multiscale Algorithm} \label{alg:VMS} 
\begin{algorithmic}[1] 
\State Input geometry, number of elements $N \times M$, coarse polynomial degree $p$, fine polynomial degree $p + k$, number of time steps $N_{dt}$, time step size $\Delta t$ 
\State \textbf{Assemble matrix} $\mathcal{\bar{L}}$ 
\State \textbf{Assemble matrix} $\mathcal{{L}}'$ \State \textbf{Solve} for Fine-Scale Greens' operator $\mathcal{G}'$ 
\State Project initial conditions $\bar{\phi}_0, \phi'_0$ onto respective meshes \For{$n = 0$ to $N_{dt} - 1$} 
\State Set iteration count \textit{i} = 0 
\State Initialise: $\bar{\phi}_{{n+1}, i} \gets \bar{\phi}_{n}$, \quad $\phi'_{{n+1},i} \gets \phi'_{n}$ 
\State \textbf{Assemble vector} $\mathcal{G}'f'_{n}$ 
\State \textbf{Assemble vector} $\bar{f}_{n}$ 
\State Set \textit{residual} = 1 
\While{\textit{residual} $\geq$ tol} 
\State \textbf{Assemble vector} $\mathcal{G}'\mathscr{C}'(\bar{\phi}_{{n+1}, i}, \phi'_{{n+1}, i})$ 
% \State \textbf{Assemble} $\bar{\mathscr{C}}(\bar{u}_{{n+1}_i}, u'_{{n+1},i})$ 
\State Update $\phi'_{n + 1, i + 1} \gets \mathcal{G}'f'_{n} - \mathcal{G}'\mathscr{C}'(\bar{\phi}_{{n+1}, i}, \phi'_{{n+1}, i})$ \State \textbf{Assemble matrix} $\bar{\mathscr{C}}(\bar{\phi}_{{n+1}, i + 1}, \phi'_{{n+1}, i + 1})$ \State \textbf{Solve} for $\bar{\phi}_{{n+1}, i + 1}$ \[ \bar{\mathcal{L}} \bar{\phi}_{n + 1, i + 1} + \bar{\mathscr{C}}\big(\bar{\phi}_{n + 1, i + 1}, \phi'_{n + 1, i + 1}\big) = \bar{f}_{n} \] \State \textit{residual} $\gets \lVert \bar{\phi}_{{n+1}, i + 1} - \bar{\phi}_{{n+1}, i} \rVert$ \State \textit{i} $\gets$ \textit{i} + 1 \EndWhile 
% \State Update $\bar{u}_{t_{n+1}}, u'_{t_{n+1}}$ 
\EndFor \end{algorithmic}
\end{algorithm}
Since the Fine-Scale Greens' operator depends only on the geometry and the symmetric part of the operator and not on the solution, it can be precomputed as a \emph{one-time offline step} for any given configuration. By approximating only the symmetric Greens' operator in a rich but finite-dimensional subspace, we circumvent the challenges posed by the infinite-dimensional nature of the unresolved-scale problem.

The non-linear system in \eqref{eq:VMS_sys} is solved at each time step using Picard iterations with an $L^2$ tolerance of $10^{-12}$. As shown in Algorithm~\ref{alg:VMS}, each iteration updates the fine scales using the Greens' operator and then solves the resolved-scale problem with this updated estimate. This approach retains non-linear multiscale interactions while avoiding the need to solve the full non-linear Navier--Stokes system on the fine mesh. The number of iterations required depends on the time-step size $\Delta t$, but in all numerical experiments, convergence was typically achieved in $\mathcal{O}(10^{1})$ iterations.

If the fine-scale contribution is omitted entirely (i.e., $\phi'_{n+1} = 0$), the method reduces to the standard Galerkin scheme on the coarse mesh, shown in Algorithm~\ref{alg:Glk}. However, this simplification fails to capture unresolved-scale effects, which are particularly important in convection-dominated regimes, and it does not recover the desired projection properties of the VMS solution.
\begin{algorithm}[htp]
\caption{Galerkin Algorithm} 
\label{alg:Glk} 
\begin{algorithmic}[1]
\State Input geometry, number of elements $N \times M$, coarse polynomial degree $p$, number of time steps $N_{dt}$, time step size $\Delta t$ 
\State \textbf{Assemble matrix} $\mathcal{\bar{L}}$ 
\State Project initial conditions $\bar{\phi}_0$ onto mesh
\For{$n = 0$ to $N_{dt} - 1$} 
\State Set iteration count \textit{i} = 0 
\State Initialise: $\bar{\phi}_{{n+1}, i} \gets \bar{\phi}_{n}$ 
\State \textbf{Assemble vector} $\bar{f}_{n}$ 
\State Set \textit{residual} = 1 \While{\textit{residual} $\geq$ tol} \State \textbf{Assemble matrix} $\bar{\mathscr{C}}(\bar{\phi}_{{n+1}, i + 1})$ 
\State \textbf{Solve} for $\bar{\phi}_{{n+1}, i + 1}$ \[ \bar{\mathcal{L}} \bar{\phi}_{n + 1, i + 1} + \bar{\mathscr{C}}\big(\bar{\phi}_{n + 1, i + 1}\big) = \bar{f}_{n} \] 
\State \textit{residual} $\gets \lVert \bar{\phi}_{{n+1}, i + 1} - \bar{\phi}_{{n+1}, i} \rVert$ 
\State \textit{i} $\gets$ \textit{i} + 1 \EndWhile
% \State Update $\bar{u}_{t_{n+1}}, u'_{t_{n+1}}$ 
\EndFor 
\end{algorithmic} 
\end{algorithm}
With the solution algorithm in place, we now assess the computational cost of the proposed methodology. Since the method requires computing the Fine-Scale Greens' operator on a refined space of polynomial degree $p+k$, the dominant cost arises from operations on this fine mesh. Consequently, the only meaningful baseline for comparison is the standard Galerkin discretisation performed directly on the $p+k$ space. The rationale is straightforward: if the VMS approach is more expensive than simply running a Galerkin method on the finer space, there is little justification for using VMS. In what follows, we present a cost analysis that identifies precisely when the proposed VMS methodology becomes computationally advantageous—namely, when it is cheaper than Galerkin on the $p+k$ space.

To compare the computational cost of the VMS method on a coarse mesh with that of the Galerkin method on a fine mesh, we adopt a simple but expressive cost model. The total cost of an algorithm is decomposed into the cost of (i) assembling all required vectors $\mathcal{A}^{vec}$ and matrices $\mathcal{A}^{mat}$, and (ii) solving the resulting linear systems $\mathcal{S}$.

The total cost of the proposed VMS methodology following \Cref{alg:VMS} is then
\begin{equation}
    C_{\mathrm{VMS}} 
    = 
    \underbrace{\mathcal{A}_{c}^{mat} + \mathcal{A}_{f}^{mat} + \mathcal{S}_c + \mathcal{S}_f}_{\text{Offline computation of }\mathcal{G}'} 
    + 
    N_{dt} (
        \underbrace{\mathcal{A}_c^{vec}}_{\bar{f}_n}
        +
        \underbrace{\mathcal{A}_f^{vec}}_{\mathcal{G}' f'_n}
    )
    +
    N_{dt} i_{ave} (
        \underbrace{\mathcal{A}_f^{vec}}_{\mathcal{G}' \mathscr{C}'}
        +
        \underbrace{\mathcal{A}_c^{mat} + \mathcal{S}_c}_{\bar{\mathcal{L}} + \bar{\mathscr{C}}}
    ),
\end{equation}
with the subscript $c$ and $f$ denoting coarse and fine spaces respectively, and $i_{ave}$ denotes the average number of non-linear iterations performed per time step. For comparison, the cost of a base Galerkin method following \Cref{alg:Glk} on the fine mesh is
\begin{equation}
    C_{\mathrm{Galerkin}} 
    =  
    \mathcal{A}_{f}^{mat} 
    + 
    N_{dt}\mathcal{A}_{f}^{vec} 
    + 
    N_{dt} i_{ave} (\mathcal{A}_f^{mat} + \mathcal{S}_f).
\end{equation}

We now define the cost ratio
\[
    R = \frac{C_{\mathrm{VMS}}}{C_{\mathrm{Galerkin}}},
\]
which expands to
\begin{align}
    R &=\frac{\mathcal{A}_{c}^{mat} + \mathcal{A}_{f}^{mat} + \mathcal{S}_c + \mathcal{S}_f + N_{dt} (\mathcal{A}_c^{vec} + \mathcal{A}_f^{vec}) + N_{dt} i_{ave} (\mathcal{A}_f^{vec} + \mathcal{A}_c^{mat} + \mathcal{S}_c)}{\mathcal{A}_{f}^{mat} + N_{dt}\mathcal{A}_{f}^{vec} + N_{dt} i_{ave} (\mathcal{A}_f^{mat} + \mathcal{S}_f)} \nonumber\\
    &= \frac{(\mathcal{A}_{f}^{mat} + N_{dt}\mathcal{A}_{f}^{vec} + (1 + (-1))N_{dt} i_{ave} (\mathcal{A}_f^{mat} + \mathcal{S}_f)) + (\mathcal{A}_c^{mat} + \mathcal{S}_c + \mathcal{S}_f + N_{dt} \mathcal{A}_c^{vec} + N_{dt} i_{ave} (\mathcal{A}_f^{vec} + \mathcal{A}_c^{mat} + \mathcal{S}_c))}{\mathcal{A}_{f}^{mat} + N_{dt}\mathcal{A}_{f}^{vec} + N_{dt} i_{ave} (\mathcal{A}_f^{mat} + \mathcal{S}_f)} \nonumber\\
    &= 1 + \frac{\mathcal{A}_c^{mat} + \mathcal{S}_c + \mathcal{S}_f + N_{dt} \mathcal{A}_c^{vec} + N_{dt} i_{ave} (\mathcal{A}_f^{vec} + \mathcal{A}_c^{mat} + \mathcal{S}_c) - N_{dt} i_{ave} (\mathcal{A}_f^{mat} + \mathcal{S}_f)}{\mathcal{A}_{f}^{mat} + N_{dt}\mathcal{A}_{f}^{vec} + N_{dt} i_{ave} (\mathcal{A}_f^{mat} + \mathcal{S}_f)}
\end{align}

For the VMS methodology to be computationally favourable, we require $R < 1$. This leads to the condition
\begin{equation}
    (\mathcal{A}_f^{mat} + \mathcal{S}_f) 
    >
    \frac{\mathcal{A}_c^{mat} + \mathcal{S}_c + \mathcal{S}_f + N_{dt} \mathcal{A}_c^{vec}}
    {N_{dt} i_{ave}}
    +
    (\mathcal{A}_f^{vec} + \mathcal{A}_c^{mat} + \mathcal{S}_c).
\end{equation}

This criterion indicates that the VMS method is advantageous when \emph{the cost of matrix assembly and linear solves on the fine mesh dominates} the combined cost of fine-scale vector assemblies, coarse-scale vector and matrix assemblies, and coarse-scale solves. In the regime of large numbers of time steps, the inequality becomes even more favourable, where the dominant terms arise from repeated coarse-scale assemblies and coarse solves. In such cases, the proposed VMS methodology provides a computational advantage by avoiding repeated fine-mesh matrix assemblies and linear solves, while still retaining the essential non-linear multiscale coupling.

\section{Numerical tests} % { \texorpdfstring{$Re = \infty$}{TEXT}}
\label{sec:num_test}
In this section, we present numerical experiments to verify the proposed framework. Two test cases are considered for the illustration of the framework. The first is the two-dimensional Taylor–Green vortex problem, which admits a smooth analytical solution and is used here to verify the order of convergence of the method. The second test case is the inviscid vortex roll-up problem \cite{Sanderse2013Energy-conservingEquations}, which is a more challenging benchmark that lacks viscosity-induced regularisation and serves to demonstrate the method’s performance in the limiting case. Together, these examples illustrate both the accuracy and robustness of the proposed approach. As noted in \Cref{rmk:under_int}, the choice between over- and under-integration affects the scheme's conservation properties. Since enstrophy conservation is not the primary concern in the convergence tests with the Taylor-Green vortex and under-integration represents the more critical scenario, we perform these tests using under-integration. For the second test case, however, enstrophy conservation is important, and we therefore employ over-integration. 
% For completeness, the under-integration results for the vortex roll-up problem are also presented in \ref{app:enstrophy}.

\subsection{Taylor-Green vortex}
\begin{figure}[htp]
    \centering
    \includegraphics[width=0.35\linewidth]{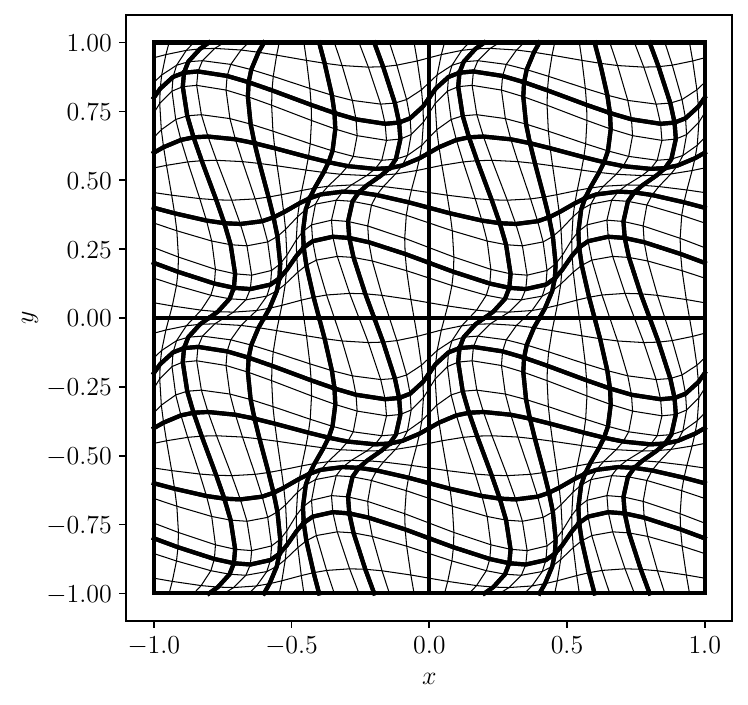}
    \caption{Example of curvilinear mesh with $10 \times 10$ elements with polynomial degree $p = 3$}
    \label{fig:msh_curv}
\end{figure}

We start with the 2D Taylor–Green vortex benchmark problem featuring the following exact solutions
\begin{align}
    u_x(x, y, t) &= -\mathrm{sin}(\pi x) \mathrm{cos}(\pi y) e^{-\frac{2 \pi^2 t}{Re}} \\
    u_y(x, y, t) &= \mathrm{cos}(\pi x) \mathrm{sin}(\pi y) e^{-\frac{2 \pi^2 t}{Re}} \\
    \pbar(x, y, t) &= \frac{1}{4} (\mathrm{cos}(2 \pi x) + \mathrm{cos}(2 \pi y)) e^{-\frac{4 \pi^2 t}{Re}} \\
    \omega(x, y, t) &= -2 \pi \mathrm{sin}(\pi x) \mathrm{sin}(\pi y) e^{-\frac{2 \pi^2 t}{Re}}.
\end{align}
We carry out all the computations on curvilinear grids over a periodic domain $(x, y) \in ]-1, 1[^2$ with the mapping ${\Phi} \: : \: (\xi, \eta) \rightarrow (x, y)$ from the reference domain $(\xi, \eta) \in [-1, 1]^2$ to the physical domain $(x, y)$ partitioned into $N \times N$ curvilinear elements given by
\begin{equation}
\left\{\begin{array} { l } 
{ \Phi_x = \hat { x } + 0.1 \sin{( 2 \pi \hat { x } )} \sin {( 2 \pi \hat { y } )}} \\
{ \Phi_y = \hat { y } - 0.1 \sin{( 2 \pi \hat { x } )} \sin {( 2 \pi \hat { y } )}}
\end{array} \text {, where } \left\{\begin{array}{l}
\hat{x} = -1 + \frac{2}{N}( i + \frac{1}{2}(1+\xi)), \quad i = 0, 1, \hdots N - 1 \\
\hat{y} = -1 + \frac{2}{N}( j + \frac{1}{2}(1+\eta)), \quad j = 0, 1, \hdots N - 1
\end{array}\right.\right. \text {. }
\end{equation}
An example from this family of curvilinear meshes is depicted in \Cref{fig:msh_curv}.
We perform a rate-of-convergence study by measuring the errors in the computed vorticity, velocity, and pressure fields using appropriate error norms, namely vorticity error in the ${H}_0(\mathrm{curl}, \Omega)$ norm, velocity error in the $H(\mathrm{div, \Omega})$ norm, and static pressure error in the $L^2(\Omega)$ norm
\begin{align}
    vorticity \: error \: &: \: \lVert \nabla \times \bar{\omega} - \nabla \times \omega \rVert_{L^2} \\
    velocity \: error \: &: \: \left( \lVert \underline{\bar{u}} - \underline{u} \rVert_{L^2}^2 + \lVert \nabla \cdot \bar{\underline{u}} - \nabla \cdot \underline{u} \rVert_{L^2}^2 \right)^{\frac{1}{2}} \\
    static \:pressure \: error \: &: \: \lVert {\bar{\pbar}} - {\pbar} \rVert_{L^2}.
\end{align}
In this test, we compare the convergence behaviour of three different solutions: the base Galerkin scheme, the exact projection of the exact solution, and the proposed VMS approach. For the VMS method, we vary the refinement parameter $k$, which defines the polynomial degree $p+k$ of the fine-scale approximation space used in the construction of the Fine-Scale Greens' function. This allows us to assess how the quality of the fine-scale approximation affects the accuracy of the VMS solution. Additionally, we evaluate the error of the VMS solution with respect to the exact projection ($\mathcal{P}\omega, \mathcal{P}\underline{u}, \mathcal{P}\pbar$) to determine how well the method recovers the desired projected solution. We also examine the convergence of the approximated unresolved scales ($\omega^{\prime}_k, \underline{u}^{\prime}_k, \pbar^{\prime}_k$) towards the exact unresolved scales ($\omega^{\prime}, \underline{u}^{\prime}, \pbar^{\prime}$). All errors are computed using the same set of norms introduced previously. The tests are conducted at a Reynolds number of $Re = 100$, using a time step of $\Delta t = 0.04$ to march the solution forward to a final time of $t = 1$, at which point the error calculations are performed. To accurately capture the exact solution, we evaluate the errors using high-order Gauss–Lobatto quadrature of order 25.

\begin{figure}[H]
    \begin{subfigure}{0.33\linewidth}
        \centering
        \includegraphics[width=\linewidth]{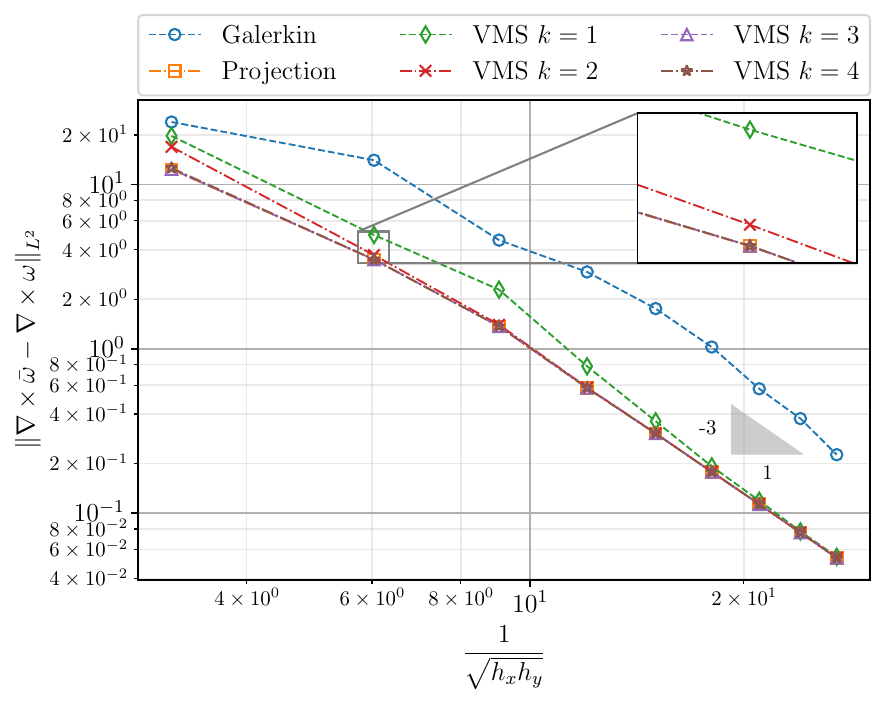}
        \caption{Vorticity error}
        \label{fig:hp_vorticity}
    \end{subfigure}
    \begin{subfigure}{0.33\linewidth}
        \centering
        \includegraphics[width=\linewidth]{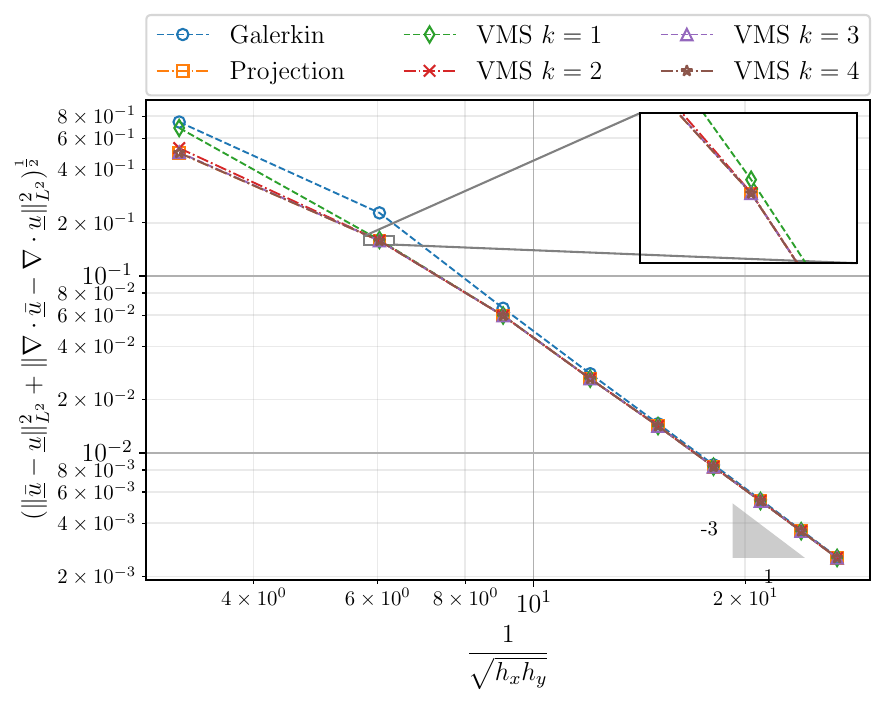}
        \caption{Velocity error}
        \label{fig:hp_u}
    \end{subfigure}
    \begin{subfigure}{0.33\linewidth}
        \centering
        \includegraphics[width=\linewidth]{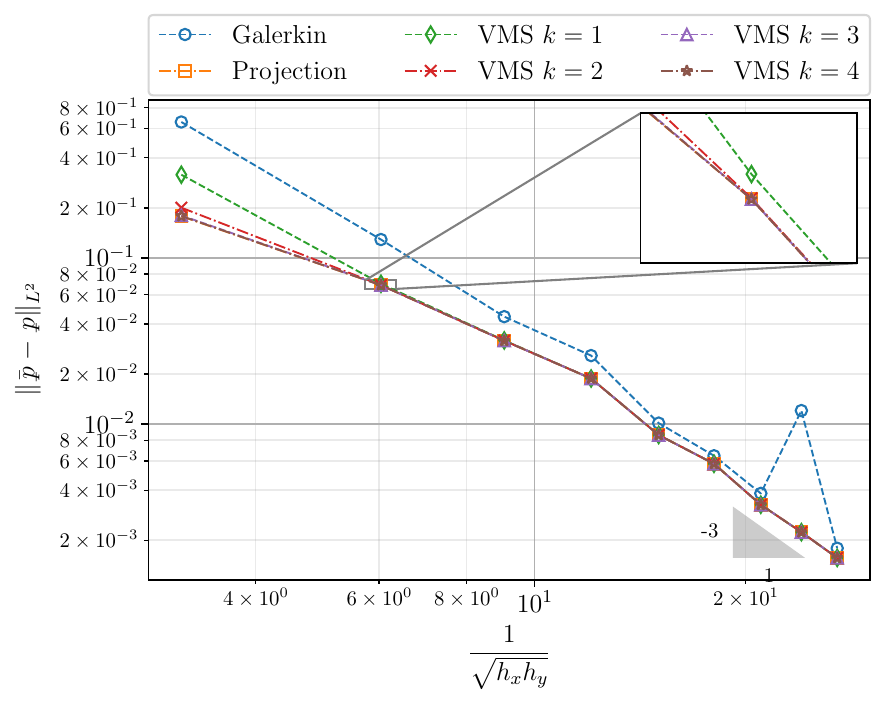}
        \caption{Static pressure error}
        \label{fig:hp_p}
    \end{subfigure}
    \caption{$h$-convergence of the different solutions for the Taylor-Green vortex problem at $Re = 100$ and $\Delta t = 0.04$ on curvilinear meshes with polynomial degree $p = 3$}
    \label{fig:hp_converg}
\end{figure}

\begin{figure}[htp]
    \begin{subfigure}{0.33\linewidth}
        \centering
        \includegraphics[width=\linewidth]{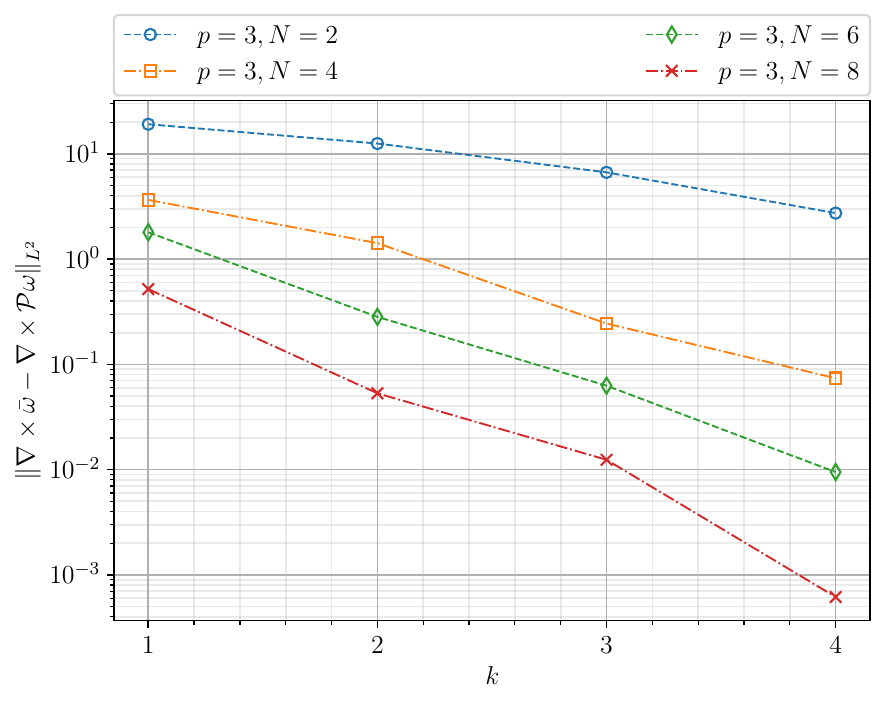}
        \caption{Vorticity error}
        \label{fig:k_vorticity_bar}
    \end{subfigure}
    \begin{subfigure}{0.33\linewidth}
        \centering
        \includegraphics[width=\linewidth]{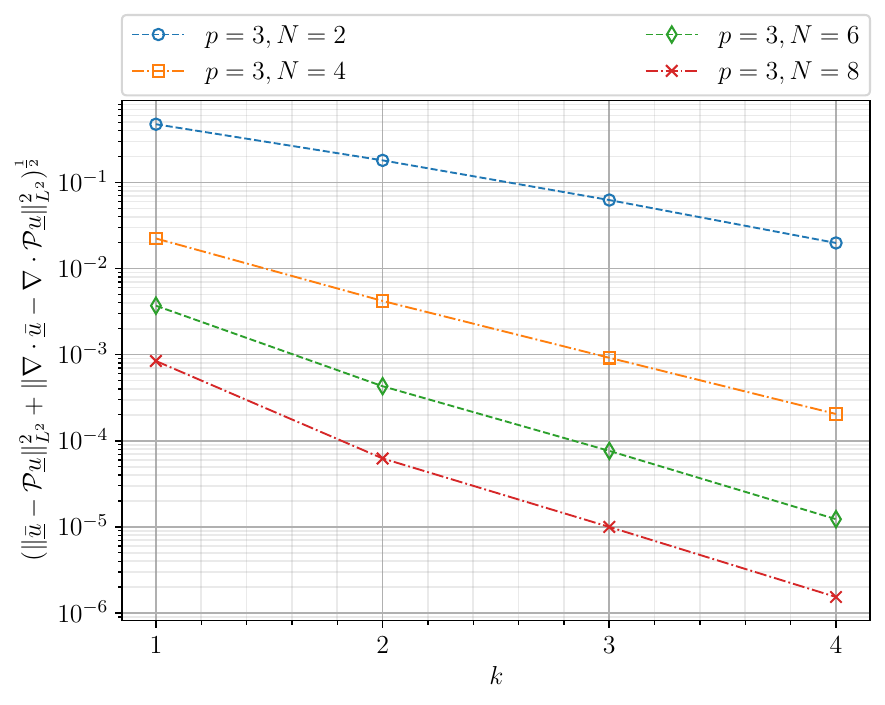}
        \caption{Velocity error}
        \label{fig:k_u_bar}
    \end{subfigure}
    \begin{subfigure}{0.33\linewidth}
        \centering
        \includegraphics[width=\linewidth]{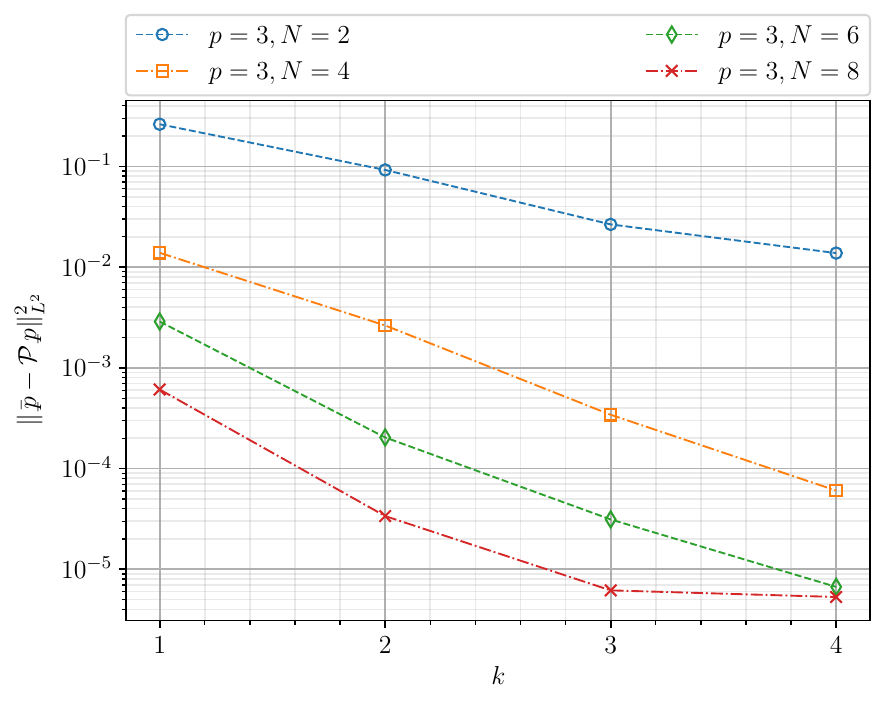}
        \caption{Static pressure error}
        \label{fig:k_p_bar}
    \end{subfigure}
    \caption{$k$-convergence of the VMS solutions with respect to the exact projection for the Taylor-Green vortex problem at $Re = 100$ and $\Delta t = 0.04$ on curvilinear meshes with polynomial degree $p = 3$}
    \label{fig:k_bar}
\end{figure}

The convergence plots for the different solutions are shown in \Cref{fig:hp_converg}, where the errors are plotted against the mesh sizes $h_x$ and $h_y$, which are analogous to the number of elements $N$. For the sake of brevity, we only present the convergence demonstration on meshes with polynomial degree $p = 3$. All three solution strategies, namely the base Galerkin scheme, the exact projection of the exact solution, and the proposed VMS approach, exhibit the optimal convergence rate of order three, as expected for the chosen discretisation (see \Cref{rem:degree}). Among the three solution strategies, the Galerkin method consistently yields the largest errors across all variables, while the exact projection, by construction, attains the lowest error. The VMS solution falls between these two, with its accuracy improving as the fine-scale refinement parameter $k$ increases. This trend demonstrates that the VMS method increasingly approximates the exact projection as the resolution of the unresolved scales is enhanced.
Building further on that point, we present the $k$-convergence of the VMS solution with respect to the exact projection in \Cref{fig:k_bar}. 
These plots demonstrate how rapidly the VMS solution converges to the target projection of the exact solution with increasing $k$. We generally observe a linear trend in the semi-log scale plots in \Cref{fig:k_bar}, which indicates exponential convergence. The curves are not exactly linear as seen for the last curve in \Cref{fig:k_p_bar}, for instance. This local ``roughness" can be attributed to the curvilinear mesh, which manifests large local errors, particularly in the pressure field as observed through \Cref{fig:hp_p} and \Cref{fig:k_p_bar}. 

\begin{figure}[H]
    \begin{subfigure}{0.33\linewidth}
        \centering
        \includegraphics[width=\linewidth]{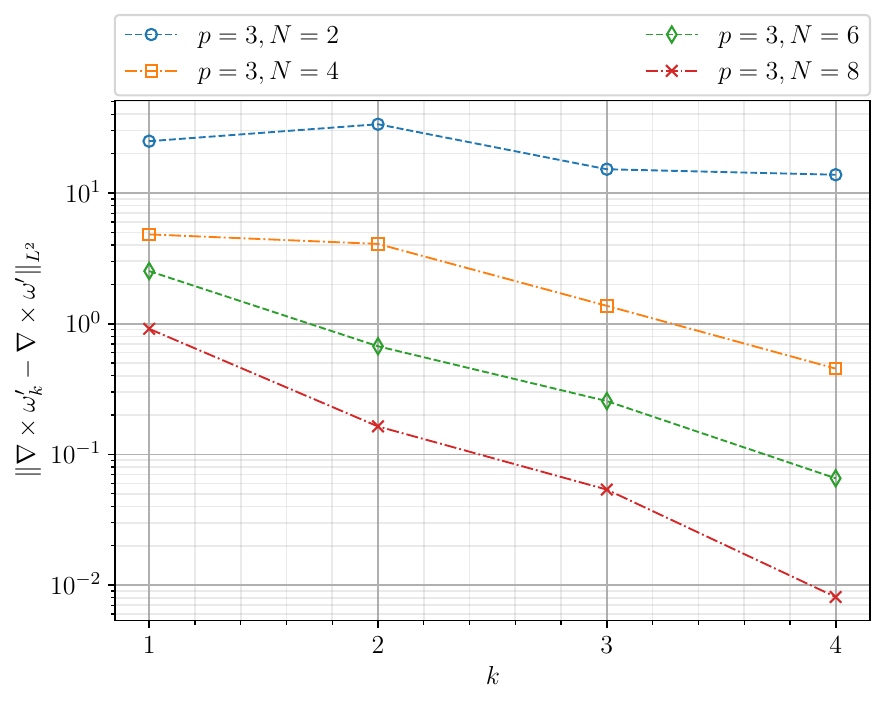}
        \caption{Vorticity error}
        \label{fig:k_vorticity_prime}
    \end{subfigure}
    \begin{subfigure}{0.33\linewidth}
        \centering
        \includegraphics[width=\linewidth]{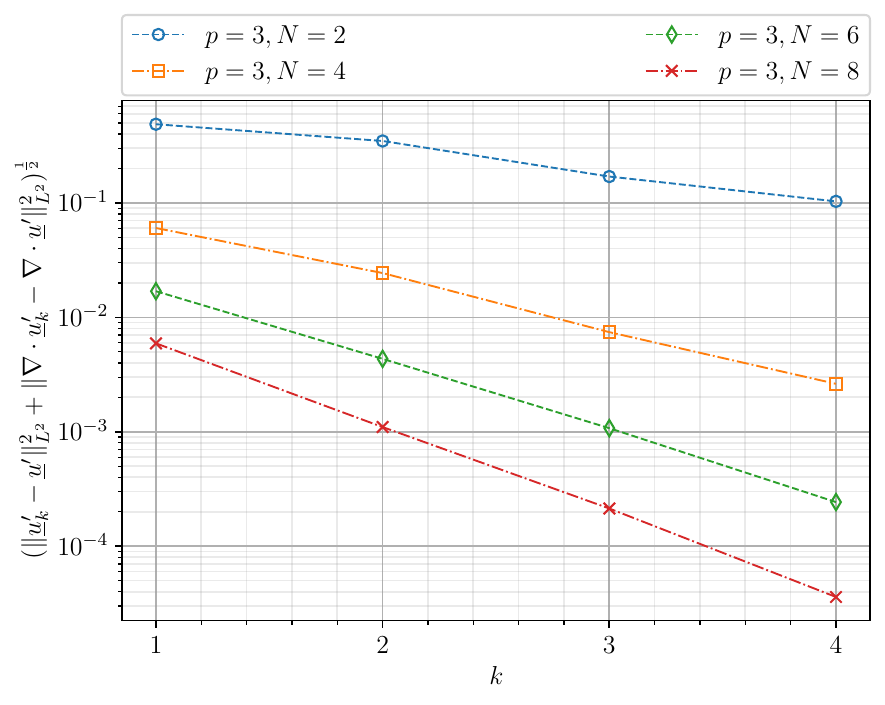}
        \caption{Velocity error}
        \label{fig:k_u_prime}
    \end{subfigure}
    \begin{subfigure}{0.33\linewidth}
        \centering
        \includegraphics[width=\linewidth]{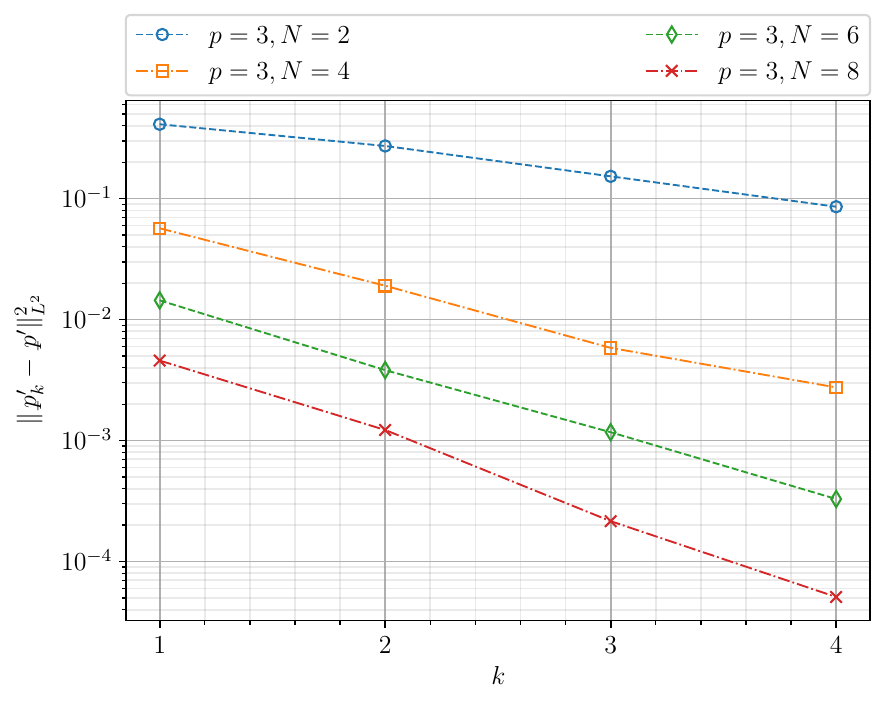}
        \caption{Static pressure error}
        \label{fig:k_p_prime}
    \end{subfigure}
    \caption{$k$-convergence of the approximate unresolved scales with respect to the exact unresolved scales for the Taylor-Green vortex problem at $Re = 100$ and $\Delta t = 0.04$ on curvilinear meshes with polynomial degree $p = 3$}
    \label{fig:k_prime}
\end{figure}

Lastly, \Cref{fig:k_prime} depicts the $k$-convergence of the computed unresolved scales with respect to the exact unresolved scales. We once again observe exponential convergence with increasing $k$ as the error curves all have a linear trend in the semi-log plot. Furthermore, we may note that the absolute error values for the unresolved scales shown in \Cref{fig:k_prime} are consistently larger than the error between the VMS solution and the exact projection. 
% We thus state that a highly accurate estimation of the fine scales is not essential for obtaining a good approximation of the projection through the proposed VMS approach.
\subsection{Vortex roll-up}
We now consider the inviscid vortex roll-up problem ($Re = \infty$) as a benchmark test case to illustrate the performance of the proposed method. This problem is known to be particularly challenging due to the development of sharp gradients and intricate small-scale structures, which must be captured without the aid of viscous regularisation. The problem is considered on a periodic domain $\Omega = ]0, 2\pi[^2$ with the following initial condition for the two velocity components
\begin{equation}
    u_x(x, y)=\left\{\begin{array}{ll}
    \tanh \left(\displaystyle \frac{y-\frac{\pi}{2}}{\delta}\right), & y \leq \pi, \\
    \tanh \left(\displaystyle \frac{\frac{3 \pi}{2}-y}{\delta}\right), & y>\pi,
    \end{array} \quad \quad \quad  u_y(x, y)=\epsilon \sin (x)\right.,
\end{equation}
where $\delta = \frac{\pi}{15}$ and $\epsilon = 0.05$.

\begin{figure}[H]
\centering
\begin{minipage}{0.35\linewidth}
    \centering
    \includegraphics[width=\linewidth]{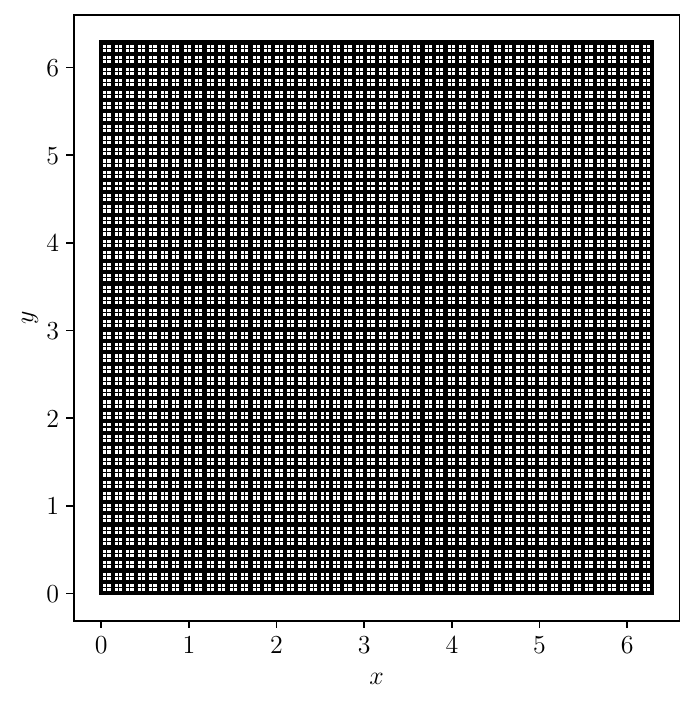}
    \caption{Highly refined mesh with $48 \times 48$ elements with polynomial degree 2 for reference}
    \label{fig:msh_fine}
\end{minipage}
\begin{minipage}{0.35\linewidth}
    \centering
    \includegraphics[width=\linewidth]{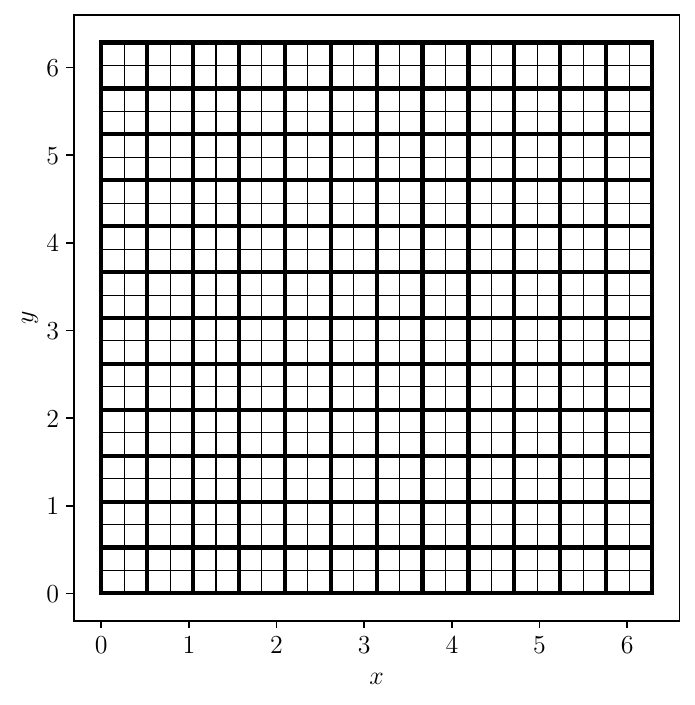}
    \caption{Coarse mesh with $12 \times 12$ elements with polynomial degree 2}
    \label{fig:msh_coarse}
\end{minipage}
\end{figure}

\begin{figure}[H]
    \centering
    \begin{subfigure}{0.32\textwidth}
        \includegraphics[width=\textwidth]{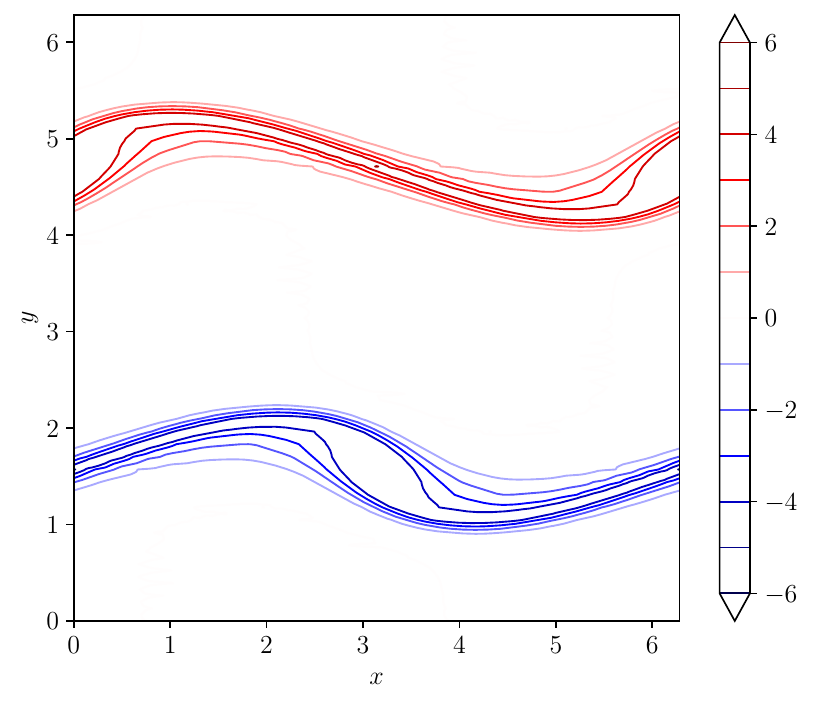}
        \caption{$t = 4$}
        % \label{fig:sub1}
    \end{subfigure}
    \begin{subfigure}{0.32\textwidth}
        \includegraphics[width=\textwidth]{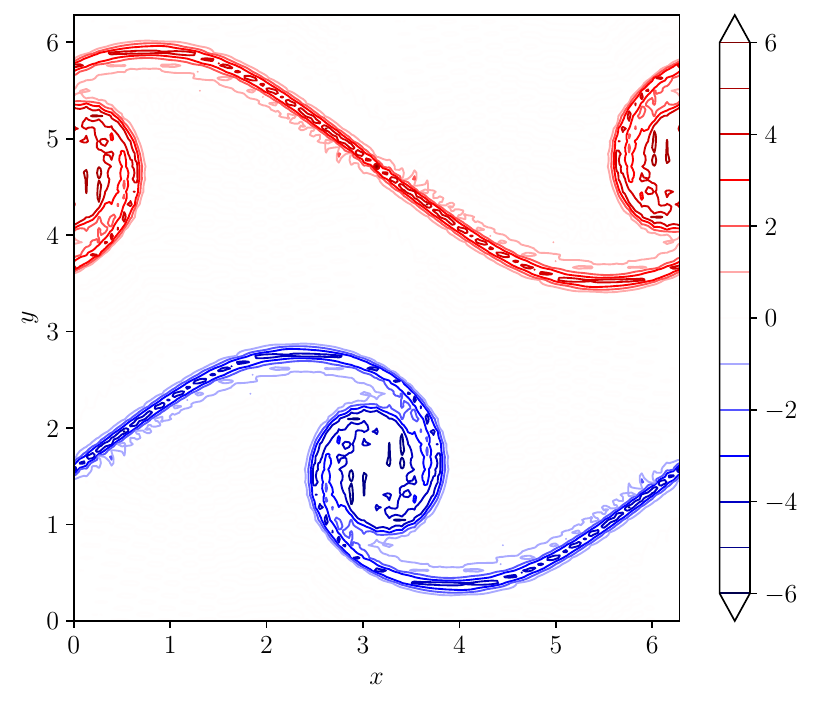}
        \caption{$t = 6$}
        % \label{fig:sub2}
    \end{subfigure}
    \begin{subfigure}{0.32\textwidth}
        \includegraphics[width=\textwidth]{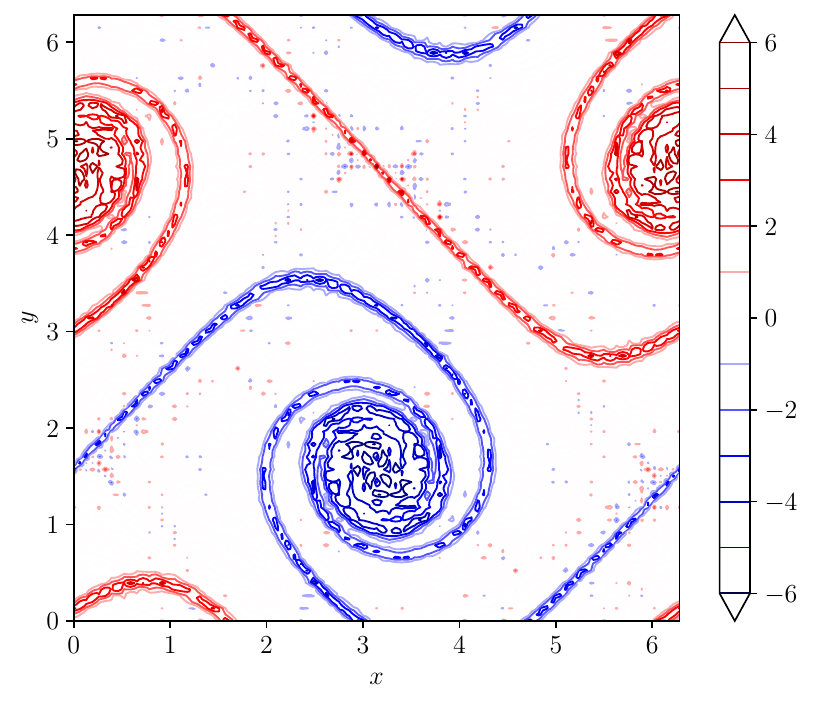}
        \caption{$t = 8$}
        % \label{fig:sub3}
    \end{subfigure}
    \caption{Reference vorticity solution on the refined mesh from \Cref{fig:msh_fine} with $\Delta t = 0.001$}
    \label{fig:omega_DNS}
\end{figure}

\begin{figure}[H]
    \centering
    \begin{subfigure}{0.32\textwidth}
        \includegraphics[width=\textwidth]{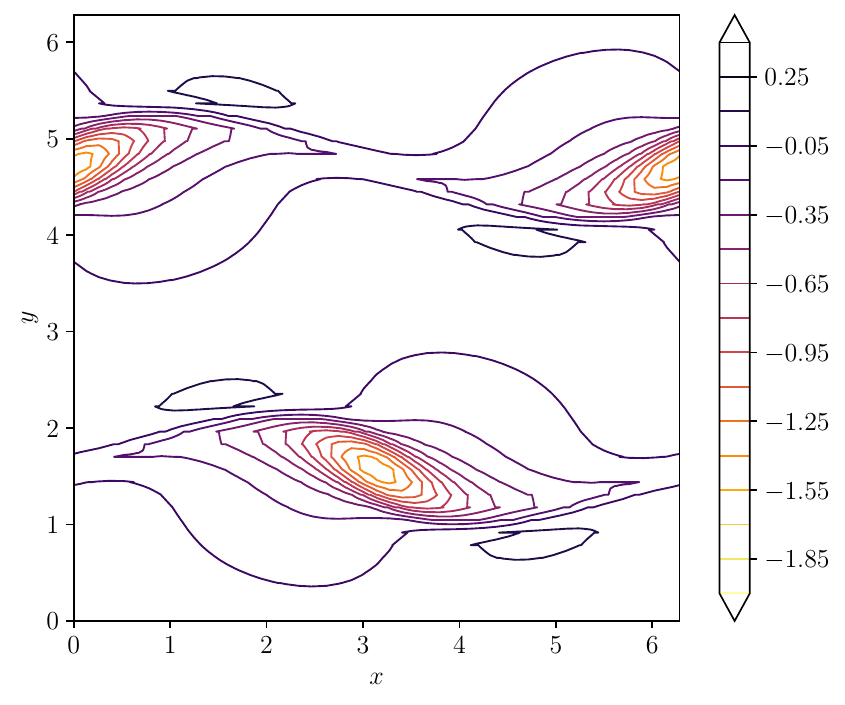}
        \caption{$t = 4$}
        % \label{fig:sub1}
    \end{subfigure}
    \begin{subfigure}{0.32\textwidth}
        \includegraphics[width=\textwidth]{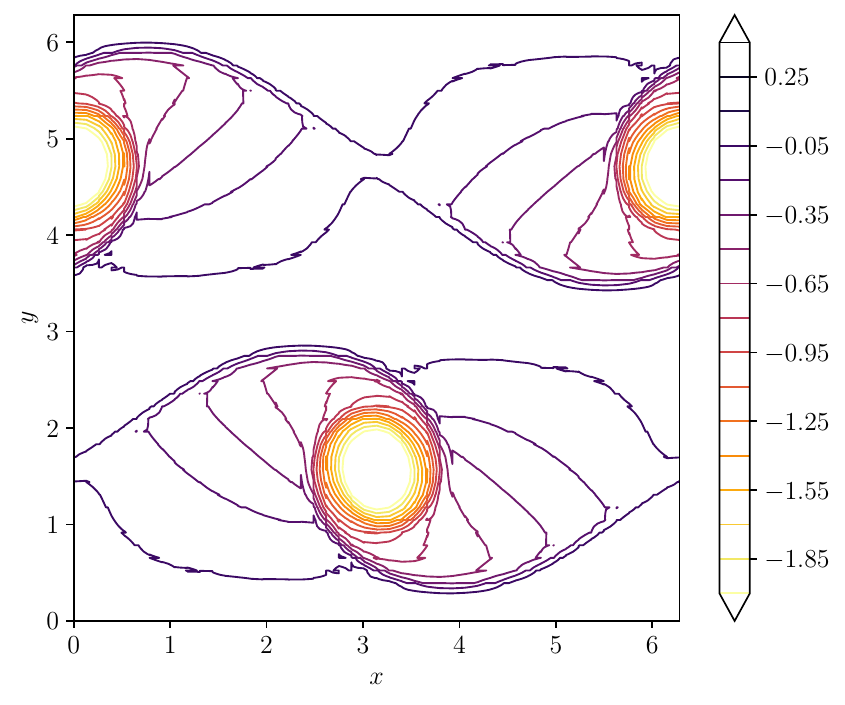}
        \caption{$t = 6$}
        % \label{fig:sub2}
    \end{subfigure}
    \begin{subfigure}{0.32\textwidth}
        \includegraphics[width=\textwidth]{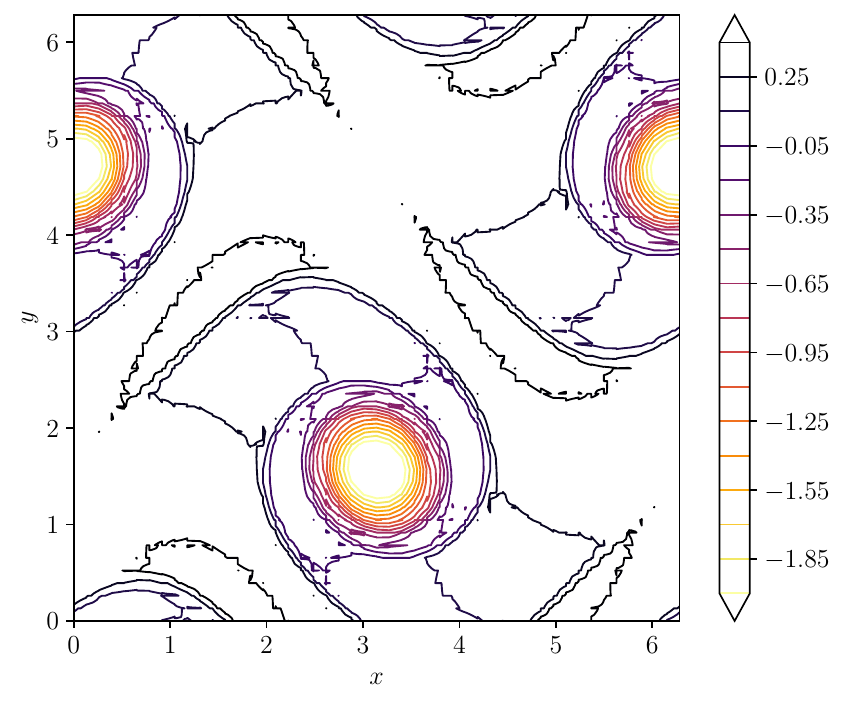}
        \caption{$t = 8$}
        % \label{fig:sub3}
    \end{subfigure}
    \caption{Reference total pressure solution on the refined mesh from \Cref{fig:msh_fine} with $\Delta t = 0.001$}
    \label{fig:pressure_DNS}
\end{figure}

As a reference solution, we employ results from a highly-resolved simulation using the base Galerkin discretisation on a fine orthogonal mesh with $48 \times 48$ elements with polynomial degree 2, shown in \Cref{fig:msh_fine}. The highly-resolved simulation captures the evolution of the flow with high fidelity, resolving the full range of active scales. We visualise the vorticity and total pressure fields at three time levels $t = 4, 6, \text{and }8$ in \Cref{fig:omega_DNS} and \Cref{fig:pressure_DNS}.
% We note that despite the mesh being highly refined, spurious oscillations persist in the vorticity field. These are attributed to the conservative nature of the scheme and the presence of sharp gradients in the solution. 
In the absence of an analytical solution, this highly-resolved simulation serves as a reference, and the objective of our VMS formulation is to recover its projection on coarser grids.
%===========================


\begin{figure}[H]
    \centering
    \begin{subfigure}{0.32\textwidth}
        \includegraphics[width=\textwidth]{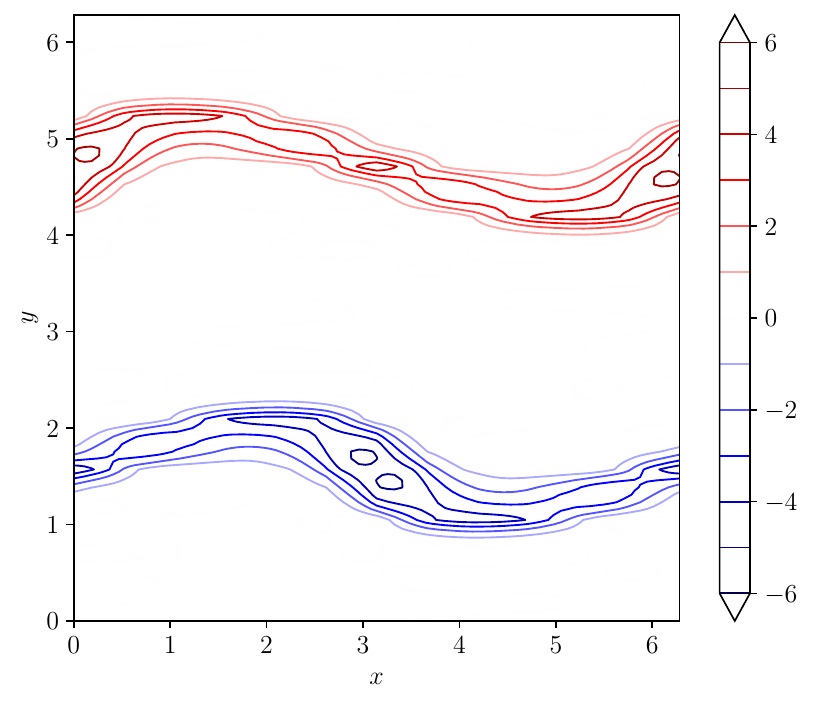}
        \caption{$t = 4$}
        % \label{fig:sub1}
    \end{subfigure}
    \begin{subfigure}{0.32\textwidth}
        \includegraphics[width=\textwidth]{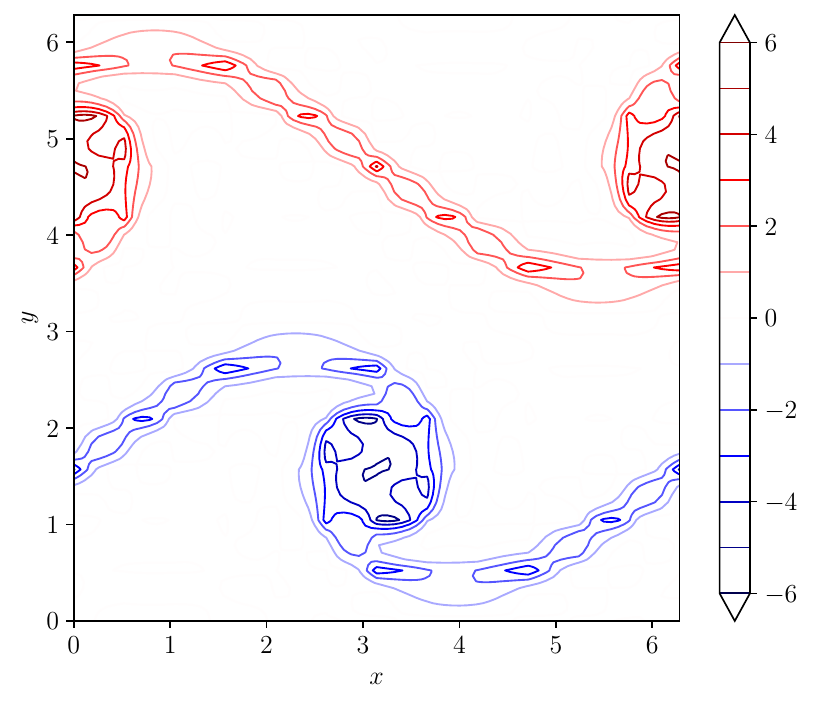}
        \caption{$t = 6$}
        % \label{fig:sub2}
    \end{subfigure}
    \begin{subfigure}{0.32\textwidth}
        \includegraphics[width=\textwidth]{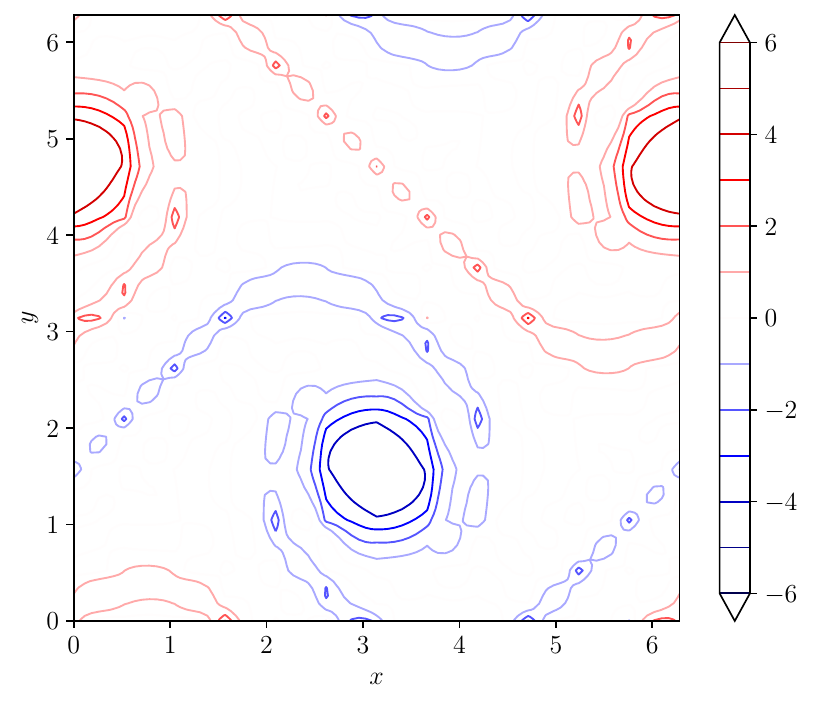}
        \caption{$t = 8$}
        % \label{fig:sub3}
    \end{subfigure}
    \caption{Projection of the highly-resolved vorticity solution onto the coarse mesh in \Cref{fig:msh_coarse}}
    \label{fig:omega_proj_ex}
\end{figure}

\begin{figure}[H]
    \centering
    \begin{subfigure}{0.32\textwidth}
        \includegraphics[width=\textwidth]{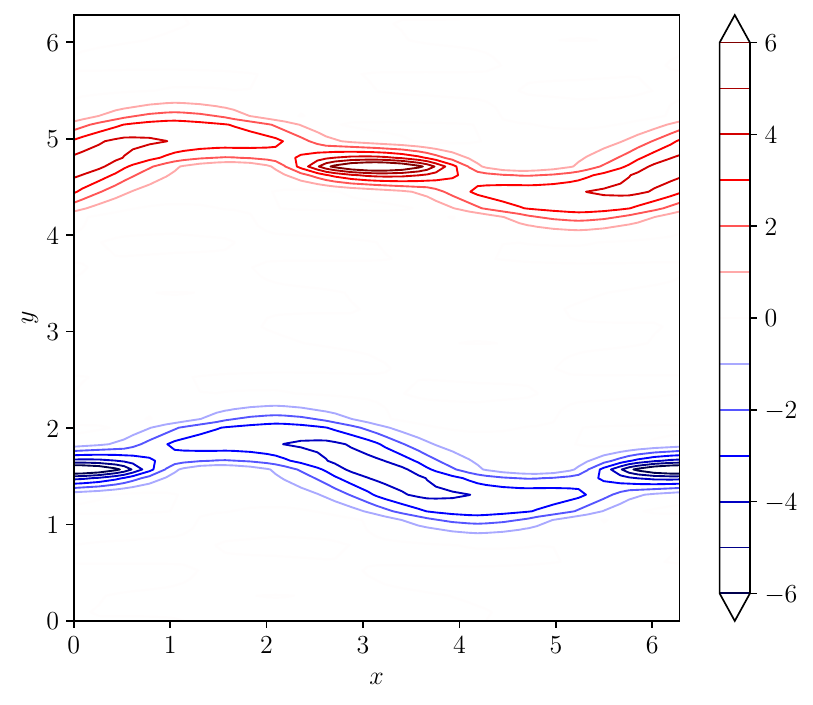}
        \caption{$t = 4$}
        % \label{fig:sub1}
    \end{subfigure}
    \begin{subfigure}{0.32\textwidth}
        \includegraphics[width=\textwidth]{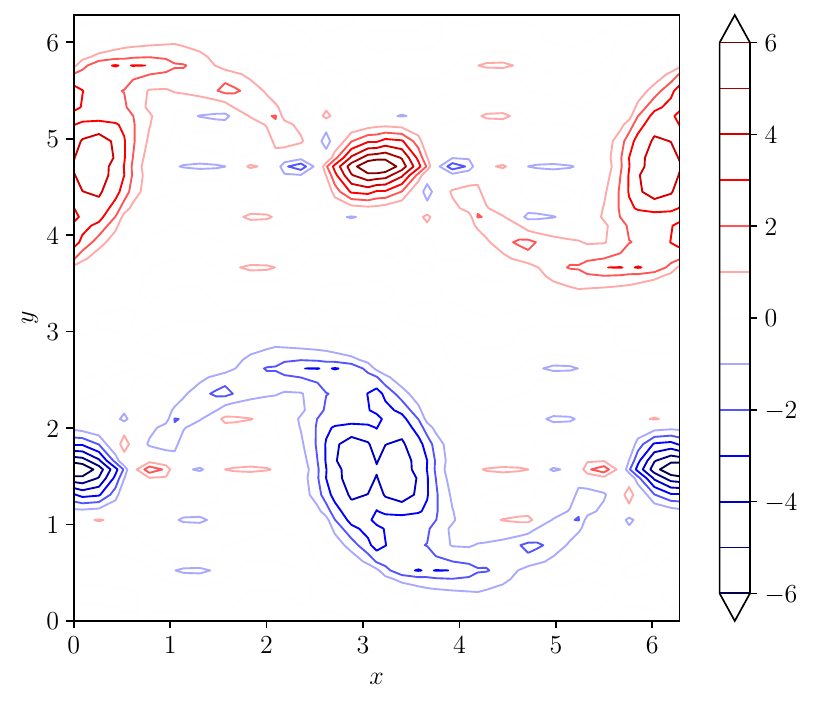}
        \caption{$t = 6$}
        % \label{fig:sub2}
    \end{subfigure}
    \begin{subfigure}{0.32\textwidth}
        \includegraphics[width=\textwidth]{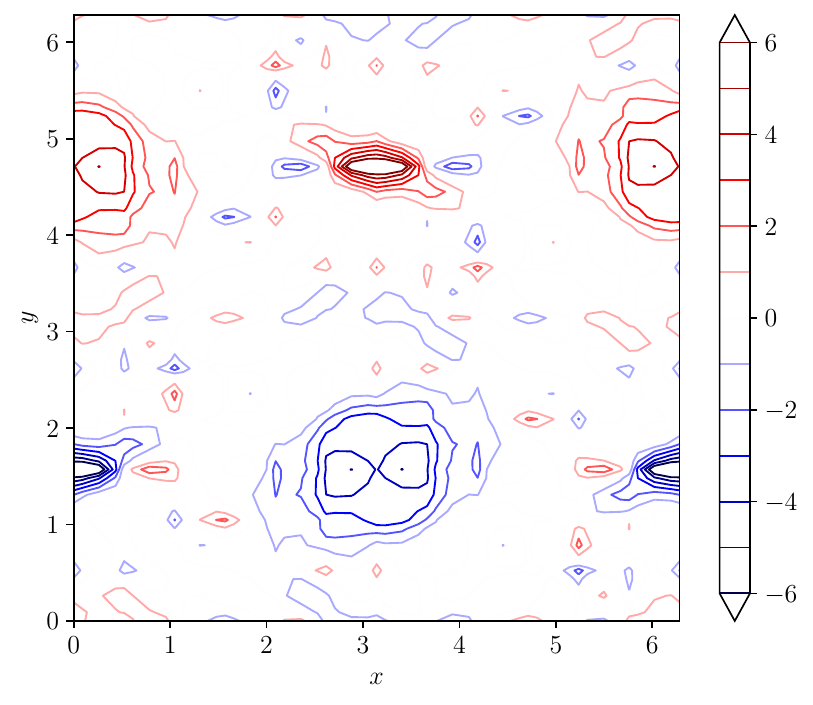}
        \caption{$t = 8$}
        % \label{fig:sub3}
    \end{subfigure}
    \caption{Vorticity field computed on the coarse mesh in \Cref{fig:msh_coarse} using the base Galerkin approach with $\Delta t = 0.001$}
    \label{fig:omega_glk}
\end{figure}

In \Cref{fig:omega_proj_ex}, we show the projection of the highly-resolved solution onto a coarse mesh with $12 \times 12$ elements with polynomial degree $2$ depicted in \Cref{fig:msh_coarse}. This projection represents the best possible approximation of the highly-resolved solution on the coarse mesh in the associated norm associated with the projector defined in \eqref{eq:NS_proj}. The projected fields serve as a target for comparison, illustrating the optimal behaviour expected from any coarse-scale method within the chosen resolution.
Next, we present the solution obtained using the base Galerkin scheme on the same coarse mesh in \Cref{fig:omega_glk}. While this solution remains numerically stable, it exhibits only a weak resemblance to the highly-resolved solution and its projection. The coarse Galerkin solution fails to capture several key flow features and develops spurious oscillations that are not present in the projected reference.

\begin{figure}[htp]
    \centering
    \begin{subfigure}{0.32\textwidth}
        \includegraphics[width=\textwidth]{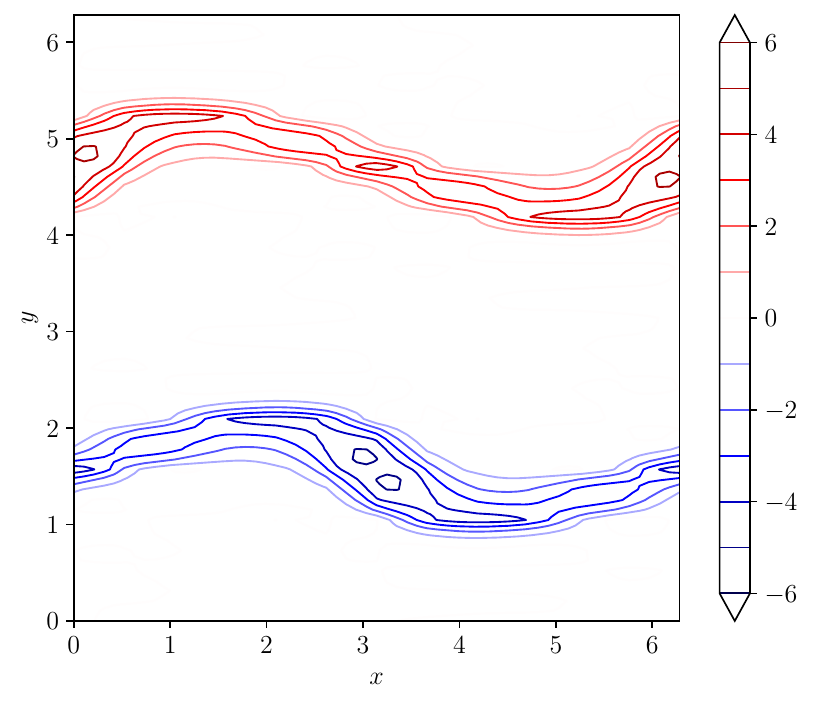}
        \caption{$t = 4$}
        % \label{fig:sub1}
    \end{subfigure}
    \begin{subfigure}{0.32\textwidth}
        \includegraphics[width=\textwidth]{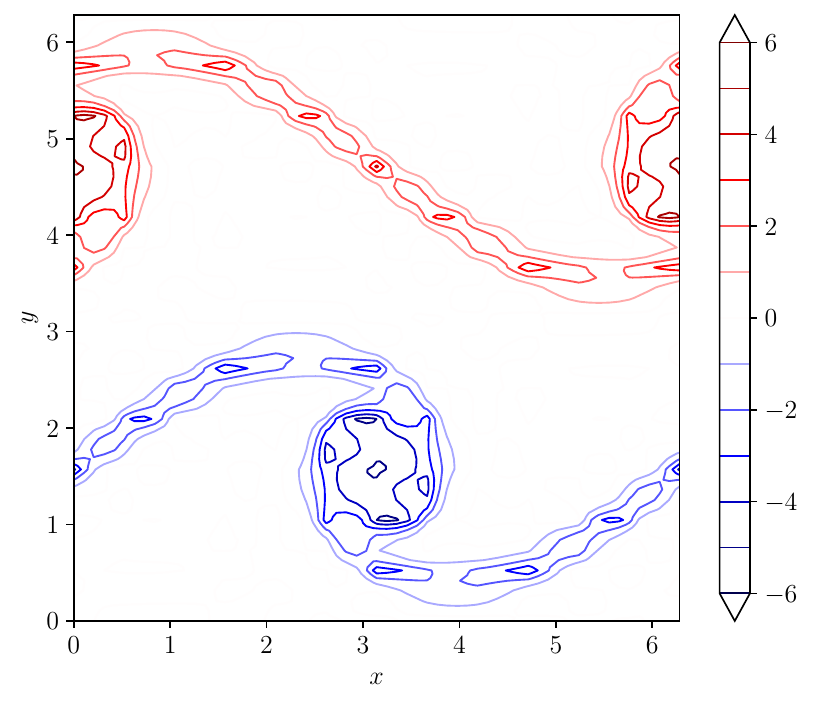}
        \caption{$t = 6$}
        % \label{fig:sub2}
    \end{subfigure}
    \begin{subfigure}{0.32\textwidth}
        \includegraphics[width=\textwidth]{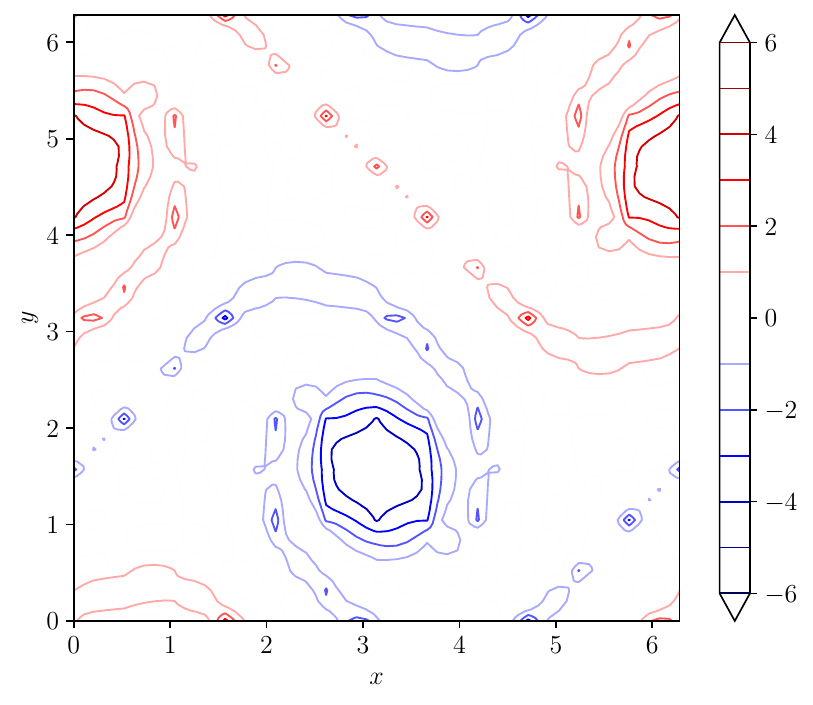}
        \caption{$t = 8$}
    \end{subfigure}
    \caption{Vorticity field computed on the coarse mesh in \Cref{fig:msh_coarse} using the VMS approach with $k = 4$ and $\Delta t = 0.001$}
    \label{fig:VMS_vorticity}
\end{figure}

In contrast, the results obtained with the proposed VMS formulation demonstrate significantly improved agreement with the projected reference solution. As observed through plots in \Cref{fig:VMS_vorticity}, the VMS solution with $k = 4$ closely tracks the target projection across all time levels. 
\begin{figure}[htp]
    \centering
    \begin{subfigure}{0.32\textwidth}
        \includegraphics[width=\textwidth]{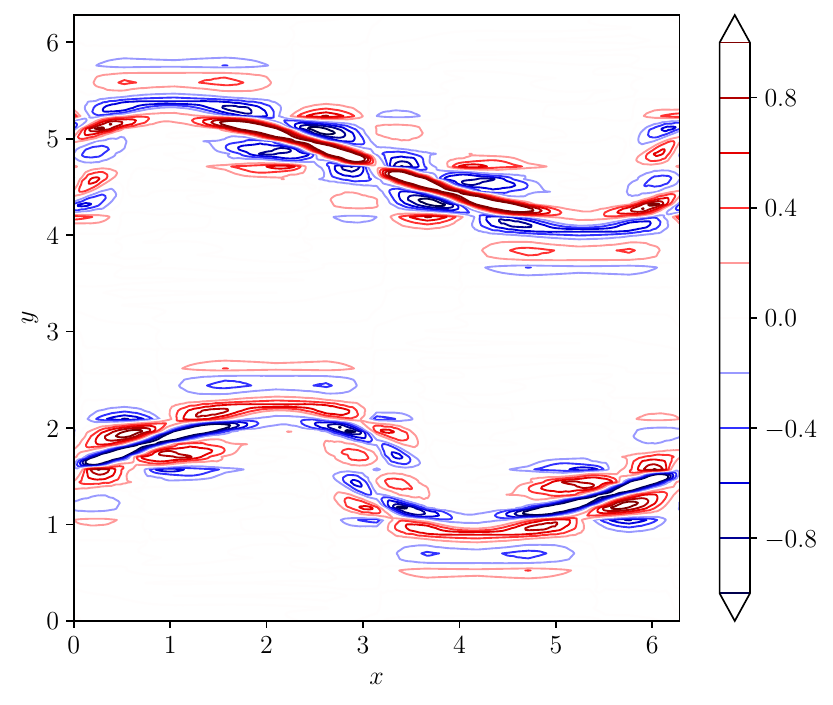}
        \caption{$t = 4$}
        % \label{fig:sub1}
    \end{subfigure}
    \begin{subfigure}{0.32\textwidth}
        \includegraphics[width=\textwidth]{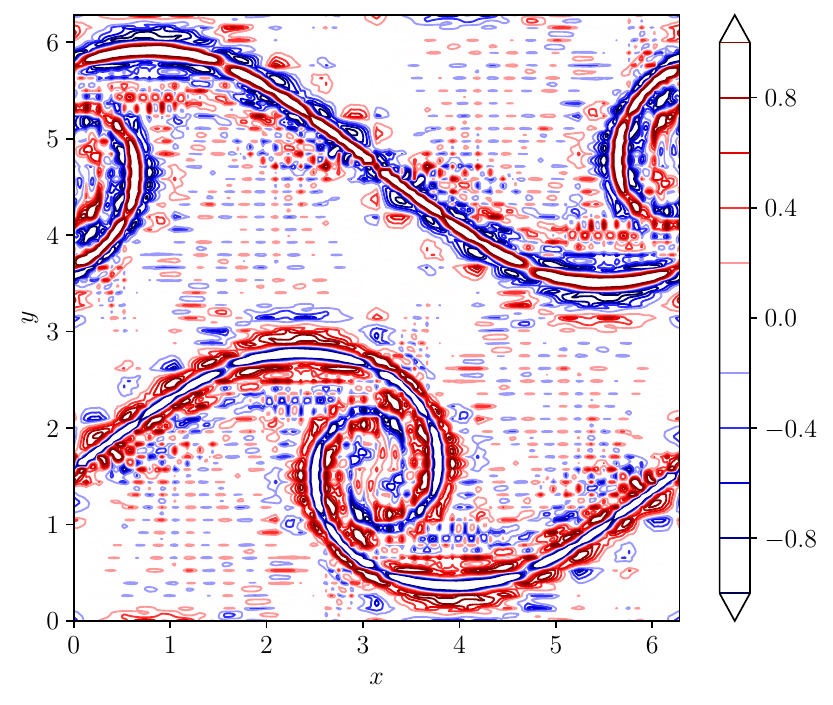}
        \caption{$t = 6$}
        % \label{fig:sub2}
    \end{subfigure}
    \begin{subfigure}{0.32\textwidth}
        \includegraphics[width=\textwidth]{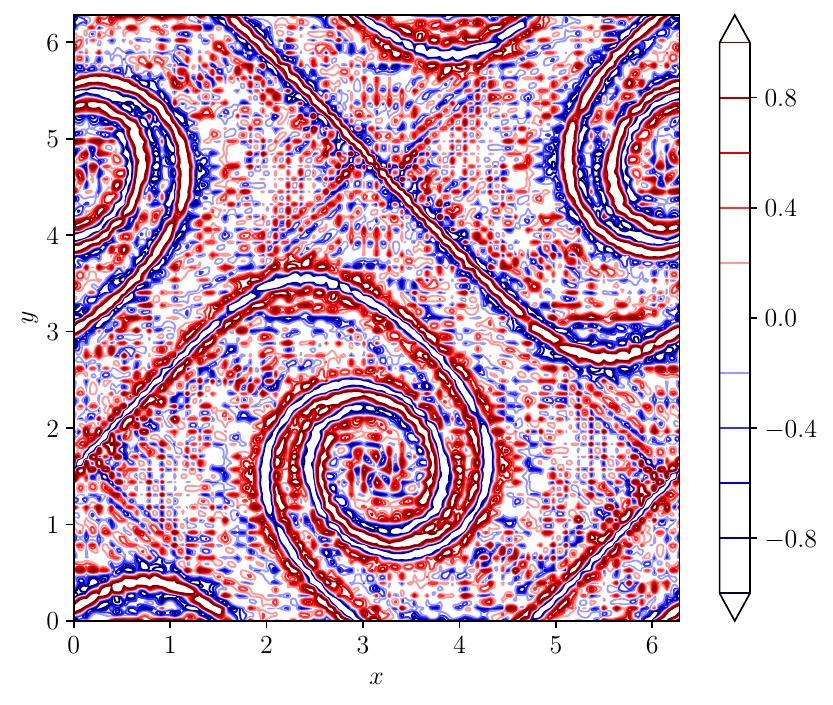}
        \caption{$t = 8$}
        \label{fig:VMS_vorticity_prime_ex_t8}
    \end{subfigure}
    \caption{Exact unresolved vorticity field computed with the highly-resolved reference solution and its projection onto the coarse mesh in \Cref{fig:msh_coarse}}
    \label{fig:omega_prime_ex}
\end{figure}
\begin{figure}[htp]
    \centering
    \begin{subfigure}{0.32\textwidth}
        \includegraphics[width=\textwidth]{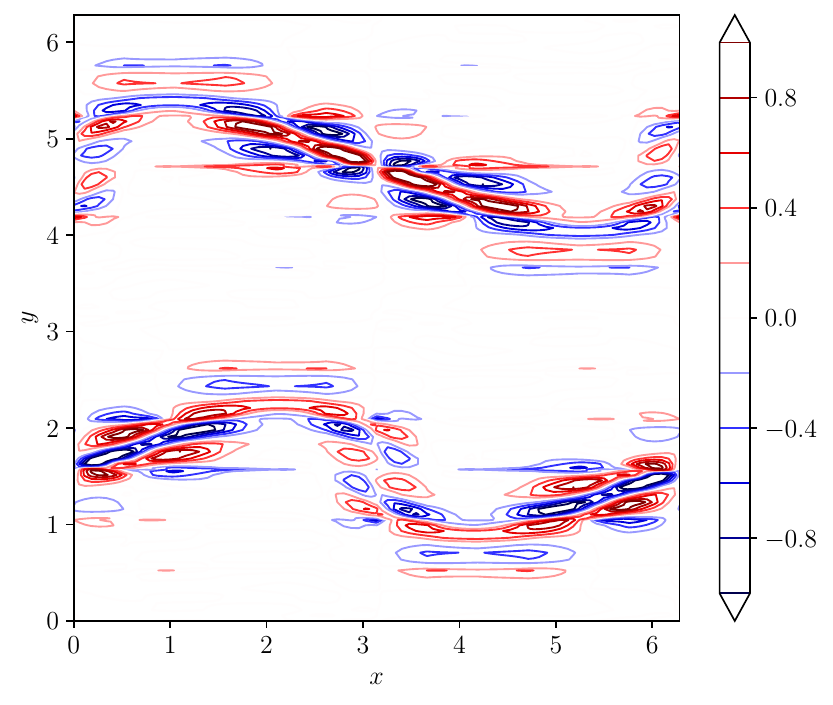}
        \caption{$t = 4$}
        % \label{fig:sub1}
    \end{subfigure}
    \begin{subfigure}{0.32\textwidth}
        \includegraphics[width=\textwidth]{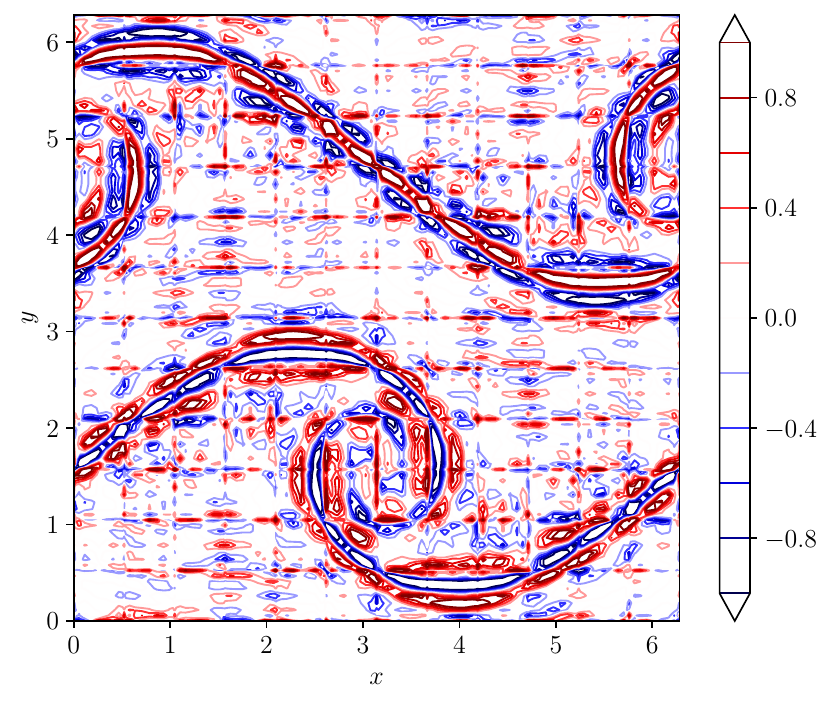}
        \caption{$t = 6$}
        % \label{fig:sub2}
    \end{subfigure}
    \begin{subfigure}{0.32\textwidth}
        \includegraphics[width=\textwidth]{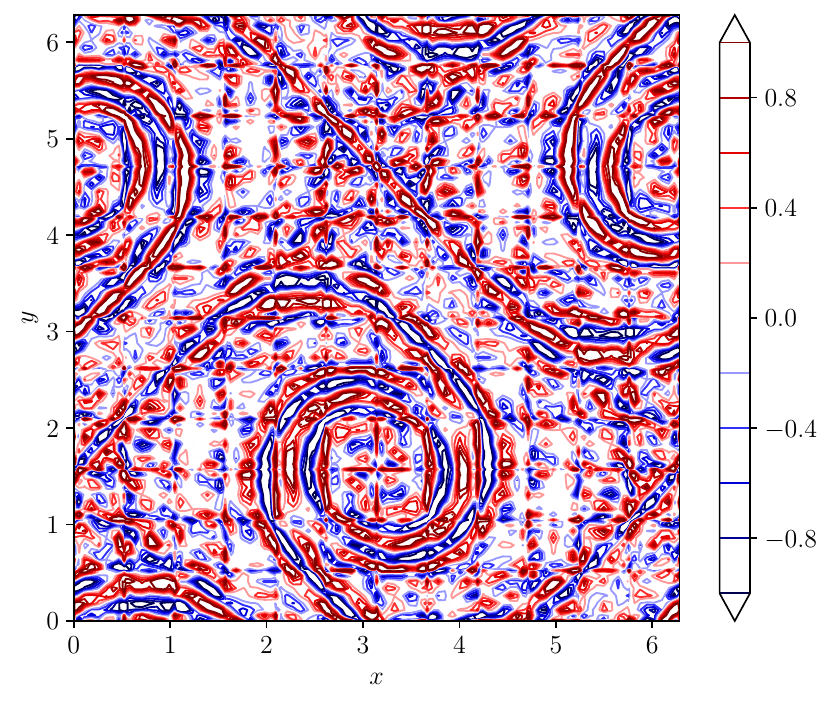}
        \caption{$t = 8$}
        \label{fig:VMS_vorticity_prime_t8}
    \end{subfigure}
    \caption{Approximate unresolved vorticity field computed using the VMS approach with $k = 4$ and $\Delta t = 0.001$}
    \label{fig:omega_prime}
\end{figure}
Moreover, when comparing the exact unresolved vorticity in \Cref{fig:omega_prime_ex} with the approximated ones in \Cref{fig:VMS_vorticity}, we note that the VMS approach provides a modest yet consistent approximation of the unresolved fine-scale content. However, as illustrated in \Cref{fig:VMS_vorticity_prime_t8}, some discrepancies emerge when the solution exhibits large, sharp gradients. In such cases, the approximation naturally becomes more challenging, and fine-scale structures are more difficult to accurately capture. Despite this, the computed unresolved scales still display a clear degree of coherence, with identifiable spiral-like structures that bear a loose resemblance to those in the exact fine-scale solution in \Cref{fig:VMS_vorticity_prime_ex_t8}. This suggests that, even in the presence of steep gradients, the VMS framework retains the ability to reconstruct relevant unresolved features with reasonable fidelity. 

\begin{figure}[H]
    \centering
    \begin{subfigure}{0.32\textwidth}
        \includegraphics[width=\textwidth]{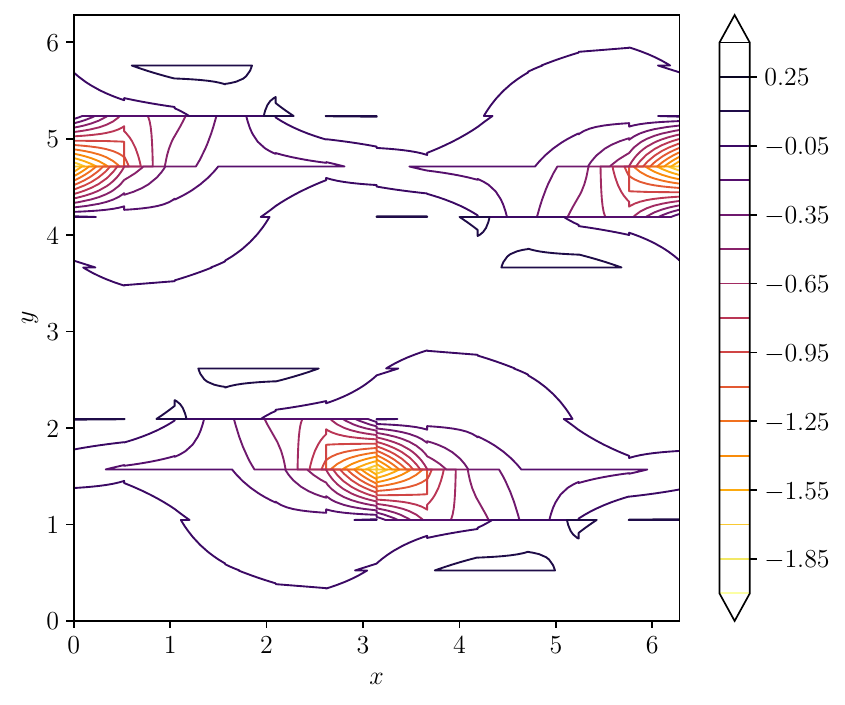}
        \caption{$t = 4$}
        % \label{fig:sub1}
    \end{subfigure}
    \begin{subfigure}{0.32\textwidth}
        \includegraphics[width=\textwidth]{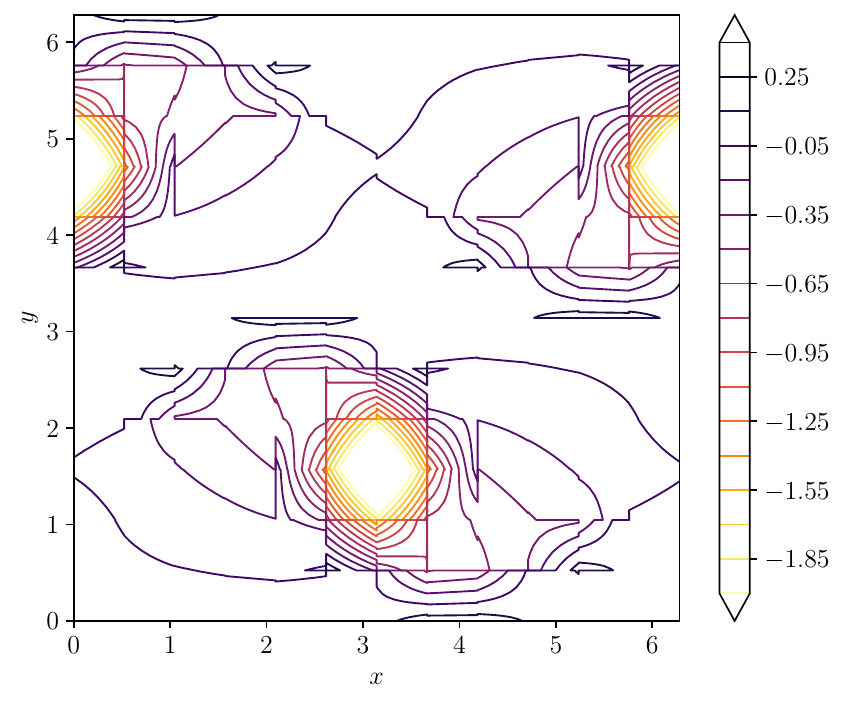}
        \caption{$t = 6$}
        % \label{fig:sub2}
    \end{subfigure}
    \begin{subfigure}{0.32\textwidth}
        \includegraphics[width=\textwidth]{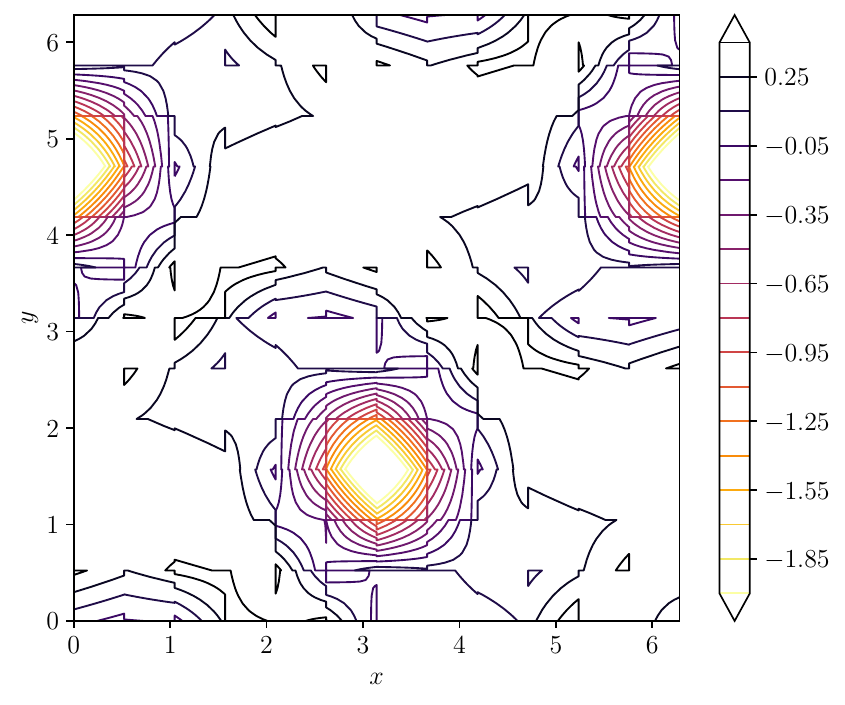}
        \caption{$t = 8$}
        % \label{fig:sub3}
    \end{subfigure}
    \caption{Projection of the highly-resolved total pressure solution onto the coarse mesh in \Cref{fig:msh_coarse}}
    \label{fig:pressure_proj_ex}
\end{figure}

\begin{figure}[H]
    \centering
    \begin{subfigure}{0.32\textwidth}
        \includegraphics[width=\textwidth]{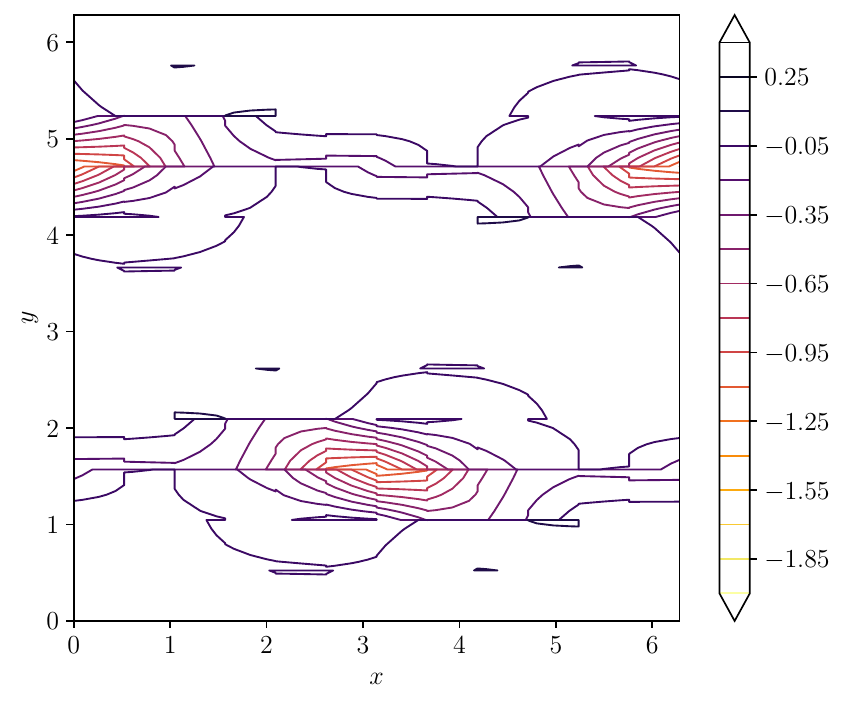}
        \caption{$t = 4$}
        % \label{fig:sub1}
    \end{subfigure}
    \begin{subfigure}{0.32\textwidth}
        \includegraphics[width=\textwidth]{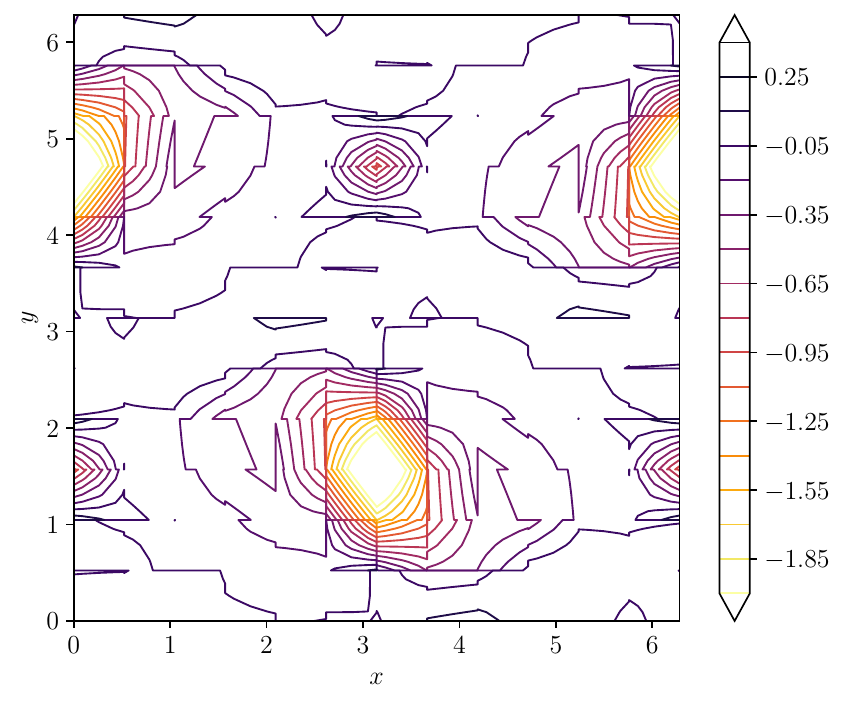}
        \caption{$t = 6$}
        % \label{fig:sub2}
    \end{subfigure}
    \begin{subfigure}{0.32\textwidth}
        \includegraphics[width=\textwidth]{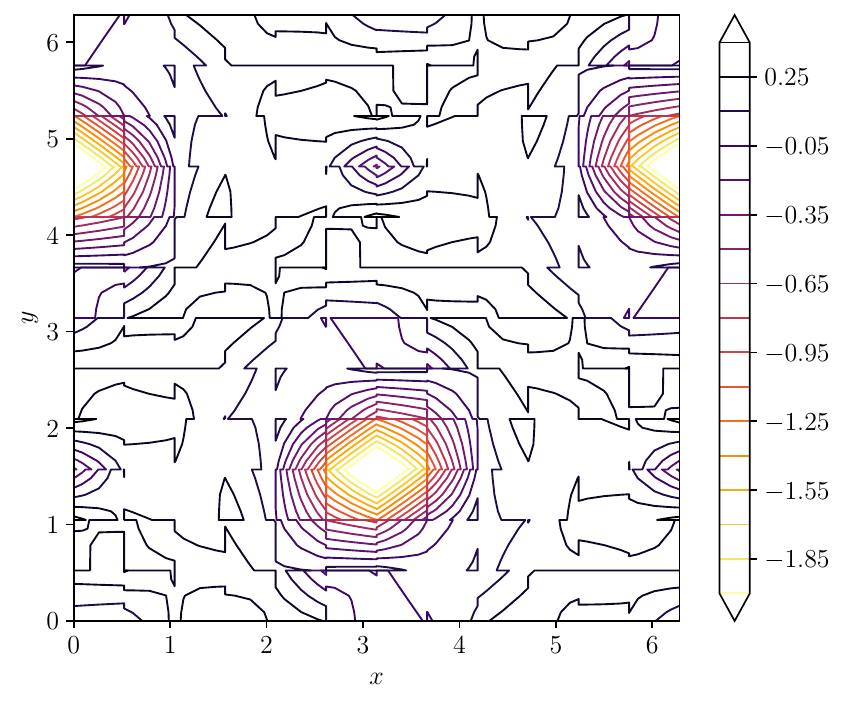}
        \caption{$t = 8$}
        % \label{fig:sub3}
    \end{subfigure}
    \caption{Total pressure field computed on the coarse mesh in \Cref{fig:msh_coarse} using the base Galerkin approach with $\Delta t = 0.001$}
    \label{fig:pressure_glk}
\end{figure}

\begin{figure}[H]
    \centering
    \begin{subfigure}{0.32\textwidth}
        \includegraphics[width=\textwidth]{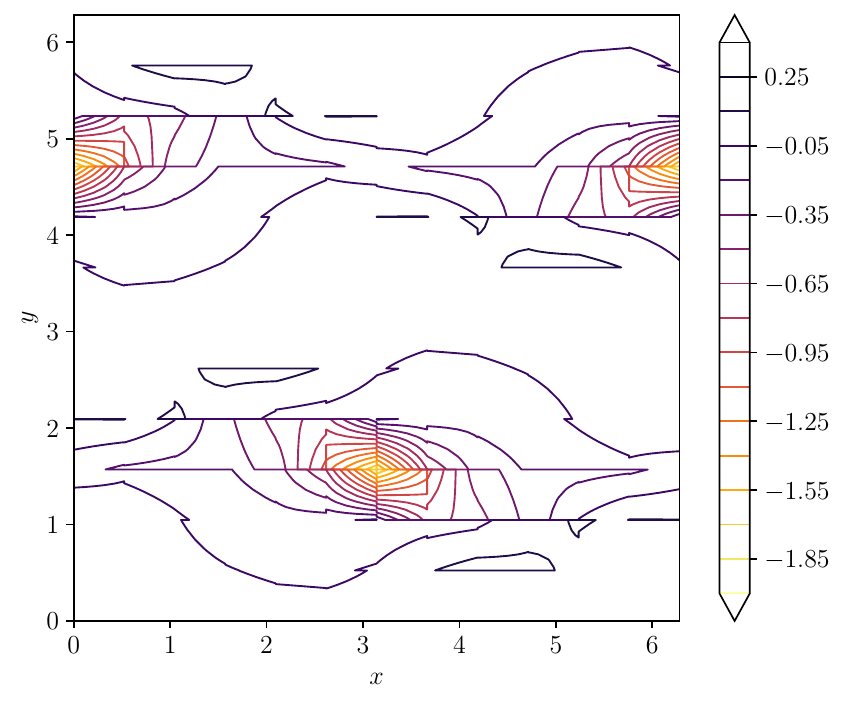}
        \caption{$t = 4$}
        % \label{fig:sub1}
    \end{subfigure}
    \begin{subfigure}{0.32\textwidth}
        \includegraphics[width=\textwidth]{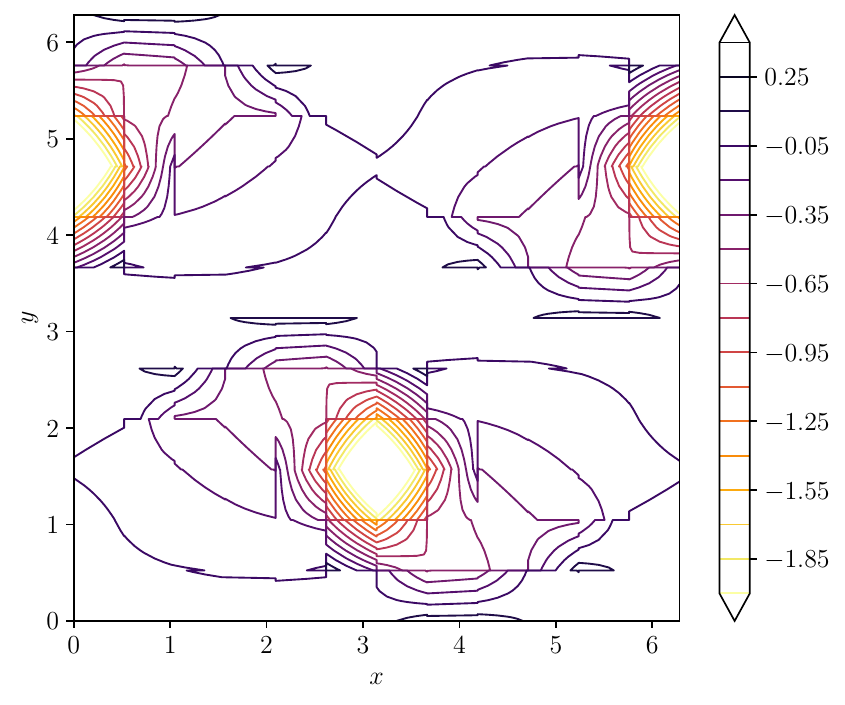}
        \caption{$t = 6$}
        % \label{fig:sub2}
    \end{subfigure}
    \begin{subfigure}{0.32\textwidth}
        \includegraphics[width=\textwidth]{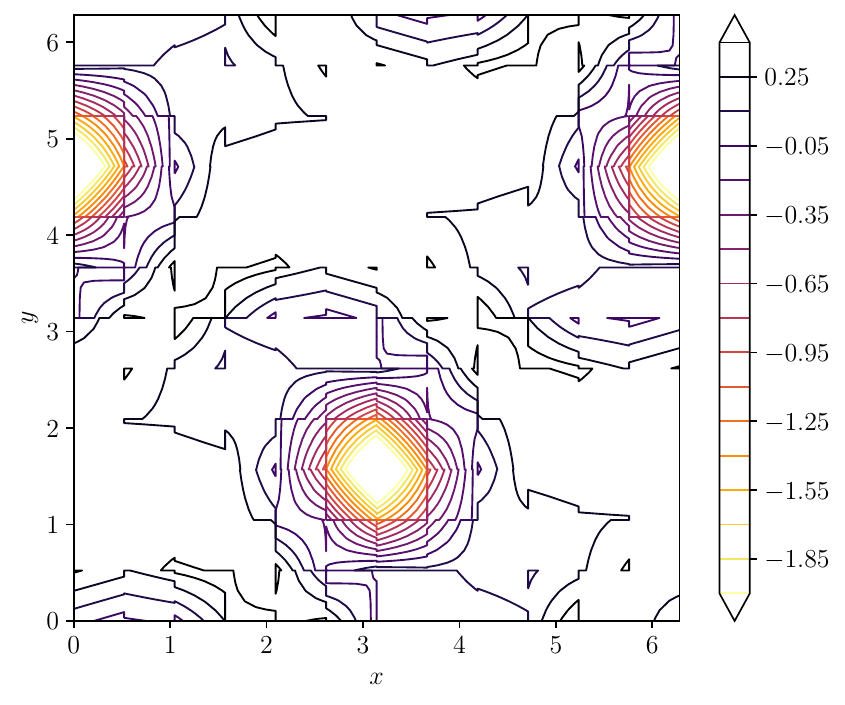}
        \caption{$t = 8$}
        % \label{fig:sub3}
    \end{subfigure}
    \caption{Total pressure field computed on the coarse mesh in \Cref{fig:msh_coarse} using the VMS approach with $k = 4$ and $\Delta t = 0.001$}
    \label{fig:pressure_VMS}
\end{figure}

\begin{figure}[H]
    \centering
    \begin{subfigure}{0.32\textwidth}
        \includegraphics[width=\textwidth]{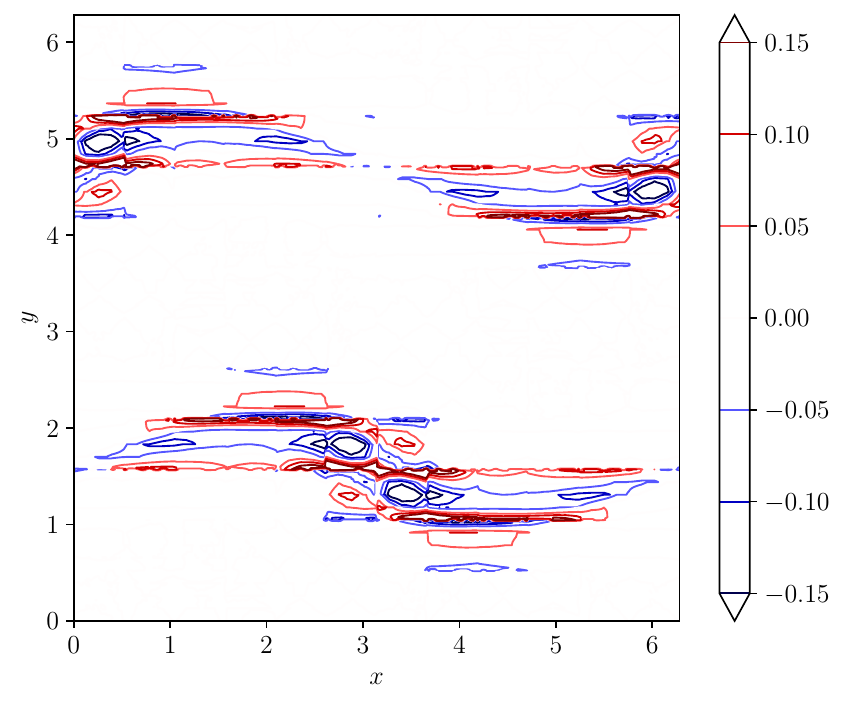}
        \caption{$t = 4$}
        % \label{fig:sub1}
    \end{subfigure}
    \begin{subfigure}{0.32\textwidth}
        \includegraphics[width=\textwidth]{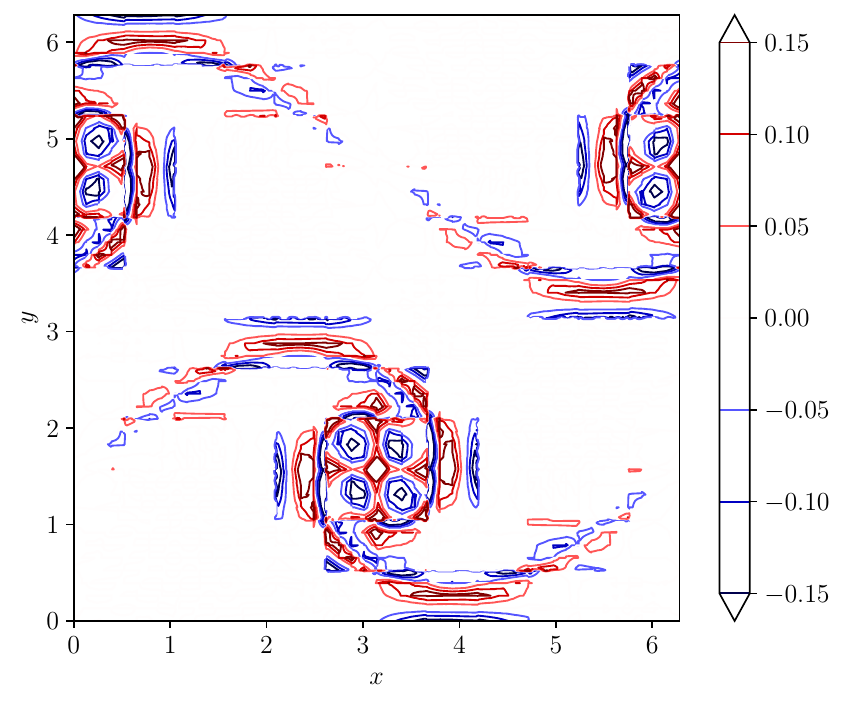}
        \caption{$t = 6$}
        % \label{fig:sub2}
    \end{subfigure}
    \begin{subfigure}{0.32\textwidth}
        \includegraphics[width=\textwidth]{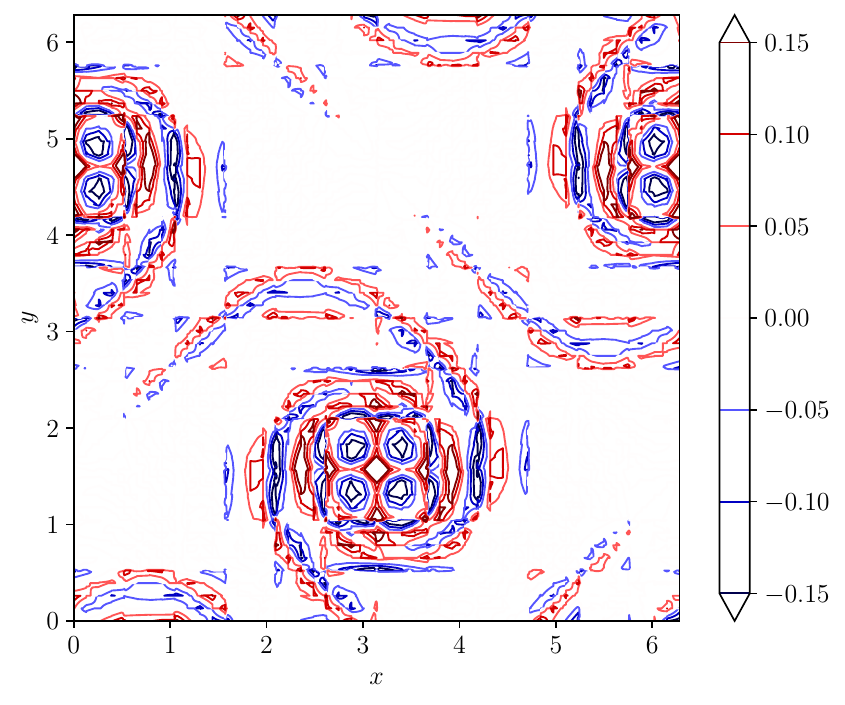}
        \caption{$t = 8$}
        % \label{fig:sub3}
    \end{subfigure}
    \caption{Exact unresolved total pressure field computed with the highly-resolved reference solution and its projection onto the coarse mesh in \Cref{fig:msh_coarse}}
    \label{fig:pressure_prime_ex}
\end{figure}

\begin{figure}[H]
    \centering
    \begin{subfigure}{0.32\textwidth}
        \includegraphics[width=\textwidth]{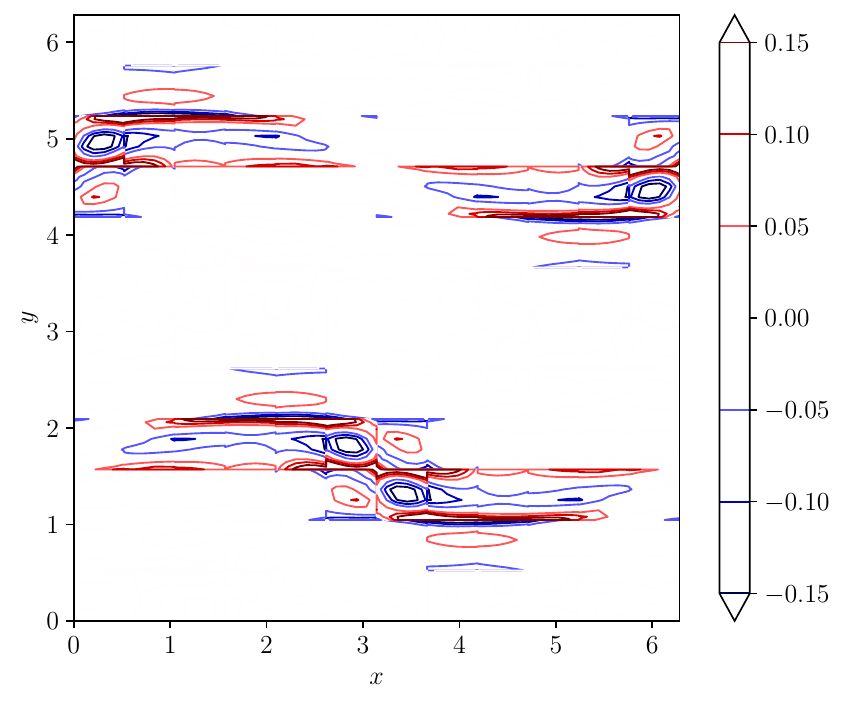}
        \caption{$t = 4$}
        % \label{fig:sub1}
    \end{subfigure}
    \begin{subfigure}{0.32\textwidth}
        \includegraphics[width=\textwidth]{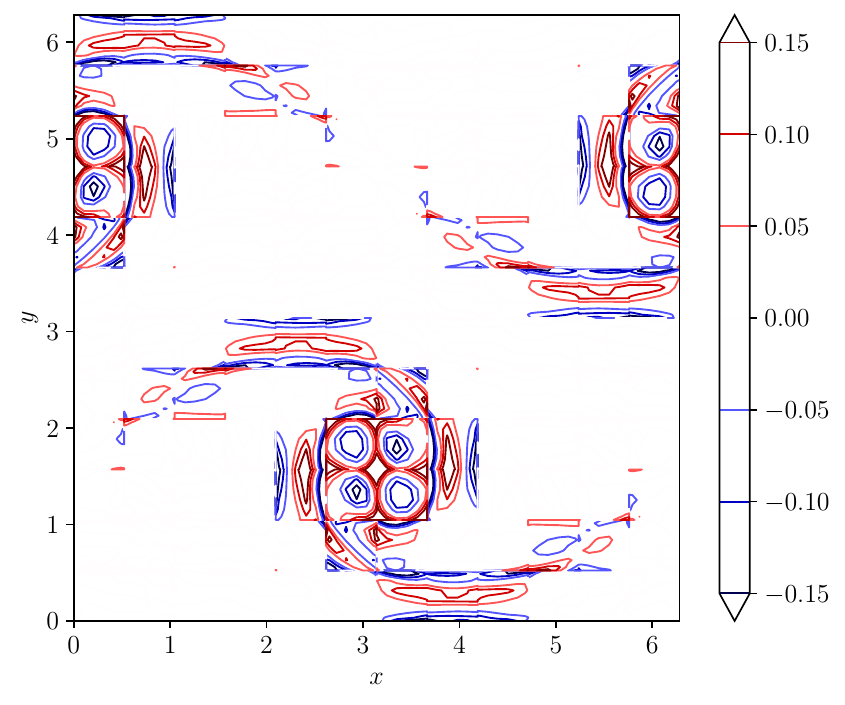}
        \caption{$t = 6$}
        % \label{fig:sub2}
    \end{subfigure}
    \begin{subfigure}{0.32\textwidth}
        \includegraphics[width=\textwidth]{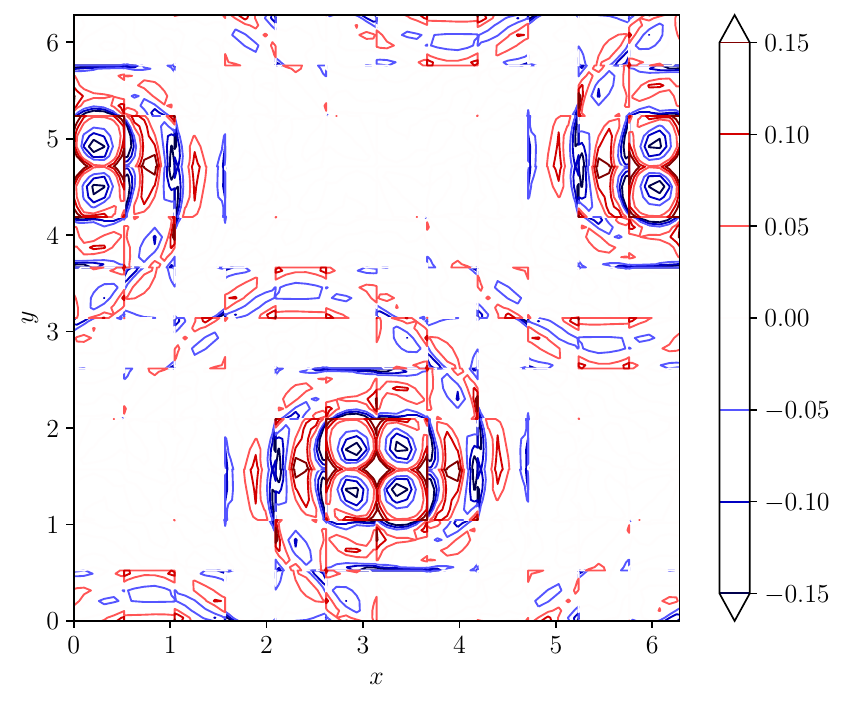}
        \caption{$t = 8$}
        % \label{fig:sub3}
    \end{subfigure}
    \caption{Approximate unresolved total pressure field computed using the VMS approach with $k = 4$ and $\Delta t = 0.001$}
    \label{fig:pressure_prime}
\end{figure}

Similar observations hold for the total pressure field. The Galerkin solution, shown in \Cref{fig:pressure_glk}, exhibits significant deviations and fails to accurately capture the features present in either the exact projection in \Cref{fig:pressure_proj_ex} or the highly-resolved reference solution. In contrast, the VMS solution shows a strong agreement with the exact projection as evident through \Cref{fig:pressure_VMS}, indicating the effectiveness of the approach in recovering pressure fields consistent with the defined optimal projection. Regarding the unresolved pressure scales, the VMS approximation continues to provide a coherent representation that aligns well with the exact unresolved component. This is apparent through the similarities between the plots of the exact unresolved total pressure field in \Cref{fig:pressure_prime_ex} and the VMS approximation thereof in \Cref{fig:pressure_prime}. Unlike the vorticity field, the total pressure field avoids the difficulties observed earlier and maintains close agreement with the exact unresolved total pressure field throughout time.

% \begin{figure}[H]
%     \centering
%     \begin{subfigure}{0.32\textwidth}
%         \includegraphics[width=\textwidth]{Images/vortex_roll_up/figs_anim/VMS_overInt_10_Re0_12x12_p2_pf6_dt_1.0e-03_omega_bar_t_idx_3.pdf}
%         \caption{$\bar{\omega}$}
%         % \label{fig:sub2}
%     \end{subfigure}
%     \begin{subfigure}{0.32\textwidth}
%         \includegraphics[width=\textwidth]{Images/vortex_roll_up/figs_anim/VMS_overInt_10_Re0_12x12_p2_pf6_dt_1.0e-03_omega_prime_t_idx_3.pdf}
%         \caption{$\omega'$}
%         % \label{fig:sub3}
%     \end{subfigure}
%     \begin{subfigure}{0.32\textwidth}
%         \includegraphics[width=\textwidth]{Images/vortex_roll_up/figs_anim/VMS_overInt_10_Re0_12x12_p2_pf6_dt_1.0e-03_omega_t_idx_3.pdf}
%         \caption{$\bar{\omega} + \omega'$}
%         % \label{fig:sub1}
%     \end{subfigure}
%     \caption{Approximate unresolved total pressure field computed using the VMS approach with $k = 4$ and $\Delta t = 0.001$}
%     \label{fig:omega_bar_prime}
% \end{figure}

\begin{figure}[H]
    \centering
    \begin{subfigure}{0.32\textwidth}
        \includegraphics[width=\textwidth]{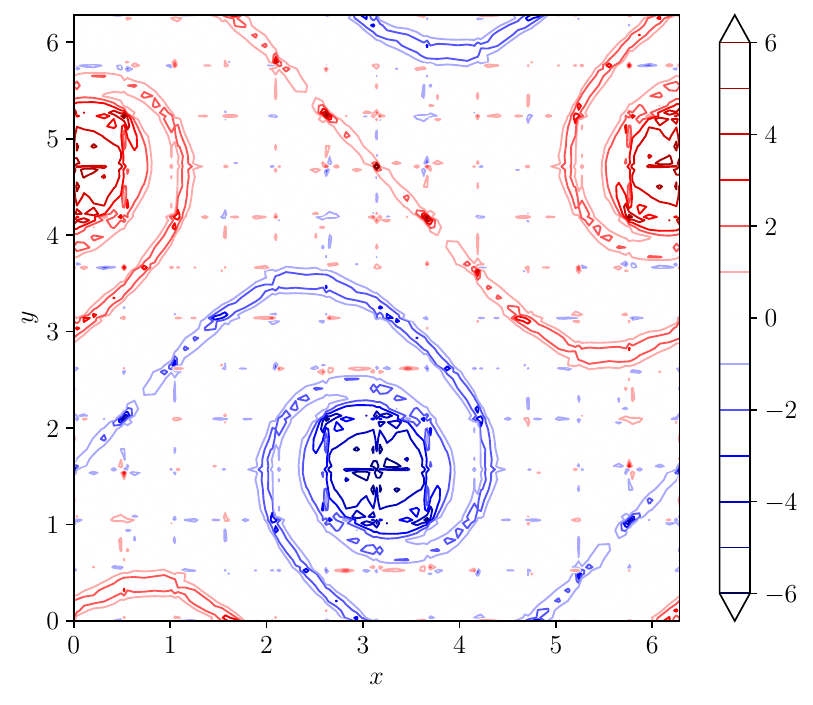}
        \caption{VMS solution $\bar{\omega} + \omega'$ with $p = 2$, $k = 4$}
        % \label{fig:sub1}
    \end{subfigure}
    \begin{subfigure}{0.32\textwidth}
        \includegraphics[width=\textwidth]{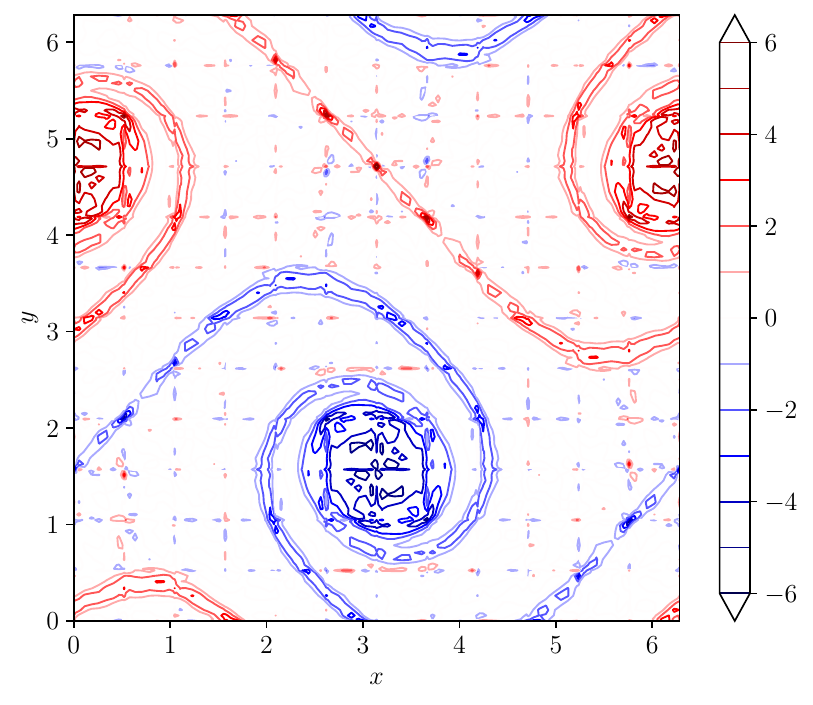}
        \caption{Galerkin solution with $p = 6$}
        % \label{fig:sub2}
    \end{subfigure}
    \begin{subfigure}{0.32\textwidth}
        \includegraphics[width=\textwidth]{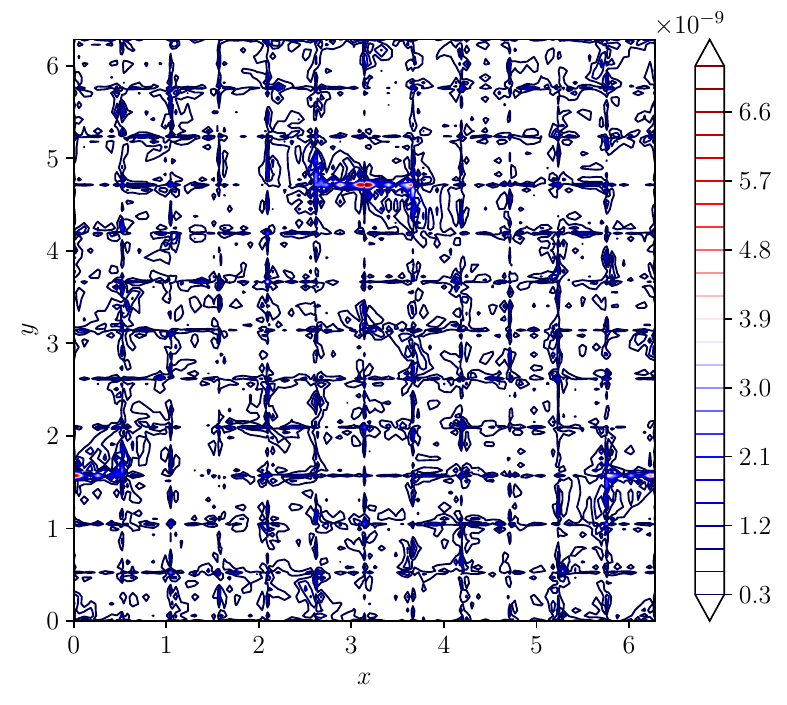}
        \caption{Difference between VMS and Galerkin}
        % \label{fig:sub3}
    \end{subfigure}
    \caption{Full VMS solution and the base Galerkin solution on a mesh with $12 \times 12$ elements at $t = 8$ with $\Delta t = 0.001$}
    \label{fig:omega_VMS_vs_Glk}
\end{figure}
Thus far, we have examined the resolved and unresolved components of the VMS solution largely in isolation. However, since computational effort is invested in approximating the unresolved scales on the fine space, it is natural to consider the combined solution. This also enables a fair comparison with a standard Galerkin method on an equivalently refined mesh. Specifically, we compare the full VMS approximation,
\(
\bar{\omega} + \omega',
\)
computed on a coarse mesh with polynomial degree \(p\) and a fine-scale enrichment of \(k\), with the Galerkin solution obtained using polynomial degree \(p + k\). A comparison of these solutions is shown in \Cref{fig:omega_VMS_vs_Glk} for the vorticity field at the final time, which reveals only negligible differences, with a maximum pointwise discrepancy of order \(\mathcal{O}(10^{-9})\). This demonstrates that the full VMS solution effectively reproduces the equivalently refined Galerkin solution, despite relying on a fundamentally different formalism and algorithmic strategy, as detailed in \Cref{sec:VMS_cost}.

Having established this agreement, we now compare the computational cost of the proposed VMS methodology with that of the standard Galerkin method. As described in \Cref{sec:VMS_cost}, we evaluate the VMS cost using coarse and fine meshes of polynomial degrees $p$ and $p+k$, respectively, and compare it against the cost of the standard Galerkin method applied directly on the fine mesh of degree $p+k$. Computational time is used as the cost metric, and the corresponding results are reported in \Cref{tab:time_vms} and \Cref{tab:time_glk}. Each data point is based on two repetitions, with variations between runs consistently below $1\,\mathrm{s}$. All cost estimations were performed in a fair and consistent manner using implementations of both the proposed VMS formulation (\Cref{alg:VMS}) and the baseline Galerkin formulation (\Cref{alg:Glk}) within the open-source package Firedrake \cite{FiredrakeUserManual}. All computations employed direct solvers for all linear systems and were executed on the same machine\footnote{PC with an AMD Ryzen 7 4800H and 32 GB DDR4 RAM} in single-core execution.

The colour coding in the tables groups cases that share the same fine-space polynomial degree: all entries in \textcolor{blue}{blue} correspond to fine meshes with $p+k=4$, all entries in \textcolor{magenta}{magenta} correspond to $p+k = 5$, and so on. Alongside the compute time, the cost ratio between the VMS and Galerkin methods, $R$, introduced in \Cref{sec:VMS_cost} is also reported in \Cref{tab:time_vms}. As noted in \Cref{sec:VMS_cost}, the proposed VMS approach becomes more favourable when the cost of assembling and solving on the fine mesh dominates the cost of the corresponding operations on the coarse mesh. This trend is clearly reflected in the data where for $k=1$, VMS is only marginally favourable and in some instances slightly unfavourable relative to the Galerkin method with $R$ being larger or at best close to 1. However, for all cases with $k>1$, the VMS formulation consistently achieves a favourable speed-up with $R$ being considerably smaller than 1. This test, and the scope of this work more broadly, is not intended to identify an optimal value of $k$, as this may vary depending on the application. However, the observed behaviour reinforces the cost analysis presented in \Cref{sec:VMS_cost} and demonstrates that the VMS approach can be favourable compared to its baseline Galerkin counterpart.
\begin{table}[htp]
    \centering
    \begin{tabular}{c|c|c|c|c}
         Compute time VMS (s) & $p = 1$ & $p = 2$ & $p = 3$ & $p = 4$ \\ \hline\hline
         $k = 1$ & \textcolor{red}{383.31} ($R = 1.51$) & \textcolor{codegreen}{578.63} ($R = 0.99$)& \textcolor{blue}{914.23} ($R = 1.01$) & \textcolor{magenta}{1608.22} ($R = 0.98$) \\ \hline
         $k = 2$ & \textcolor{codegreen}{453.44} ($R = 0.78$) & \textcolor{blue}{643.41} ($R = 0.71$) & \textcolor{magenta}{1156.18} ($R = 0.71$) & \textcolor{orange}{2335.31} ($R = 0.74$) \\ \hline
         $k = 3$ & \textcolor{blue}{588.26} ($R = 0.65$) & \textcolor{magenta}{895.80} ($R = 0.55$) & \textcolor{orange}{1786.11} ($R = 0.56$) & \textcolor{olive}{3322.17} ($R = 0.58$) \\ \hline
         $k = 4$ & \textcolor{magenta}{792.85} ($R = 0.49$) & \textcolor{orange}{1447.89} ($R = 0.46$) & \textcolor{olive}{2757.74} ($R = 0.48$) & \textcolor{violet}{4988.08} ($R = 0.48$)
         \end{tabular}
    \caption{Compute time for the proposed VMS methodology for the inviscid vortex roll-up test case for varying polynomial degrees with uniform $12 \times 12$ elements and $\Delta t = 0.01, N_{dt} = 800, i_{ave} \approx 10$}
    \label{tab:time_vms}
\end{table}

\begin{table}[htp]
    \centering
    \begin{tabular}{c|c|c|c|c|c|c|c}
         Compute time Galerkin (s) & $p = 2$ & $p = 3$ & $p = 4$ & $p = 5$ & $p = 6$ & $p = 7$ & $p = 8$ \\ \hline\hline
          & \textcolor{red}{253.01} & \textcolor{codegreen}{579.64} & \textcolor{blue}{909.11} & \textcolor{magenta}{1634.50} & \textcolor{orange}{3162.10} & \textcolor{olive}{5689.34} & \textcolor{violet}{10390.61} \\ 
         \end{tabular}
    \caption{Compute time for the base Galerkin methodology for the inviscid vortex roll-up test case for varying polynomial degrees with uniform $12 \times 12$ elements and $\Delta t = 0.01, N_{dt} = 800, i_{ave} \approx 10$}
    \label{tab:time_glk}
\end{table}

% t = 1e-3; VMS p = 2, k = 2: 4102.81 Galerkin p = 4: 5523.60

\begin{figure}[H]
    \begin{minipage}{0.49\linewidth}
        \centering
        \includegraphics[width=\linewidth]{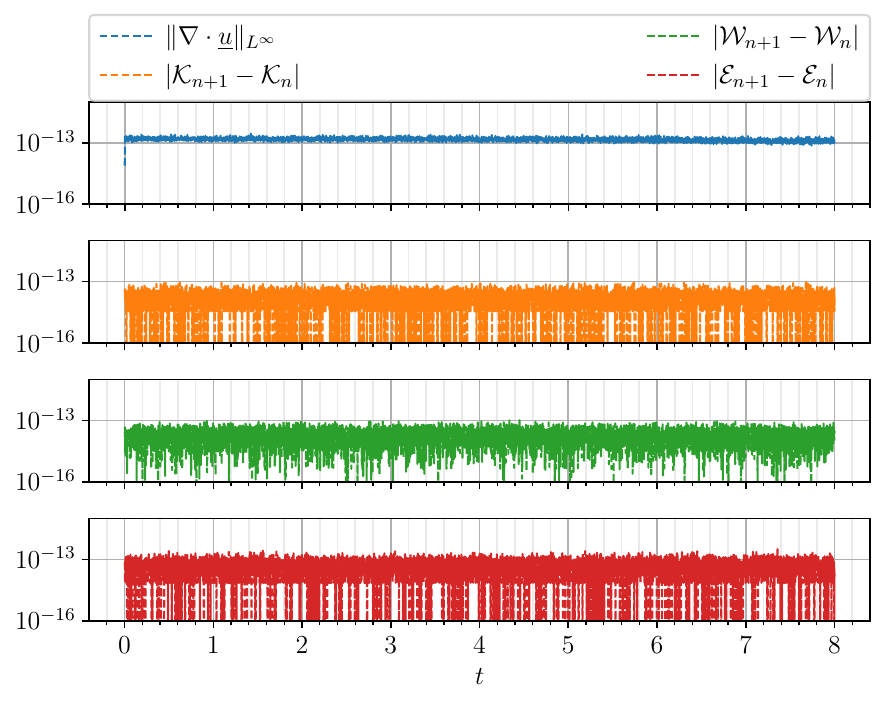}
        \caption{Conservation errors for the base Galerkin scheme}
        \label{fig:conservation_glk}
    \end{minipage}
    \begin{minipage}{0.49\linewidth}
        \centering
        \includegraphics[width=\linewidth]{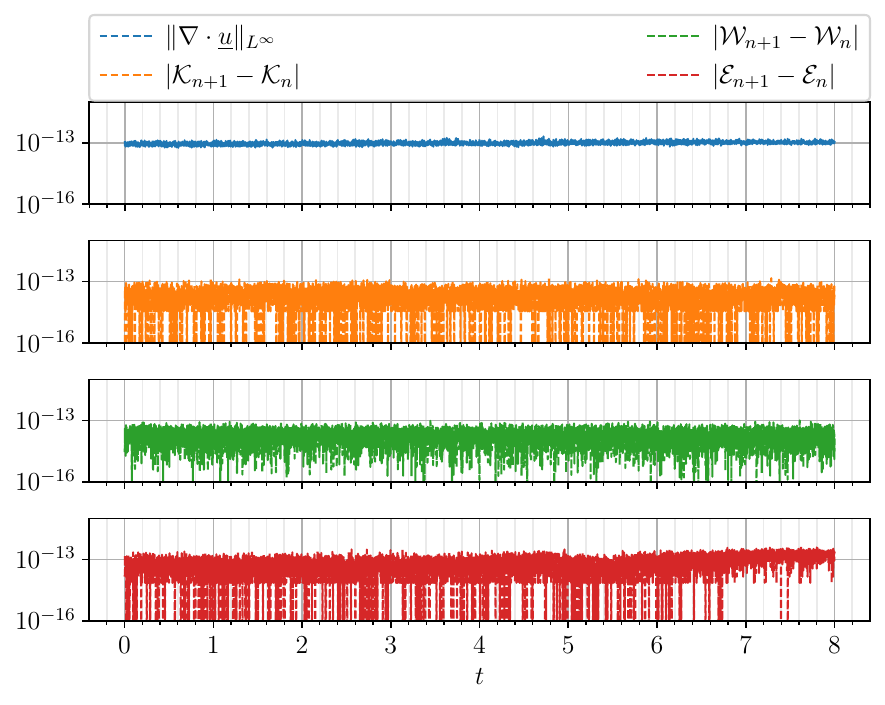}
        \caption{Conservation errors for the VMS approach}
        \label{fig:conservation_VMS}
    \end{minipage}
\end{figure}

Lastly, we close off the section by considering the conservation properties. Both the Galerkin and VMS schemes satisfy the same conservation properties as described by the MEEVC paper \cite{Zhang2024AConditions}. In the case of the VMS approach, conservation is naturally inherited from the orthogonality condition imposed between the resolved and unresolved scales by the projector, as shown in \Cref{sec:VMS_prop}. This ensures that the introduction of the fine-scale model does not disrupt the fundamental balance properties of the numerical scheme. \Cref{fig:conservation_glk} and \Cref{fig:conservation_VMS} show the mass, kinetic energy and total vorticity conservation errors for the base Galerkin scheme and the proposed VMS approach for the inviscid vortex roll-up test case, where we see that all the errors are in the order of machine precision. 

\section{Summary}
\label{sec:summary}
This work presented an algebraic VMS formulation tailored for high-order discretisations of the incompressible Navier-Stokes equations, designed to enhance coarse-scale accuracy while preserving key conservation properties. Building on a Galerkin foundation, the method introduces a scale separation via an orthogonal projector based on the underlying symmetric problem. Employing this projector in combination with the explicit construction of the Fine-Scale Greens' function of the symmetric operator enables a consistent treatment of unresolved scales without introducing artificial dissipation.

To evaluate the proposed method, we first performed convergence tests using the Taylor-Green vortex. These tests confirmed the expected convergence rates, including exponential convergence of the VMS solution to the exact projection, as well as exponential convergence of the unresolved scales to their exact counterparts. Next, we considered the inviscid vortex roll-up problem, which is a challenging benchmark characterised by sharp gradients and dynamic multiscale features. A highly-resolved solution obtained with the baseline Galerkin scheme served as a reference. By projecting this reference solution onto a coarse mesh, we established a benchmark representing the best-possible resolved-scale behaviour in the chosen norm. Comparisons against coarse-mesh simulations show that the proposed VMS approach closely reproduces this projection, outperforming the coarse Galerkin solution, which poorly approximated the exact projection of the reference. A cost analysis demonstrated that the proposed VMS methodology becomes computationally favourable relative to the baseline Galerkin method on the fine mesh when the cost of matrix assemblies and solves on the fine mesh dominates that on the coarse mesh. This was further confirmed through practical comparisons of computational costs, which highlighted the efficiency of the VMS approach.

Several avenues remain for future investigation. A more detailed study of the treatment of pressure projections within the VMS framework is needed. Additionally, the choice of $k$, which controls the fidelity of the Fine-Scale Greens' function approximation, warrants further exploration. While accurate reconstructions of the projected solution were achieved for relatively low values of $k$ (typically $\leq 4$) and computational favourability was observed for $k>1$ in the cases considered, a systematic determination of the optimal $k$ that balances computational cost and numerical accuracy remains application-dependent and is left for future work.

Crucially, both the Galerkin and VMS schemes were shown to satisfy the same set of conservation properties under the MEEVC framework, with conservation in the VMS case ensured by the orthogonality between resolved and unresolved scales. These findings demonstrate the capability of the proposed approach to produce accurate, stable, and conservative approximations even in under-resolved regimes, providing a strong foundation for future extensions to more complex problems.

\section*{Acknowledgements}
Suyash Shrestha and Esteban Ferrer acknowledge the funding received by Clean Aviation Joint Undertaking under the European Union’s Horizon Europe research and innovation programme under Grant Agreement HERA (Hybrid-Electric Regional Architecture) no. 101102007. Views and opinions expressed are, however, those of the author(s) only and do not necessarily reflect those of the European Union or CAJU. Neither the European Union nor the granting authority can be held responsible for them. Esteban Ferrer and Gonzalo Rubio acknowledge the funding received by the Grant DeepCFD (Project No. PID2022-137899OB-I00) funded by MICIU/AEI/10.13039/501100011033 and by ERDF, EU. We would like to thank Artur Palha for helping with the implementation in Firedrake.


\bibliographystyle{elsarticle-num} 
\bibliography{references, ref_old}

\newpage
\appendix

\section{Proof of well-posedness}
\label{sec:well_posed}
For the mixed formulation of the Stokes problem to be well-posed, we require that the bilinear form $(\varepsilon^h, \bar{\omega})_{L^2(\Omega)}$ is bounded and coercive in the null space of the curl and the bilinear forms $(\nabla \times \varepsilon^h,  \underline{\bar{u}})_{L^2(\Omega)}$ and $(\nabla \cdot \underline{v}^h, \bar{p})_{L^2(\Omega)}$ satisfy the \emph{inf-sup} conditions
    \[\adjustlimits\inf_{\underline{\bar{u}} \in \bar{V}} \sup_{\varepsilon^h \in \bar{U}} 
     \frac{(\nabla \times \varepsilon^h,  \underline{\bar{u}})_{L^2(\Omega)}}{\lVert \varepsilon^h \rVert_{H(\mathrm{curl})} \lVert \underline{\bar{u}}\rVert_{H(\mathrm{div})}} > \beta_{\omega} > 0, \quad \quad \quad \quad 
     \adjustlimits\inf_{{\bar{\pbar}} \in \bar{W}} \sup_{\underline{v}^h \in \bar{V}} 
     \frac{(\nabla \cdot \underline{v}^h, \bar{\pbar})_{L^2(\Omega)}}{\lVert \underline{v}^h \rVert_{H(\mathrm{div})} \lVert {\bar{\pbar}}\rVert_{L^2}} > \beta_{\underline{u}} > 0
     \]
We make use of the following assumptions/lemmas for the proof
\begin{itemize}
    \item The finite-dimensional Hilbert spaces $(\bar{U}$, $\bar{V}$, $\bar{W})$ are subsets of the infinite-dimensional spaces $(H(\mathrm{curl}, \Omega)$, $H(\mathrm{div}, \Omega)$, $L^2(\Omega))$ where each space may be decomposed into a null-space/range and its $L^2$-orthogonal complement
    \[\bar{U} = \mathrm{Ker}^h(\nabla \times) \oplus \mathrm{Ker}^h(\nabla \times)^{\perp}, \quad \bar{V} = \mathrm{Ker}^h(\nabla \cdot) \oplus \mathrm{Ker}^h(\nabla \cdot)^{\perp}, \quad \bar{W} = \mathrm{Range}^h(\nabla \cdot) \oplus \mathrm{Range}^h(\nabla \cdot)^{\perp}, \quad \]
    \item On a contractable domain, the curl operator is a bijective map between $\mathrm{Ker}^h(\nabla \times)^{\perp}$ and $\mathrm{Ker}^h(\nabla \cdot)$, and the divergence operator is a bijective map between $\mathrm{Ker}^h(\nabla \cdot)^{\perp}$ and $\mathrm{Range}^h(\nabla \cdot)$
    \[\nabla \times \: : \: \mathrm{Ker}^h(\nabla \times)^{\perp} \xrightarrow[]{} \mathrm{Ker}^h(\nabla \cdot), \quad \quad \nabla \cdot \: : \: \mathrm{Ker}^h(\nabla \cdot)^{\perp} \xrightarrow[]{} \mathrm{Range}^h(\nabla \cdot) \]
    \item Poincar\'{e}'s lemma states that $\exists C_{P_{curl}}, C_{P_{div}} > 0$ which are only functions of the topology of $\Omega$ such that
    \[\lVert \varepsilon^h \rVert_{H(\mathrm{curl})} \leq C_{P_{curl}} \lVert \nabla \times \varepsilon^h \rVert_{L^2}, \quad \forall \varepsilon^h \in \mathrm{Ker}^h(\nabla \times)^{\perp}, \quad \quad \lVert \underline{v}^h \rVert_{H(\mathrm{div})} \leq C_{P_{div}} \lVert \nabla \cdot \underline{v}^h \rVert_{L^2}, \quad \forall \underline{v}^h \in \mathrm{Ker}^h(\nabla \cdot)^{\perp}.\]
\end{itemize}

We start the proof by showing boundedness and coercivity. The bilinear form $\bilnear{\varepsilon^h}{\bar{\omega}}$ is trivially bounded in the null space $\mathrm{Ker}^h(\nabla \times)$
\[\bilnear{\varepsilon^h_0}{\bar{\omega}_0} \leq \lVert \varepsilon_0^h \rVert_{L^2} \lVert \bar{\omega}_0 \rVert_{L^2}, \quad \forall \varepsilon_0^h, \bar{\omega}_0 \in \mathrm{Ker}^h(\nabla \times).\]
Moreover the bilinear form is coercive in $\mathrm{Ker}^h(\nabla \times)$ with a coercivity constant of 1
\[\bilnear{\bar{\omega}_0}{\bar{\omega}_0} = \lVert \bar{\omega}_0\rVert_{H(\mathrm{curl})}^2, \quad \forall \bar{\omega}_0 \in \mathrm{Ker}^h(\nabla \times).\]

We now look at the first \emph{inf-sup} condition. We start by noting that $\underline{\bar{u}} \in \mathrm{Ker}^h(\nabla \cdot)$ and we first take the supremum over all $\varepsilon^h \in \bar{U}$ and thus we have the following inequality
\begin{align*}
     \sup_{\varepsilon^h \in \bar{U}} 
     \frac{\bilnear{\nabla \times \varepsilon^h}{\underline{\bar{u}}}}{\lVert \varepsilon^h \rVert_{H(\mathrm{curl})}} \geq \frac{\bilnear{\nabla \times \gamma^h}{\underline{\bar{u}}}}{\lVert \gamma^h \rVert_{H(\mathrm{curl})}}, \quad \forall \gamma^h \in \bar{U}, \forall \underline{\bar{u}} \in \mathrm{Ker}^h(\nabla \cdot).
\end{align*}
Based on the second point listed above we can note that there exists exactly one element in $\mathrm{Ker}^h(\nabla \times)^{\perp}$ which gets mapped to an element of $\mathrm{Ker}^h(\nabla \cdot)$ by the curl operator
\[ \exists \gamma^h \in \mathrm{Ker}^h(\nabla \times)^{\perp} \:  \mathrm{s.t} \: \nabla \times \gamma^h = \underline{\bar{u}}, \quad \forall \underline{\bar{u}} \in \mathrm{Ker}^h(\nabla \cdot).\]
We can thus fill in the particular $\gamma^h$ whose curl gives $\underline{\bar{u}}$ and get the following
\begin{align*}
     \sup_{\varepsilon^h \in \bar{U}} 
     \frac{\bilnear{\nabla \times \varepsilon^h}{\underline{\bar{u}}}}{\lVert \varepsilon^h \rVert_{H(\mathrm{curl})}} &\geq \frac{\bilnear{\underline{\bar{u}}}{\underline{\bar{u}}}}{\lVert \gamma^h \rVert_{H(\mathrm{curl})}} {\geq}  \frac{1}{C_{P_{curl}}}\frac{\lVert \underline{\bar{u}} \rVert_{L^2}^2}{\lVert \underline{\bar{u}} \rVert_{L^2}} \geq \frac{1}{C_{P_{curl}}} \lVert \underline{\bar{u}} \rVert_{H(\mathrm{div})}.
     % \frac{\bilnear{\nabla \times \gamma^h}{\underline{\bar{u}}}}{\lVert \gamma^h \rVert_{H(\mathrm{curl})}}, \quad \forall \gamma^h \in \bar{U} \backslash \{\varepsilon^h\}, \forall \underline{\bar{u}} \in \mathrm{Ker}^h(\nabla \cdot) \\
     % & = \frac{\bilnear{\underline{\bar{u}}}{\underline{\bar{u}}}}{\lVert \gamma^h \rVert_{H(\mathrm{curl})}} \stackrel{\S \Cref{lem:poincare}}{\geq}  \frac{1}{C_{P_{curl}}}\frac{\lVert \underline{\bar{u}} \rVert_{L^2}^2}{\lVert \underline{\bar{u}} \rVert_{L^2}} \geq \frac{1}{C_{P_{curl}}} \lVert \underline{\bar{u}} \rVert_{H(\mathrm{div})}.
\end{align*}
Then taking the infimum over $\underline{\bar{u}}$ yields the \emph{inf-sup} constant $\frac{1}{C_{P_{curl}}}$. Since the Poincar\'{e} constant is known to be positive and only a function of the topology of the domain, we can conclude that this \emph{inf-sup} condition is satisfied. 

We perform a similar analysis for the second \emph{inf-sup} condition. We take the supremum over all $\underline{v}^h \in \bar{V}$ and thus we have the following inequality
\begin{align*}
     \sup_{\underline{v}^h \in \bar{V}} 
     \frac{\bilnear{\nabla \cdot \underline{v}^h}{{\bar{\pbar}}}}{\lVert \underline{v}^h \rVert_{H(\mathrm{div})}} \geq \frac{\bilnear{\nabla \cdot \underline{\varphi}^h}{{\bar{\pbar}}}}{\lVert \underline{\varphi}^h \rVert_{H(\mathrm{div})}}, \quad \forall \underline{\varphi}^h \in \bar{V}, \forall {\bar{\pbar}} \in \bar{W}.
\end{align*}
Once again, we note that there exists exactly one element in $\mathrm{Ker}^h(\nabla \cdot)^{\perp}$ which gets mapped to an element of $\mathrm{Range}^h(\nabla \cdot)$ by the divergence operator
\[ \exists \underline{\varphi}^h \in \mathrm{Ker}^h(\nabla \cdot)^{\perp} \:  \mathrm{s.t} \: \nabla \cdot \underline{\varphi}^h = {\bar{\pbar}_0}, \quad \forall {\bar{\pbar}_0} \in \mathrm{Range}^h(\nabla \cdot).\]
% However, we define $\bar{\pbar} \in L^2(\Omega)$, so the fulfilment of this \emph{inf-sup} condition relies on $\mathrm{Range}(\nabla \cdot)^{\perp}$ which is ultimately determined by the boundary conditions on the normal velocity along $\Gamma_{n}$ and the boundary conditions on the pressure along $\Gamma_p$ with $\partial \Omega = \Gamma_n \cup \Gamma_p$. If the pressure is prescribed on part of the boundary i.e. $\Gamma_p \neq \varnothing$ and the normal velocity is prescribed on the rest of the boundary, then the functions in $\bar{V}$ are only non-zero along $\Gamma_p$ and the divergence operator maps to the full $L^2(\Omega)$ space, hence $\mathrm{Range}(\nabla \cdot)^{\perp} = \varnothing$. In that case, we have the following
However, we define $\bar{\pbar} \in \bar{W}$, so the fulfilment of this \emph{inf-sup} condition depends on the orthogonal complement of the range of the divergence operator, $\mathrm{Range}^h(\nabla \cdot)^{\perp}$. This, in turn, is determined by the boundary conditions, specifically, the normal component of the velocity along $\Gamma_n$ and the pressure along $\Gamma_p$, where $\partial \Omega = \Gamma_n \cup \Gamma_p$. If the pressure is prescribed on a portion of the boundary (i.e., $\Gamma_p \neq \varnothing$) and the normal velocity is prescribed on the remainder, then the functions in $\bar{V}$ are only non-zero on $\Gamma_p$. In this case, the divergence operator maps onto the full space $L^2(\Omega)$, implying that $\mathrm{Range}(\nabla \cdot)^{\perp} = \varnothing$. Under such conditions, we have the following
\begin{align*}
     \sup_{\underline{v}^h \in \bar{V}} 
     \frac{\bilnear{\nabla \cdot \underline{v}^h}{{\bar{\pbar}}}}{\lVert \underline{v}^h \rVert_{H(\mathrm{div})}} & \geq \frac{\bilnear{{\bar{\pbar}}}{{\bar{\pbar}}}}{\lVert \underline{\varphi}^h \rVert_{H(\mathrm{div})}} {\geq}  \frac{1}{C_{P_{div}}}\frac{\lVert {\bar{\pbar}} \rVert_{L^2}^2}{\lVert {\bar{\pbar}} \rVert_{L^2}} \geq \frac{1}{C_{P_{div}}} \lVert {\bar{\pbar}} \rVert_{L^2},
     % \frac{\bilnear{\nabla \cdot \underline{\varphi}^h}{{\bar{\pbar}}}}{\lVert \underline{\varphi}^h \rVert_{H(\mathrm{div})}}, \quad \forall \underline{\varphi}^h \in \bar{V}, \forall {\bar{\pbar}} \in L^2(\Omega) \\
     % & \geq \frac{\bilnear{{\bar{\pbar}}}{{\bar{\pbar}}}}{\lVert \underline{\varphi}^h \rVert_{H(\mathrm{div})}} \stackrel{\S \Cref{lem:poincare}}{\geq}  \frac{1}{C_{P_{div}}}\frac{\lVert {\bar{\pbar}} \rVert_{L^2}^2}{\lVert {\bar{\pbar}} \rVert_{L^2}} \geq \frac{1}{C_{P_{div}}} \lVert {\bar{\pbar}} \rVert_{L^2},
\end{align*}
where taking the infimum over $\bar{\pbar}$ yields the \emph{inf-sup} constant $\frac{1}{C_{P_{div}}}$, thereby satisfying the \emph{inf-sup} condition. If, on the other hand, the pressure is not prescribed on any part of the boundary (i.e., $\Gamma_p = \varnothing$), and only the normal velocity is prescribed along $\partial \Omega$, then the functions in $\bar{V}$ vanish on the entire boundary. In this case, the divergence operator maps onto a subspace of $L^2(\Omega)$ with $\mathrm{Range}(\nabla \cdot)^{\perp} = \mathbb{R}$, corresponding to the space of constant functions. As a result, the \emph{inf-sup} condition is only satisfied with the pressure in $L^2(\Omega) \backslash \mathbb{R}$, with the previously stated \emph{inf-sup} constant $\frac{1}{C_{P_{{div}}}}$. The singular constant pressure mode is a well-known physical singular mode for incompressible flows and is not a numerical artifact.

% If, on the other hand, the pressure is not prescribed anywhere on the boundary i.e. $\Gamma_p = \varnothing$, and only the normal velocity is prescribed on $\partial \Omega$, then the functions in $\bar{V}$ are all zero along $\partial \Omega$ and we have $\mathrm{Range}(\nabla \cdot)^{\perp} = \mathbb{R}$ which is the space of constant functions. For this case, the \emph{inf-sup} condition is strictly violated with the singular mode being the constant pressure. We can, of course, remedy this by considering pressure in $L^2(\Omega) \backslash \{\mathbb{R}\}$ in which case the aforementioned \emph{inf-sup} constant of $\frac{1}{C_{P_{div}}}$ holds.

\section{Norm optimality}
\label{app:optimal}
\begin{proof}
    We have the following orthogonality conditions for the projector
\begin{equation}
    \begin{array}{ccccccl}
        \bilnear{\varepsilon^h}{\bar{\omega} - {\omega}} & - & \bilnear{\nabla \times \varepsilon^h}{\bar{\underline{u}} - {\underline{u}}} & & = &0,  & \forall \varepsilon^h \in \bar{\mathcal{C}} \\[1.5ex]
        -\bilnear{\underline{v}^h}{\nabla \times (\bar{\omega} - {\omega})}
        & & + & \bilnear{\nabla \cdot \underline{v}^h}{\bar{\pbar} - {\pbar}}
        &= & \underline{0}, 
        & \forall \underline{v}^h \in \bar{\mathcal{D}} \\[1.5ex]
        & & \bilnear{\eta^h}{\nabla \cdot (\bar{\underline{u}} - \underline{{u}})} & & = & 0, 
        &  \forall \eta^h \in \bar{\mathcal{S}}
    \end{array}.
    \label{eq:stokes_projector_app}
\end{equation}
% We want to show that the solution of \eqref{eq:stokes_projector_app} achieves optimality in the $H_0(\mathrm{curl}, \Omega)$ norm 
% \textcolor{red}{How restrictive is it that you set $\omega$ to zero on the boundary? Note that you need to mention that this is the 2D case, because in 3D there are two possible traces you can set to zero, making the expression $H_0(\mathrm{curl}, \Omega)$ ambiguous. Looking at the proof, it suffices to consider $H(\mathrm{curl}, \Omega)$} . 
We start with the error estimate in the null space, where we have
\begin{align*}
    \lVert  {\bar{\omega}}_0 - {{\omega}} \rVert_{L^2}^2 &= \bilnear{ {\bar{\omega}}_0 - {{\omega}}}{{\bar{\omega}}_0 - {{\omega}}} \\
    &= ({\underbrace{{\bar{\omega}}_0 - {{\varepsilon}}^h_0}_{\in \mathrm{Ker}^h(\nabla \times)}}, {{\bar{\omega}_0 - {{\omega}}}})_{L^2(\Omega)} + \bilnear{{{\varepsilon}}^h_0 - {{\omega}}}{{\bar{\omega}}_0 - {{\omega}}}, \quad \forall {{\varepsilon}}^h_0 \in \mathrm{Ker}^h(\nabla \times).
\end{align*}
Making a particular choice of $\varepsilon^h = {\bar{\omega}}_0 - {{\varepsilon}}^h_0$ as the test function in the first equation of \eqref{eq:stokes_projector_app}, we get $({\bar{\omega}}_0 - {{\varepsilon}}^h_0, {{\bar{\omega}_0 - {{\omega}}}})_{L^2(\Omega)} = 0$. Which leads to
\begin{align*}
    \lVert  {\bar{\omega}}_0 - {{\omega}} \rVert_{L^2}^2 &= \bilnear{{{\varepsilon}}^h_0 - {{\omega}}}{{\bar{\omega}}_0 - {{\omega}}}, \quad \forall {{\varepsilon}}^h_0 \in \mathrm{Ker}^h(\nabla \times) \\
    &\leq \lVert {{\varepsilon}}^h_0 - {{\omega}} \rVert_{L^2}  \lVert {\bar{\omega}}_0 - {{\omega}} \rVert_{L^2} , \quad \forall {{\varepsilon}}^h_0 \in \mathrm{Ker}^h(\nabla \times) \\
    \lVert  {\bar{\omega}}_0 - {{\omega}} \rVert_{L^2} &\leq  \lVert {\bar{\omega}}_0 - {{\omega}} \rVert_{L^2} , \quad \forall {{\varepsilon}}^h_0 \in \mathrm{Ker}^h(\nabla \times),
\end{align*}
and hence we have
\begin{equation}
    \lVert  {\bar{\omega}}_0 - {{\omega}} \rVert_{L^2} = \inf_{\varepsilon^h_0 \in \mathrm{Ker}^h(\nabla \times)} \lVert \varepsilon^h_0 - \omega  \rVert_{L^2}.
\end{equation}

Similarly, the error estimate for the components orthogonal to the closed fields follows as
\begin{align*}
    \lVert \nabla \times {\bar{\omega}}_{\perp} - \nabla \times {{\omega}} \rVert_{L^2}^2 &= \bilnear{\nabla \times {\bar{\omega}}_{\perp} - \nabla \times {{\omega}}}{\nabla \times {\bar{\omega}}_{\perp} - \nabla \times {{\omega}}} \\
    \begin{split}
        &= ({\underbrace{\nabla \times {\bar{\omega}}_{\perp} - \nabla \times {{\varepsilon}}^h_{\perp}}_{\in \mathrm{Ker}^h(\nabla \cdot)}}, {\nabla \times {\bar{\omega}_{\perp} - \nabla \times {{\omega}}}})_{L^2(\Omega)} \\ &+ \bilnear{\nabla \times {{\varepsilon}}^h_{\perp} - \nabla \times {{\omega}}}{\nabla \times {\bar{\omega}}_{\perp} - \nabla \times {{\omega}}}
    \end{split} \:, \quad \forall {{\varepsilon}}^h_{\perp} \in \mathrm{Ker}^h(\nabla \times)^{{\perp}}.
\end{align*}
We can make a particular choice of $\underline{v}^h = \nabla \times {\bar{\omega}}_{\perp} - \nabla \times {{\varepsilon}}^h_{\perp}$ in the second equation of \eqref{eq:stokes_projector_app} and find that $({\nabla \times {\bar{\omega}}_{\perp} - \nabla \times {{\varepsilon}}^h_{\perp}}, {\nabla \times {\bar{\omega} - \nabla \times {{\omega}}}})_{L^2(\Omega)} = 0$. We thus have
\begin{align}
    \lVert \nabla \times {\bar{\omega}}_{\perp} - \nabla \times {{\omega}} \rVert_{L^2}^2 &= \bilnear{\nabla \times {{\varepsilon}}^h_{\perp} - \nabla \times {{\omega}}}{\nabla \times {\bar{\omega}}_{\perp} - \nabla \times {{\omega}}}, \quad \forall {{\varepsilon}}^h_{\perp} \in \mathrm{Ker}^h(\nabla \times)^{{\perp}} \nonumber \\
    &\leq \lVert \nabla \times {{\varepsilon}}^h_{\perp} - \nabla \times {{\omega}} \rVert_{L^2} \lVert \nabla \times {\bar{\omega}}_{\perp} - \nabla \times {{\omega}} \rVert_{L^2}, \quad \forall {{\varepsilon}}^h_{\perp} \in \mathrm{Ker}^h(\nabla \times)^{{\perp}} \nonumber \\
    % &\leq \lVert {\varepsilon}^h - {{\omega}} \rVert_{\H_0(\mathrm{curl}, \Omega)}^2\lVert {\bar{\omega}} - {{\omega}} \rVert_{\H_0(\mathrm{curl}, \Omega)} \\
   \lVert \nabla \times {\bar{\omega}}_{\perp} - \nabla \times {{\omega}} \rVert_{L^2} & \leq \lVert \nabla \times {\varepsilon}^h_{\perp} - \nabla \times {{\omega}} \rVert_{L^2}, \quad \forall {{\varepsilon}}^h_{\perp} \in \mathrm{Ker}^h(\nabla \times)^{{\perp}} \nonumber 
\end{align}
\begin{equation}
    \lVert \nabla \times {\bar{\omega}}_{\perp} - \nabla \times {{\omega}} \rVert_{L^2} = \inf_{{\varepsilon}^h_{\perp} \in \mathrm{Ker}^h(\nabla \times)^{{\perp}}} \lVert \nabla \times {\varepsilon}^h_{\perp} - \nabla \times {{\omega}} \rVert_{L^2}.
\end{equation}
\end{proof}
\end{document}